\documentclass[letterpaper,onecolumn]{IEEEtran}
\usepackage{amsmath,amsfonts}
\usepackage{amsthm}
\usepackage{amssymb}
\usepackage{algorithmic}
\usepackage{algorithm}
\usepackage{array}
\usepackage[caption=false,font=normalsize,labelfont=sf,textfont=sf]{subfig}
\usepackage{textcomp}
\usepackage{stfloats}
\usepackage{url}
\usepackage{verbatim}
\usepackage{graphicx}
\usepackage{cite}
\usepackage{mathrsfs}
\usepackage{mathtools}
\usepackage{color}
\usepackage{enumitem}
\usepackage{etoc}
\usepackage{hyperref}

\theoremstyle{plain}
\newtheorem{theorem}{Theorem}

\newtheorem{corollary}{Corollary}
\newtheorem{proposition}{Proposition}
\newtheorem{definition}{Definition}

\allowdisplaybreaks[1]

\begin{document}

\title{Improved Nearly Minimax Prediction for Independent Poisson Processes Under Kullback--Leibler Loss}

\author{
\IEEEauthorblockN{Xiao Li\IEEEauthorrefmark{1}\IEEEauthorrefmark{3} and Fumiyasu Komaki\IEEEauthorrefmark{1}\IEEEauthorrefmark{2}}
\IEEEauthorblockA{\IEEEauthorrefmark{1}Department of Mathematical Informatics, 
Graduate School of Information Science and Technology, The University of Tokyo, 7-3-1 Hongo, Bunkyo-ku, Tokyo, 113-0033, Japan}
\IEEEauthorblockA{\IEEEauthorrefmark{2}RIKEN Center for Brain Science, 
2-1 Hirosawa, Wako City, Saitama, 351-0198, Japan}
\IEEEauthorblockA{\IEEEauthorrefmark{3}Corresponding author: lixiaoms@163.com}
\thanks{}
}

\markboth{IEEE Transactions on Information Theory,~Vol.~XX, No.~X, Month~2026}%
{Li and Komaki: Improved Nearly Minimax Prediction for Independent Poisson Processes}


\maketitle

\noindent\textit{This work has been submitted to the IEEE for possible publication. 
Copyright may be transferred without notice, after which this version may no longer be accessible.}

\begin{abstract}
Simultaneous predictive distributions for independent Poisson observables are investigated, and the performance of predictive distributions is evaluated using the Kullback--Leibler (K-L) loss. This study introduces intuitive sufficient conditions, based on superharmonicity of priors, to improve the Bayesian predictive distribution based on the Jeffreys prior. The sufficient conditions exhibit a certain analogy with those known for the multivariate normal distribution. Additionally, this study examines the case where the observed data and target variables to be predicted are independent Poisson processes with different durations. Examples that satisfy the sufficient conditions are provided, including point and subspace shrinkage priors. The K-L risk of the improved predictions is demonstrated to be less than $1.04$ times a minimax lower bound.
\end{abstract}


\section{Introduction}
\IEEEPARstart{T}{he} problem of predicting Poisson random variables $y = (y_1,y_2,...,y_d)$ using observations $x = (x_1,x_2,...,x_d)$ arises in various statistical applications. We consider the fundamental problem of predictive distribution estimation for independent Poisson models. It is assumed that $x = (x_1,\dots,x_d)$ and $y = (y_1,\dots,y_d)$ are independently distributed according to
\begin{equation*}
p(x\mid\lambda)=\prod_{i=1}^d\frac{(r\lambda_i)^{x_i}e^{-r\lambda_i}}{x_i!},
\quad
p(y\mid\lambda)=\prod_{i=1}^d\frac{(s\lambda_i)^{y_i}e^{-s\lambda_i}}{y_i!},
\end{equation*}
where $r,s>0$ are known and $\lambda=(\lambda_1,\dots,\lambda_d)$ is the unknown parameter vector.

The Bayesian prediction distribution based on prior $\pi(\lambda)$ is expressed as:
\begin{equation*}
p_{\pi}(y\mid x)=\frac{\int p(x,y\mid \lambda)\pi(\lambda)\mathrm{d} \lambda}{\int p(x\mid \lambda)\pi(\lambda)\mathrm{d} \lambda}=\frac{\int p(x\mid \lambda)p(y\mid \lambda)\pi(\lambda)\mathrm{d} \lambda}{\int p(x\mid \lambda)\pi(\lambda)\mathrm{d} \lambda}.
\end{equation*}
The Kullback--Leibler (K-L) loss of the predictive distribution $p_{\pi}(y\mid x)$ is used in this study, which is expressed as:
\begin{equation*}
D(p(y \mid \lambda),p_{\pi}(y\mid x))
=\sum_yp(y\mid\lambda)
\log \frac{p(y\mid\lambda)}{p_{\pi}(y\mid x)}.
\end{equation*}

A natural prior is the Jeffreys prior: 
\begin{equation*}
\pi_{\mathrm{J}}(\lambda)\mathrm{d}\lambda_1\mathrm{d}\lambda_2\cdots \mathrm{d}\lambda_d\propto \lambda_1^{-1/2}\lambda_2^{-1/2}\cdots\lambda_d^{-1/2}\mathrm{d}\lambda_1\mathrm{d}\lambda_2\cdots \mathrm{d}\lambda_d,
\end{equation*}
which is a frequently used noninformative prior. The prior $\pi_{\mathrm{J}}(\lambda)$ has an advantage in that the K-L risk of the Bayesian predictive distribution $p_{\mathrm{J}}(y\mid x)$ has a small upper bound for any $\lambda$. The corresponding theorem is provided in this study. Therefore, a main focus of this study is the construction of a prior that is superior to the prior $\pi_{\mathrm{J}}(\lambda)$. This study provides sufficient conditions for the prior $\pi(\lambda)$ to make $p_{\pi}(y\mid x)$ dominate the Bayesian predictive distribution $p_{\mathrm{J}}(y\mid x)$ based on $\pi_{\mathrm{J}}(\lambda)$. Additionally, the role of superharmonic functions in improving the prediction for a multivariate Poisson vector is revealed as a counterpart to the case of a multivariate normal vector (George et al.~\cite{george2006improved}).

Numerous studies have been performed on the simultaneous estimation of Poisson parameters. Clevenson and Zidek~\cite{clevenson1975simultaneous} proposed generalized Bayes estimators dominating the maximum likelihood estimator when $d\ge2$ under the standardized squared error loss $\sum\lambda_i^{-1}(\hat\lambda_i-\lambda_i)^2.$ Tsui and Press~\cite{tsui1982simultaneous} studied the estimation under the generalized loss function $\sum(\hat\lambda_i-\lambda_i)^2/\lambda_i^{k},$ where $k$ is a given positive integer.
Ghosh and Yang~\cite{ghosh1988simultaneous} characterized admissible linear estimators of multiple Poisson parameters under K-L loss.

Estimation of parameters under K-L loss can be generalized to a predictive distribution problem, which is important for several statistical scenarios. The predictive method was shown to be preferable in Aitchison~\cite{aitchison1975goodness}. Noninformative priors or vague prior distributions are often used for constructing Bayesian predictive distributions. The Jeffreys prior has been widely used in various problems, such as in Akaike~\cite{akaike1978new} and Clarke and Barron~\cite{clarke1994jeffreys}.

Compared with the large number of estimation studies, decision theory regarding predictive distributions on the Poisson model has been developed relatively recently. A class of prior distributions, 
\begin{equation*}
\pi_{\alpha,\mathrm{J}}(\lambda)\mathrm{d}\lambda_1\mathrm{d}\lambda_2\cdots \mathrm{d}\lambda_d\propto \frac{\lambda_1^{-1/2}\lambda_2^{-1/2}\cdots\lambda_d^{-1/2}}
{(\lambda_1+\lambda_2+\cdots+\lambda_d)^\alpha}\mathrm{d}\lambda_1\mathrm{d}\lambda_2\cdots \mathrm{d}\lambda_d
\end{equation*}
was proposed in Komaki~\cite{komaki2004simultaneous}. Let $\theta_i\coloneqq \sqrt{\lambda_i},\ i=1,\dots,d.$ The Bayesian predictive distribution based on the prior $\Vert\theta\Vert^{2-d}\mathrm{d}\theta$ was shown to dominate that based on the Jeffreys prior when $d\ge3$ (Komaki~\cite{komaki2004simultaneous}). Komaki~\cite{komaki2015simultaneous} considered the problem of independent Poisson processes with different durations and introduced a class of prior densities that is a generalization of $\pi_{\alpha,\mathrm{J}}(\lambda)$. The corresponding Bayesian predictive distribution was shown to dominate that based on the Jeffreys prior. A class of proper priors was proposed, and Bayesian predictive distributions and estimators based on the priors were shown to dominate the Bayesian predictive distribution and estimator based on the Jeffreys prior under K-L loss (Komaki~\cite{komaki2006class}). The proper priors with respect to $\theta$ coincides with the function of Strawderman's prior in a normal model (Strawderman~\cite{strawderman1971proper}). Recently, Hamura and Kubokawa~\cite{hamura2020bayesian} studied the Bayesian predictive distribution for a Poisson model with parametric restriction under K-L loss. Yano et al.~\cite{yano2021minimax} presented a class of Bayesian predictive distributions that attain asymptotic minimaxity in sparse Poisson sequence models.

It is natural that similar results hold simultaneously for the multivariate normal and Poisson models from the viewpoint of a model manifold with the Fisher metric (Komaki~\cite{komaki2006class}). There are several counterparts for these two models. The Bayesian predictive distribution for a multivariate normal model $\mathrm{N}_d(\mu,\sigma^2I)$ based on Stein's harmonic prior (Stein~\cite{stein1974estimation}), i.e., 
\begin{equation*}
\pi(\mu)=\Vert\mu\Vert^{2-d},
\end{equation*}
dominates that based on the Jeffreys prior (Komaki~\cite{komaki2001shrinkage}). This result is similar to that reported in Komaki~\cite{komaki2004simultaneous}. Johnstone~\cite{johnstone1984admissibility} studied the admissibility and recurrence in estimating a Poisson mean under the standardized squared error loss, which is a counterpart to the diffusion characterization of admissibility in the estimation of a multivariate normal mean that was introduced by Brown~\cite{brown1971admissible}.

George et al.~\cite{george2006improved} generalized the result presented in Komaki~\cite{komaki2001shrinkage} for the multivariate normal model and proved that Bayesian predictive distributions based on superharmonic priors dominate those based on the Jeffreys prior. Thus, it is natural to speculate that a similar result exists for the independent Poisson observable model. This speculation is confirmed in the present study, which demonstrates the relationship between the superharmonic function and improved Bayesian prediction in the Poisson model. Although superharmonic priors have been frequently studied in the context of the normal model (Okudo and Komaki~\cite{okudo2024predictive}; Marchand and Strawderman~\cite{marchand2025on}), to our knowledge, this is the first non-asymptotic result concerning superharmonic functions for priors in count models.

In this paper, Theorem~\ref{Theorem 1} provides sufficient conditions for prior $\pi(\lambda)$ to make $p_{\pi}(y\mid x)$ dominate the Bayesian predictive distribution $p_{\mathrm{J}}(y\mid x)$ based on $\pi_{\mathrm{J}}(\lambda)$. Sufficient conditions are also provided for prior $\pi(\lambda)$ to make $p_{\pi}(y\mid x)$ dominate $p_{\mathrm{J}}(y\mid x)$ for the simultaneous prediction of independent Poisson processes with different durations (Theorem~\ref{Theorem 2}). Let $f$ denote $\pi(\lambda)/\pi_{\mathrm{J}}(\lambda)$; it is shown that the key is to study the function $f(\sqrt{\lambda_1},\cdots,\sqrt{\lambda_d})$. The proof of the theorems may be applicable to the study of improved Bayesian predictive distribution for other distributions. 

In Corollaries~\ref{Corollary 1} and \ref{Corollary 3}, the Bayesian predictive distribution $p_{\mathrm{J}}(y\mid x)$ based on the Jeffreys prior is improved using superharmonic functions. The result is a counterpart to George et al.~\cite{george2006improved}, indicating the similarity and difference between the predictive distribution theory of the multivariate Poisson and normal models. Using the results, classes of priors are introduced to improve the Bayesian predictive distribution $p_{\mathrm{J}}(y\mid x)$. In Theorem~\ref{Theorem 3}, the K-L risk of $p_{\mathrm{J}}(y\mid x)$ is shown to be less than $1.04$ times a minimax lower bound. Therefore, these improved Bayesian predictive distributions are nearly minimax.

This paper is organized as follows. Theorem~\ref{Theorem 1}, i.e., the result of sufficient conditions for prior $\pi(\lambda)$ to improve the Bayesian predictive distribution $p_{\mathrm{J}}(y\mid x)$ is presented in Section~\ref{sec:superharmonic}. Then Corollary~\ref{Corollary 1}, i.e., the main result on using superharmonic functions to improve the Bayesian predictive distribution $p_{\mathrm{J}}(y\mid x)$, is provided. Similar results for independent Poisson processes with different durations are presented in Section~\ref{sec:different}. Examples are provided in Section~\ref{sec:examples}, including point and subspace shrinkage priors. Here, the Bayesian predictive distributions based on the example priors are proven to dominate those based on the Jeffreys prior. Numerical experiments and a real data application are presented in Section~\ref{sec:numerical}. Finally, the results and example priors are discussed in Section~\ref{sec:discussion}. The detailed proofs of the theorems and propositions are provided in the Supplemental Material.

\section{Improved Bayesian Prediction Using Superharmonic Function}
\label{sec:superharmonic}

In this section, we define $\theta_i\coloneqq \sqrt{\lambda_i},\ i=1,\dots,d$. The goal is to find function $f(\theta_1,\dots,\theta_d)$ such that for prior 
\begin{equation*}
\pi_{f,\mathrm{J}}(\lambda)\mathrm{d}\lambda_1\mathrm{d}\lambda_2\cdots \mathrm{d}\lambda_d\propto f(\theta_1,\dots,\theta_d)\lambda_1^{-1/2}\lambda_2^{-1/2}\cdots\lambda_d^{-1/2}\mathrm{d}\lambda_1\mathrm{d}\lambda_2\cdots \mathrm{d}\lambda_d,
\end{equation*}
the Bayesian predictive distribution $p_{f,\mathrm{J}}(y\mid x)$ based on $\pi_{f,\mathrm{J}}(\lambda)$ dominates the Bayesian predictive distribution $p_{\mathrm{J}}(y\mid x)$ based on $\pi_{\mathrm{J}}(\lambda)$ under K-L loss. Prior $\pi_{f,\mathrm{J}}(\lambda)\mathrm{d}\lambda$ is equivalent to prior $f(\theta)\mathrm{d}\theta$.

In Theorem~\ref{Theorem 1}, we show that if $f$ satisfies certain conditions, the Bayesian predictive distribution ${p_{f,\mathrm{J}}(y\mid x)}$ dominates the Bayesian predictive distribution $p_{\mathrm{J}}(y\mid x)$. These conditions include three regularity and two essential conditions. As a corollary, for any superharmonic function satisfying certain regularity conditions, we can construct $f$ such that $p_{f,\mathrm{J}}(y\mid x)$ based on prior $f(\theta)\mathrm{d}\theta$ dominates $p_{\mathrm{J}}(y\mid x)$ based on the Jeffreys prior. The specific description of the result is presented in Corollary~\ref{Corollary 1}.

In brief, $p_{f,\mathrm{J}}(y\mid x)$ dominates $p_{\mathrm{J}}(y\mid x)$ if $f(\theta)=\sum_{a\in\{1,-1\}^d}h(a_1\theta_1,a_2\theta_2,\dots,a_d\theta_d)$ and $h(\theta)$ is a superharmonic function. Here,
\begin{align*}
    &\sum_{a\in\{1,-1\}^d}h(a_1\theta_1,a_2\theta_2,\dots,a_d\theta_d)\\
    &\coloneqq h(\theta_1,\theta_2,\dots,\theta_d)+h(-\theta_1,\theta_2,\dots,\theta_d)+\cdots+h(\theta_1,-\theta_2,\dots,-\theta_d)+h(-\theta_1,-\theta_2,\dots,-\theta_d).
\end{align*}
However, $p_{f,\mathrm{J}}(y\mid x)$ based on the superharmonic prior $f(\theta)=h(\theta)$ does not necessarily dominate $p_{\mathrm{J}}(y\mid x)$, which is different from the result in the normal model (George et al.~\cite{george2006improved}). For example, the K-L risk of $p_{\mathrm{J}}(y\mid x)$ is $0.56$ when $r=s=1$, $d=3$, and $\lambda=(0.4,0.4,0.4)$. The K-L risk of $p_{f,\mathrm{J}}(y\mid x)$ with superharmonic prior
\begin{align*}
    f(\theta)=\Big((\theta_1-2)^2+(\theta_2-2)^2+(\theta_3-2)^2\Big)^{-1/2}
\end{align*}
is $0.62$, which is larger than that of $p_{\mathrm{J}}(y\mid x)$. The K-L risk of $p_{f,\mathrm{J}}(y\mid x)$ is $0.56$ if
\begin{align*}
    f(\theta)=\sum_{a\in\{1,-1\}^d}\Big((a_1\theta_1-2)^2+(a_2\theta_2-2)^2+(a_3\theta_3-2)^2\Big)^{-1/2}.
\end{align*}

\begin{figure}[!t]
\centering
\includegraphics[width=0.6\columnwidth]{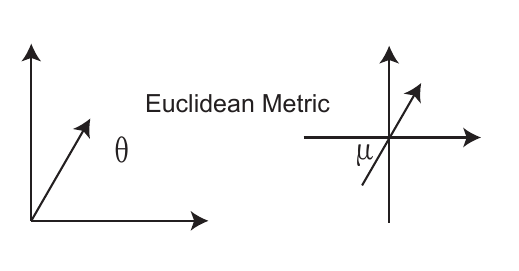}
\caption{Poisson model with parameter $\theta$ and normal model with parameter $\mu$ exhibit the same Fisher metric even though their domains of definition are different.}
\label{euclidean}
\end{figure}

Note that the domain of definition for the mean parameter $\mu$ of the normal model is $\mathbb{R}^d$, whereas the domain of definition for parameter $\theta$ of the Poisson model is $\mathbb{R}_+^d$. In addition, the Fisher metric of the normal model with parameter $\mu$ and that of the Poisson model with parameter $\theta$ are Euclidean metric (Fig.~\ref{euclidean}). The similarity and difference between the Poisson and normal models are consistent with the results of this study. The Bayesian predictive distribution based on the Jeffreys prior could be improved using superharmonic functions in both models because the Fisher metric is Euclidean. ${f(\theta)=\sum_{a\in\{1,-1\}^d}h(a_1\theta_1,a_2\theta_2,\dots,a_d\theta_d)}$ is set to ensure that the prior $f(\theta)\mathrm{d}\theta$ is superior to the Jeffreys prior. This is consistent with the domain of definition for the parameter $\theta$, which is $\mathbb{R}_+^d$.

We first introduce the regularity conditions required for our theoretical development:
\begin{subequations}
\begin{align}
&\int f(\theta) \prod_{j=1}^{d}\theta_j^{2z_j}\exp(-r\theta_j^2)\mathrm{d}\theta<\infty, \label{eq:t1-condition-a}\\
&\int \Bigl\lvert \frac{\partial f}{\partial\theta_i}(\theta) \Bigr\rvert \prod_{j=1}^{d}\theta_j^{2z_j}\exp(-r\theta_j^2)\mathrm{d}\theta<\infty, \label{eq:t1-condition-b}\\
&\int \Bigl\lvert \frac{\partial^2 f}{\partial\theta_i^2}(\theta) \Bigr\rvert \prod_{j=1}^{d}\theta_j^{2z_j}\exp(-r\theta_j^2)\mathrm{d}\theta<\infty, \label{eq:t1-condition-c}
\end{align}
for every $z\in\mathbb{N}^d$, $r>0$, and $i$.
\end{subequations}
Conditions \eqref{eq:t1-condition-a}, \eqref{eq:t1-condition-b}, and \eqref{eq:t1-condition-c} guarantee that all terms remain integrable when integration by parts is used in the proof. These regularity conditions hold if $f(\theta)$ and its first two partial derivatives are all bounded by $\exp(g(\theta))$ for some function $g(\theta)$ that grows slower than $\sum_{j=1}^d \theta_j^2$ in the sense that
\begin{equation*}
\frac{g(\theta)}{\sum_{j=1}^d \theta_j^2} \to 0 \quad \text{as} \quad \sum_{j=1}^d \theta_j^2 \to \infty.
\end{equation*}
This ensures that the product $\prod_{j=1}^{d}\theta_j^{2z_j}\exp\bigl(g(\theta)-r\sum_{j=1}^d\theta_j^2\bigr)$ is integrable.

We say that a predictive distribution $p_1$ weakly dominates $p_2$ if its risk is nowhere larger than that of $p_2$. For $z=(z_1,z_2,\dots,z_d)\in \mathbb{N}^d$, define
\begin{equation*}
F(z,r)\coloneqq F(z_1,z_2,\dots,z_d,r)= \int f(\theta_1,\dots,\theta_d)\prod_{i=1}^d\frac{r^{z_i+1/2}\lambda_i^{z_i-1/2}\exp(-r\lambda_i)}{\Gamma(z_i+1/2)}\mathrm{d}\lambda.
\end{equation*}

\begin{theorem}\label{Theorem 1}\leavevmode
\begin{enumerate}[label=\arabic*)]
\item
Let $\theta_i\coloneqq \sqrt{\lambda_i},\ i=1,\dots,d$. The Bayesian predictive distribution $p_{f,\mathrm{J}}(y\mid x)$ weakly dominates the Bayesian predictive distribution $p_{\mathrm{J}}(y\mid x)$ if for every $z\in \mathbb{N}^d$, $r>0$,
\begin{equation}
    \sum_{i=1}^dz_i(F(z,r)-F(z-\delta_i,r))+\sum_{i=1}^d(z_i+1/2)(F(z,r)-F(z+\delta_i,r))\ge 0, \label{eq:Fineq}
\end{equation}
where $\delta_i$ denotes the $d$-dimensional vector whose $i^{\text{th}}$ element is 1, all other elements are 0, and $F(z-\delta_i,r)$ is defined as $1$ if $z_i=0$. 
Furthermore, if $F(z,r)$ is not a constant function of $z$ on $\mathbb{N}^d$, then $p_{f,\mathrm{J}}(y\mid x)$ dominates $p_{\mathrm{J}}(y\mid x)$.

\item
Condition \eqref{eq:Fineq} is satisfied if $f\in\mathcal{C}^2([0,\infty)^d)$ satisfies the regularity conditions \eqref{eq:t1-condition-a}, \eqref{eq:t1-condition-b}, \eqref{eq:t1-condition-c} and the following conditions:
\begin{subequations}
\begin{align}
&\sum_{i=1}^{d}\frac{\partial^2 f}{\partial\theta_i^2}(\theta) \le0;
\label{eq:t1-condition-d}
\end{align}
and
\begin{align}
&\left.\frac{\partial f}{\partial\theta_i}(\theta)\right|_{\theta_i=0}\le0, \label{eq:t1-condition-e}
\end{align}
for every $i$.
\end{subequations}
\end{enumerate}
\end{theorem}

The proof of Theorem~\ref{Theorem 1} is presented in the Appendix.

Condition \eqref{eq:t1-condition-d} indicates that $f$ is a superharmonic function. Condition \eqref{eq:t1-condition-e} indicates that the derivative of $f$ on the boundary is nonpositive.
For any superharmonic function $h\in\mathcal{C}^2(\mathbb{R}^d)$, we note that 
\begin{align*}
    f(\theta)&=\sum_{a\in\{1,-1\}^d}h(a_1\theta_1,a_2\theta_2,\dots,a_d\theta_d)
\end{align*} satisfies the two conditions. Therefore, the relationship between the superharmonic function and improvement of the Bayesian predictive distribution based on the Jeffreys prior is obtained, as provided below.

\begin{corollary}\label{Corollary 1}
Let $\theta_i\coloneqq \sqrt{\lambda_i},\ i=1,\dots,d$. The Bayesian predictive distribution $p_{f,\mathrm{J}}(y\mid x)$ weakly dominates the Bayesian predictive distribution $p_{\mathrm{J}}(y\mid x)$ based on the Jeffreys prior if 
\begin{equation*}
f(\theta)=\sum_{a\in\{1,-1\}^d}h(a_1\theta_1,a_2\theta_2,\dots,a_d\theta_d),
\end{equation*}
where $h\in\mathcal{C}^2(\mathbb{R}^d)$ denotes a superharmonic function and $f$ satisfies regularity conditions \eqref{eq:t1-condition-a}, \eqref{eq:t1-condition-b}, and \eqref{eq:t1-condition-c}. Furthermore, if $F(z,r)$ is not a constant function of $z$ on $\mathbb{N}^d$, then $p_{f,\mathrm{J}}(y\mid x)$ dominates $p_{\mathrm{J}}(y\mid x)$.
\end{corollary}

For the multivariate normal model, George et al.~\cite{george2006improved} proved that Bayesian predictive distributions based on the superharmonic prior $h(\mu_1,\dots,\mu_d)\mathrm{d}\mu$ dominate those based on the Jeffreys prior. Here, we proved that the Bayesian predictive distributions based on ${\sum_{a\in\{1,-1\}^d}h(a_1\theta_1,a_2\theta_2,\dots,a_d\theta_d)\mathrm{d}\theta}$ dominate those based on the Jeffreys prior for the multivariate Poisson model. We have shown that Bayesian predictive distributions based on $h(\theta_1,\dots,\theta_d)\mathrm{d}\theta$ does not necessarily dominate those based on the Jeffreys prior using a numerical example provided at the beginning of this section. Hence, the results indicate the similarity and difference between the multivariate Poisson and normal models.

Theorem~\ref{Theorem 1} is divided into two parts to make it applicable to prior $f(\theta)$ that is not differentiable. For example, we can construct a series of differentiable priors
\begin{align*}
    f_n(\theta)=\sum_{a\in\{1,-1\}^3}\Big((a_1\theta_1-2)^2+(a_2\theta_2-2)^2+(a_3\theta_3-2)^2+\frac{1}{n}\Big)^{-1/2},
\end{align*} for the prior 
\begin{align*}
    f(\theta)=\sum_{a\in\{1,-1\}^3}\Big((a_1\theta_1-2)^2+(a_2\theta_2-2)^2+(a_3\theta_3-2)^2\Big)^{-1/2}.
\end{align*}  $f_n$ satisfies the conditions of the second part, thus $f_n$ satisfies the conditions of the first part. Subsequently, it is proven that $f$ also satisfies the conditions of the first part by using $f_n\overset{a.e.}\to f$.

The K-L risk of estimator $\hat\lambda$ is defined as the K-L divergence of $p(x\mid\lambda)$ and plug-in density $p(x\mid\hat\lambda)$, which is expressed as:
\begin{equation*}
\sum_x \biggl[ \sum_{i=1}^d
\Bigl\{ r\lambda_i\log \Bigl(\frac{\lambda_i}{\hat\lambda_i} \Bigr)
-r\lambda_i+r\hat\lambda_i \Bigr\} \biggr]
\prod_{j=1}^d\frac{(r\lambda_j)^{x_j}\exp(-r\lambda_j)}{x_j!}.
\end{equation*}
The Bayesian estimators based on $\pi_{\mathrm{J}}(\lambda)$ and $\pi_{f,\mathrm{J}}(\lambda)$ are known to be 
\begin{equation*}
\Bigl(\frac{x_1+1/2}{r},\ \frac{x_2+1/2}{r},\dots,\ \frac{x_d+1/2}{r}\Bigr)
\end{equation*}
and
\begin{equation*}
\Bigl(\frac{x_1+1/2}{r}\frac{F(x+\delta_1,r)}{F(x,r)},\ 
\frac{x_2+1/2}{r}\frac{F(x+\delta_2,r)}{F(x,r)},\dots,\ \frac{x_d+1/2}{r}\frac{F(x+\delta_d,r)}{F(x,r)}\Bigr),
\end{equation*}
respectively.
Therefore, the difference between the K-L risks of the Bayesian estimators based on $\pi_{\mathrm{J}}(\lambda)$ and $\pi_{f,\mathrm{J}}(\lambda)$ has the same sign as \eqref{eq:differential2}, which is provided in the Appendix. Using the proof of Theorem~\ref{Theorem 1}, the following can be obtained.

\begin{corollary}\label{Corollary 2}\leavevmode
\begin{enumerate}[label=\arabic*)]
\item
Let $\theta_i\coloneqq \sqrt{\lambda_i},\ i=1,\dots,d$. The Bayesian estimator based on $\pi_{f,\mathrm{J}}$ weakly dominates the Bayesian estimator based on $\pi_{\mathrm{J}}$ if for every $z\in \mathbb{N}^d$, $r>0$,
\begin{equation}
 \sum_{i=1}^dz_i \bigl\{ F(z,r)-F(z-\delta_i,r) \bigr\}
 +\sum_{i=1}^d(z_i+1/2) \bigl\{ F(z,r)-F(z+\delta_i,r) \bigr\} \ge 0. \label{Fineq2-2}
\end{equation}
Furthermore, if $F(z,r)$ is not a constant function of $z$ on $\mathbb{N}^d$, then the Bayesian estimator based on $\pi_{f,\mathrm{J}}$ dominates the Bayesian estimator based on $\pi_{\mathrm{J}}$.

\item
Condition \eqref{Fineq2-2} is satisfied if $f\in\mathcal{C}^2([0,\infty)^d)$ satisfies the regularity conditions \eqref{eq:t1-condition-a}, \eqref{eq:t1-condition-b}, \eqref{eq:t1-condition-c} and 
\begin{equation*}
f(\theta)=\sum_{a\in\{1,-1\}^d}h(a_1\theta_1,a_2\theta_2,\dots,a_d\theta_d),
\end{equation*}
where $h\in\mathcal{C}^2(\mathbb{R}^d)$ denotes a superharmonic function.
\end{enumerate}
\end{corollary}

From the corollary, the Bayesian estimator based on the Jeffreys prior can also be improved using superharmonic functions.

\section{Improved Prediction for Independent Poisson Processes with Different Durations}
\label{sec:different}

In this section, we consider the case of independent Poisson processes with different durations. Suppose that $x_i$ and $y_i\ (i = 1,\dots,d)$ are independently distributed according to Poisson distributions
\begin{equation*}
p(x\mid\lambda)=\prod_{i=1}^d\frac{(r_i\lambda_i)^{x_i}}{x_i!}e^{-r_i\lambda_i}
\end{equation*}
and
\begin{equation*}
p(y\mid\lambda)=\prod_{i=1}^d\frac{(s_i\lambda_i)^{y_i}}{y_i!}e^{-s_i\lambda_i}
\end{equation*}
with mean $r_i\lambda_i$ and $s_i\lambda_i$, respectively. If $r_1=r_2=\cdots=r_d$ and $s_1=s_2=\cdots=s_d$, it is a case of same duration which is discussed in the previous section.

Define $\gamma_i\coloneqq 1/r_i-1/(r_i+s_i)$ and $\theta_i\coloneqq \sqrt{\lambda_i/\gamma_i},\ i=1,\dots,d.$ Prior $\pi_{f,\mathrm{J}}(\lambda)$ is still defined as 
\begin{equation*}
\pi_{f,\mathrm{J}}(\lambda)\mathrm{d}\lambda_1\mathrm{d}\lambda_2\cdots \mathrm{d}\lambda_d\propto f(\theta_1,\dots,\theta_d)\lambda_1^{-1/2}\lambda_2^{-1/2}\cdots\lambda_d^{-1/2}\mathrm{d}\lambda_1\mathrm{d}\lambda_2\cdots \mathrm{d}\lambda_d.
\end{equation*}
Therefore, prior $\pi_{f,\mathrm{J}}(\lambda)\mathrm{d}\lambda$ is equivalent to prior $f(\theta)\mathrm{d}\theta$. Note that $\gamma_i$ is given and does not change in the following results.

We show that if $f$ satisfies certain conditions that are similar to the conditions in Section~\ref{sec:superharmonic}, the Bayesian predictive distribution $p_{f,\mathrm{J}}(y\mid x)$ dominates the Bayesian predictive distribution $p_{\mathrm{J}}(y\mid x)$. In this section, we define 
\begin{equation*}
F(z,r)\coloneqq F(z_1,z_2,\dots,z_d,r)= \int f(\theta_1,\dots,\theta_d)\prod_{i=1}^d\frac{r_i^{z_i+1/2}\lambda_i^{z_i-1/2}\exp(-r_i\lambda_i)}{\Gamma(z_i+1/2)}\mathrm{d}\lambda.
\end{equation*}

\begin{theorem}\label{Theorem 2}\leavevmode
\begin{enumerate}[label=\arabic*)]
\item
Let $\theta_i\coloneqq \sqrt{\lambda_i/\gamma_i},\ i=1,\dots,d$. The Bayesian predictive distribution $p_{f,\mathrm{J}}(y\mid x)$ weakly dominates the Bayesian predictive distribution $p_{\mathrm{J}}(y\mid x)$ if for every $z\in \mathbb{N}^d$, $r\in\mathbb{R}_+^d$,
\begin{equation}
    \sum_{i=1}^d\gamma_ir_iz_i \Bigl\{ F(z,r)-F(z-\delta_i,r)\Bigr\} +\sum_{i=1}^d\gamma_ir_i(z_i+1/2) \Bigl\{F(z,r)-F(z+\delta_i,r) \Bigr\} \ge 0,
\label{Fineq3}
\end{equation} 
where $F(z-\delta_i,r)\coloneqq 1$ if $z_i=0$.
Furthermore, if $F(z,r)$ is not a constant function of $z$ on $\mathbb{N}^d$, then $p_{f,\mathrm{J}}(y\mid x)$ dominates $p_{\mathrm{J}}(y\mid x)$.

\item
Condition \eqref{Fineq3} is satisfied if $f\in\mathcal{C}^2([0,\infty)^d)$ satisfies regularity conditions \eqref{eq:t1-condition-a}, \eqref{eq:t1-condition-b}, and \eqref{eq:t1-condition-c}, and conditions \eqref{eq:t1-condition-d} and \eqref{eq:t1-condition-e}.
\end{enumerate}
\end{theorem}

The proof of Theorem~\ref{Theorem 2} is a generalization of that of Theorem~\ref{Theorem 1}, which is presented in the Supplemental Material. The relationship between the superharmonic function and improvement of the Bayesian predictive distribution based on the Jeffreys prior is similar to Corollary~\ref{Corollary 1}, as provided below.

\begin{corollary}\label{Corollary 3}
Let $\theta_i\coloneqq \sqrt{\lambda_i/\gamma_i},\ i=1,\dots,d$.
The Bayesian predictive distribution $p_{f,\mathrm{J}}(y\mid x)$ weakly dominates the Bayesian predictive distribution $p_{\mathrm{J}}(y\mid x)$ based on the Jeffreys prior if 
\begin{equation*}
f(\theta)=\sum_{a\in\{1,-1\}^d}h(a_1\theta_1,a_2\theta_2,\dots,a_d\theta_d),
\end{equation*}
where $h\in\mathcal{C}^2(\mathbb{R}^d)$ denotes a superharmonic function and $f$ satisfies regularity conditions \eqref{eq:t1-condition-a}, \eqref{eq:t1-condition-b}, and \eqref{eq:t1-condition-c}. Furthermore, if $F(z,r)$ is not a constant function of $z$ on $\mathbb{N}^d$, then $p_{f,\mathrm{J}}(y\mid x)$ dominates $p_{\mathrm{J}}(y\mid x)$.
\end{corollary}

Therefore, the Bayesian predictive distribution based on prior $\sum_{a\in\{1,-1\}^d}h(a_1\theta_1,a_2\theta_2,\dots,a_d\theta_d)\mathrm{d}\theta$ dominates that based on the Jeffreys prior. We focus on the construction of a prior that is superior to the Jeffreys prior because of the results provided below.

\begin{theorem}\label{Theorem 3}\leavevmode
\begin{enumerate}[label=\arabic*)]
\item
For any $\lambda$, the K-L risk of $p_{\mathrm{J}}(y\mid x)$ is less than $0.52\sum_{i=1}^d\log((r_i+s_i)/r_i)$.
\item
For any predictive distribution $q(y\mid x)$ and positive number $\epsilon$, there exists $\lambda$ such that the K-L risk of $q(y\mid x)$ is greater than $0.5\sum_{i=1}^d\log((r_i+s_i)/r_i)-\epsilon$.
\end{enumerate}
\end{theorem}

The result is a generalization of Theorems 1 and 2 in Li \cite{li2024nearly}. The proof is presented in the Supplemental Material. According to the theorem, the upper bound of the K-L risk of $p_{\mathrm{J}}(y\mid x)$ is less than $1.04$ times the minimax lower bound. This motivates the following definition of a nearly minimax predictive distribution.

\begin{definition}\label{Def 1}
A predictive distribution $q(y\mid x)$ is labeled as nearly minimax if for any $\lambda$, the K-L risk of $q(y\mid x)$ is less than $1.04$ times the minimax lower bound.
\end{definition}

Hence, the Bayesian predictive distribution based on the Jeffreys prior is nearly minimax.

\section{Examples}
\label{sec:examples}

We provide examples that satisfy the conditions presented in Sections~\ref{sec:superharmonic} and \ref{sec:different}, including point and subspace shrinkage priors. The proofs of the propositions are presented in the Supplemental Material.

\subsection{Point Shrinkage Prior} 

The class of priors in Komaki \cite{komaki2004simultaneous} and \cite{komaki2015simultaneous} are considered in this study. These are
\begin{equation*}
\pi_{\alpha,\mathrm{J}}(\lambda)\mathrm{d}\lambda_1\mathrm{d}\lambda_2\cdots \mathrm{d}\lambda_d\propto\frac{\lambda_1^{-1/2}\lambda_2^{-1/2}\cdots\lambda_d^{-1/2}}{(\lambda_1+\lambda_2+\cdots+\lambda_d)^{\alpha}}\mathrm{d}\lambda_1\mathrm{d}\lambda_2\cdots \mathrm{d}\lambda_d
\end{equation*}
for the prediction of independent Poisson processes with the same duration, and
\begin{equation*}
\pi_{\alpha,\mathrm{J},\gamma}(\lambda)\mathrm{d}\lambda_1\mathrm{d}\lambda_2\cdots \mathrm{d}\lambda_d\propto\frac{\lambda_1^{-1/2}\lambda_2^{-1/2}\cdots\lambda_d^{-1/2}}{(\lambda_1/\gamma_1+\lambda_2/\gamma_2+\cdots+\lambda_d/\gamma_d)^{\alpha}}\mathrm{d}\lambda_1\mathrm{d}\lambda_2\cdots \mathrm{d}\lambda_d
\end{equation*}
for the prediction of independent Poisson processes with different durations, where $0<\alpha\le (d-2)/2$. The two priors are the same as $\pi_{f, \mathrm{J}}(\lambda)$ in this study, where $f(\theta)=\big(\sum_{i=1}^d\theta_i^2\big)^{-\alpha}$.

The Bayesian predictions based on prior $\pi_{f, \mathrm{J}}(\lambda)$ with $f(\theta)=\big(\sum_{i=1}^d\theta_i^2\big)^{-\alpha}$ shrink $\theta$ toward the origin. Therefore, it is natural to investigate the priors that shrink $\theta$ toward a general point $(\eta_1,\dots,\eta_d).$ Note that function $\big(\sum_{i=1}^d(\theta_i-\eta_i)^2\big)^{-\alpha}$ is superharmonic when $0<\alpha\le (d-2)/2$. Proposition~\ref{Proposition 2} is obtained using Theorem~\ref{Theorem 2}.

\begin{proposition}\label{Proposition 2} The Bayesian predictive distribution based on $\pi_{f,\mathrm{J}}(\lambda)$ dominates that based on the Jeffreys prior and is thus nearly minimax for the prediction of independent Poisson processes with the same duration or different durations when 
\begin{equation*}
f(\theta)=\sum_{a\in\{1,-1\}^d}\Big(\sum_{i=1}^d(a_i\theta_i-\eta_i)^2\Big)^{-\alpha},
\end{equation*}
where $0<\alpha\le (d-2)/2$.
\end{proposition}

Prior $\pi_{f,\mathrm{J}}(\lambda)\mathrm{d}\lambda$ is equivalent to prior $f(\theta)\mathrm{d}\theta$. Thus, prior $\pi_{f,\mathrm{J}}(\lambda)$ shrinks $\theta$ toward point $(\eta_1,\dots,\eta_d)$ if $\eta_i\ge0,\ i=1,\dots,d.$ From Corollary~\ref{Corollary 2}, it follows that the Bayesian estimator based on the point shrinkage prior dominates that based on the Jeffreys prior for a case with the same duration.

\subsection{Subspace Shrinkage Prior}  
In the previous examples, the Bayesian predictions shrink $\theta$ toward a point. Therefore, it is natural to investigate the subspace shrinkage prior that is constructed by the function $(s_V(\theta))^{-\alpha}$, where $s_V(\theta)$ represents the squared distance from $\theta$ to a linear subspace $V\subset \mathbb{R}^d$.

The complementary space is assumed to be $V^\perp=\text{span}(v_1,v_2,\dots,v_{d-k})$, where $\{v_1,v_2,\dots,v_{d-k}\}$ denotes a standard orthonormal basis and $k$ denotes the dimension of $V$. We have 
\begin{equation*}
s_V(\theta)=\sum_{i=1}^{d-k}\langle \theta,v_i \rangle ^2.
\end{equation*}
Note that function $(s_V(\theta))^{-\alpha}$ is superharmonic when $0<\alpha\le (d-k-2)/2$. Proposition~\ref{proposition 3} is obtained using Theorem~\ref{Theorem 2}.

\begin{proposition}\label{proposition 3} The Bayesian predictive distribution based on $\pi_{f,\mathrm{J}}(\lambda)$ dominates that based on the Jeffreys prior and is thus nearly minimax for the prediction of independent Poisson processes with the same duration or different durations when 
\begin{equation*}
f(\theta)=\sum_{a\in\{1,-1\}^d}(s_V(a_1\theta_1,a_2\theta_2,\dots,a_d\theta_d))^{-\alpha}\eqqcolon\sum_{a\in\{1,-1\}^d}(s_V(a\theta))^{-\alpha},
\end{equation*}
where $0<\alpha\le (d-k-2)/2$.
\end{proposition}

Prior $\pi_{f,\mathrm{J}}(\lambda)\mathrm{d}\lambda$ is equivalent to prior $f(\theta)\mathrm{d}\theta$. Thus, prior $\pi_{f,\mathrm{J}}(\lambda)$ shrinks $\theta$ toward subspace $V$ if $V\cap\mathbb{R}_+^d\neq\emptyset$.

\subsection{Mix Subspace Shrinkage Prior}  
The function $F$ in the first half of Theorem~\ref{Theorem 2} has additivity, which implies that $F$ equals the sum of $F_1$ and $F_2$ if $f=f_1+f_2$, which is expressed as:
\begin{equation*}
F(z,r)=\int (f_1+f_2)\prod_{i=1}^d\frac{r_i^{z_i+1/2}\lambda_i^{z_i-1/2}\exp(-r_i\lambda_i)}{\Gamma(z_i+1/2)}\mathrm{d}\lambda=F_1(z,r)+F_2(z,r).
\end{equation*}
Therefore, we can take the sum over different priors to obtain a new prior and ensure that the inequality condition \eqref{Fineq3} in the first half of Theorem~\ref{Theorem 2} is still satisfied. For example, if there are $n$ pairs of subspace $V_i$ and $\alpha_i$ satisfying the conditions of Proposition~\ref{proposition 3}, Proposition~\ref{proposition 4} is obtained using Theorem~\ref{Theorem 2}.

\begin{proposition}\label{proposition 4} The Bayesian predictive distribution based on $\pi_{f,\mathrm{J}}(\lambda)$ dominates that based on the Jeffreys prior and is thus nearly minimax for the prediction of independent Poisson processes with the same duration or different durations when 
\begin{equation*}
f(\theta)=\sum_{a\in\{1,-1\}^d}(s_{V_1}(a\theta))^{-\alpha_1}+\sum_{a\in\{1,-1\}^d}(s_{V_2}(a\theta))^{-\alpha_2}+\cdots+\sum_{a\in\{1,-1\}^d}(s_{V_n}(a\theta))^{-\alpha_n},
\end{equation*}
where $0<\alpha_i\le (d-k_i-2)/2$ and $k_i$ denotes the dimension of $V_i$, $i=1,2,\dots,n$.
\end{proposition}

We call this type of prior a ``mix subspace shrinkage prior." George \cite{george1986minimax} and George et al.~\cite{george2012minimax} studied a similar prior distribution for the normal model. From Corollary~\ref{Corollary 2}, it follows that the Bayesian estimators based on the subspace and mix subspace shrinkage priors dominate that based on the Jeffreys prior for a case with the same duration.

\section{Numerical Experiments and Real Data Application}
\label{sec:numerical}

We perform numerical experiments to demonstrate the difference between the risk of the Bayesian predictive distribution based on the Jeffreys prior and those based on the shrinkage priors discussed in Section~\ref{sec:examples}. We then use an application involving real data to compare the performance of different types of priors. Only cases with the same duration are considered in the experiments.

\subsection{Experiment 1}

We set $r=s=1$, $d=3$, and $\theta_i=\sqrt{\lambda_i},\ i=1,\dots,3$. The first prior is ``point shrinkage prior 1'' with 
\begin{equation*}
f(\theta)=\Big(\sum_{i=1}^{3}\theta_i^2\Big)^{-(3-2)/2},
\end{equation*}
which shrinks $\theta$ toward the origin. The second prior is ``point shrinkage prior 2'' with 
\begin{equation*}
f(\theta)=\sum_{a\in\{1,-1\}^d}\Big(\sum_{i=1}^{3}(a_i\theta_i-2)^2\Big)^{-(3-2)/2},
\end{equation*}
which shrinks $\theta$ toward point $(2,2,2)$. The third prior is a harmonic prior with 
\begin{equation*}
f(\theta)=\Big(\sum_{i=1}^{3}(\theta_i-2)^2\Big)^{-(3-2)/2}.
\end{equation*}

\begin{figure*}[!t]
\centering
\includegraphics[width=0.9\columnwidth]{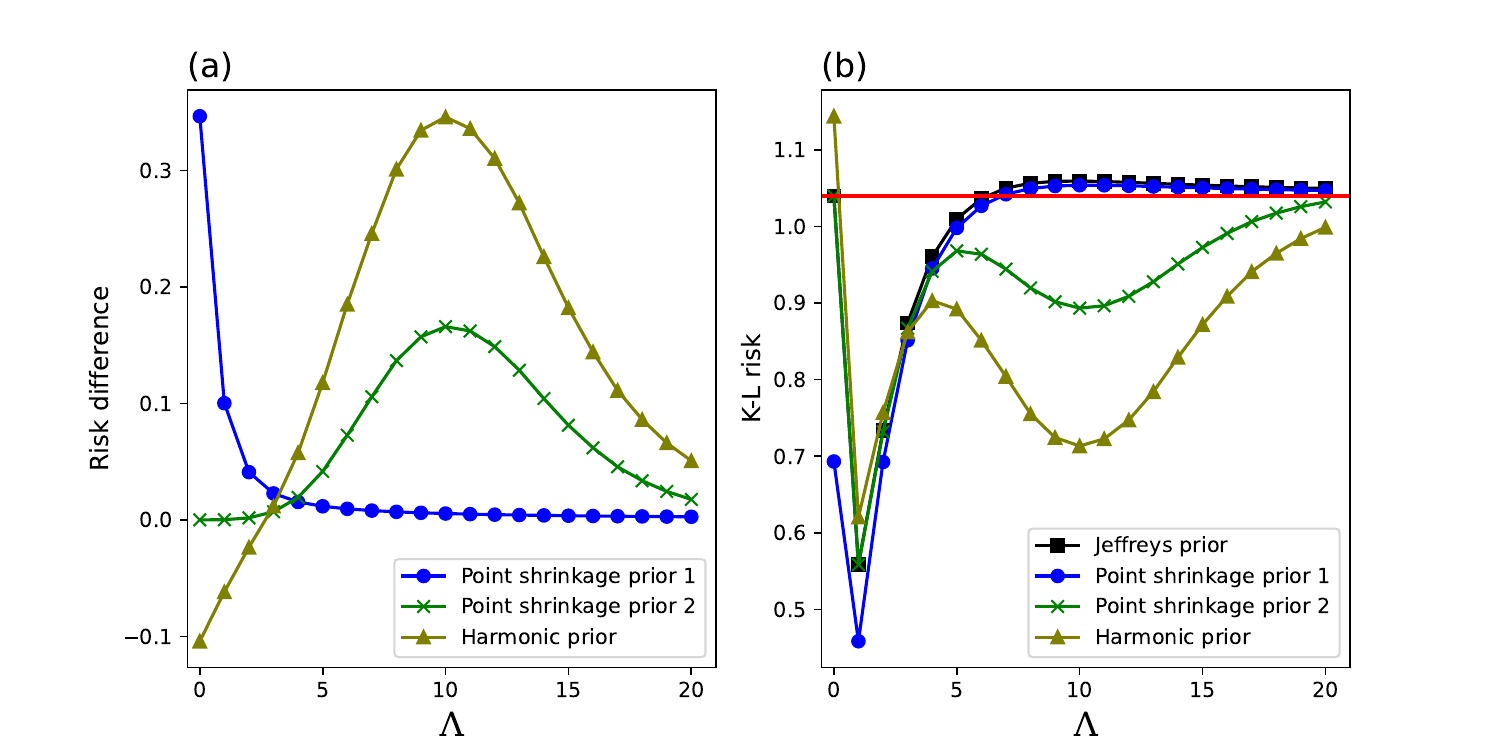}
\caption{(a) The K-L risk difference between $p_{\mathrm{J}}(y\mid x)$ and $p_{f,\mathrm{J}}(y\mid x)$. (b) The K-L risks of Bayesian predictive distributions based on different priors. The red horizontal line denotes the minimax lower bound $0.5d\log((r+s)/r)$.}
\label{point}
\end{figure*}

Fig.~\ref{point} shows the differences between the risks of the Bayesian predictive distributions based on the three priors and the Jeffreys prior when $\lambda=\Lambda\times(0.4,0.4,0.4)$. When $\Lambda$ is small, $\theta$ is close to the origin and the Bayesian predictive distribution based on the point shrinkage prior 1 performs better. When $\Lambda$ is close to $10$, $\theta$ is close to point $(2,2,2)$ and the Bayesian predictive distribution based on the point shrinkage prior 2 performs well.

It can be observed that the Bayesian predictive distribution based on the harmonic prior does not dominate that based on the Jeffreys prior. Specifically, the harmonic prior performs worse than the Jeffreys prior when $\lambda$ is close to the origin. However, Bayesian predictive distributions based on the harmonic prior dominate that based on the Jeffreys prior in the multivariate normal model. This example demonstrates the difference between the multivariate Poisson and normal models.

\subsection{Experiment 2}

We set $r=s=1$, $d=4$, and $\theta_i=\sqrt{\lambda_i},\ i=1,\dots,3$. The first prior is a point shrinkage prior with 
\begin{equation*}
f(\theta)=\Big(\sum_{i=1}^{4}\theta_i^2\Big)^{-(4-2)/2},
\end{equation*}
which shrinks $\theta$ toward the origin. The second prior is a subspace shrinkage prior with 
\begin{equation*}
f(\theta)=\sum_{a\in\{1,-1\}^d}(s_V(a\theta))^{-(4-3)/2},
\end{equation*}
which shrinks $\theta$ toward subspace $V=\text{span}((1,1,1,1))$.

\begin{figure}[!t]
\centering
\includegraphics[width=0.6\columnwidth]{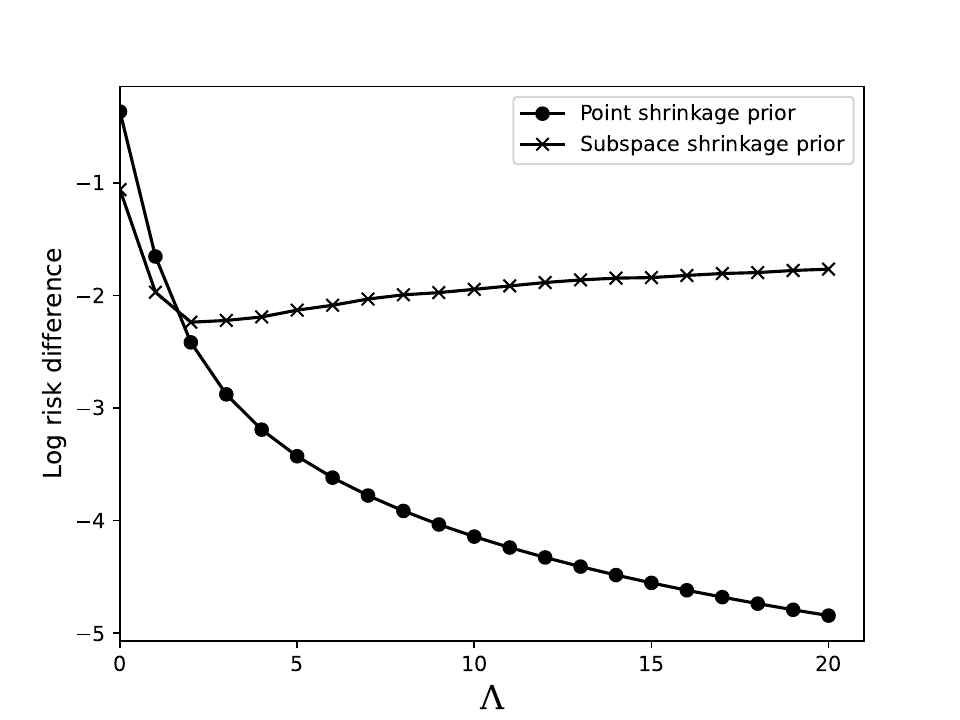}
\caption{Log value of the difference between the K-L risks of $p_{\mathrm{J}}(y\mid x)$ and $p_{f,\mathrm{J}}(y\mid x)$.}
\label{subspace}
\end{figure}

Fig.~\ref{subspace} shows the differences between the risks of the Bayesian predictive distributions based on the two priors and the Jeffreys prior when $\lambda=\Lambda\times(0.4,0.4,0.4,0.4)$. When $\Lambda$ is small, $\theta$ is close to the origin and the Bayesian predictive distribution based on the point shrinkage prior performs better. Note that $\theta$ here is always in the subspace $V=\text{span}((1,1,1,1))$. Therefore, the Bayesian predictive distribution based on the subspace shrinkage prior still performs well when $\Lambda$ is large.

\subsection{Experiment 3}

We set $r=s=1$, $d=4$, and $\theta_i=\sqrt{\lambda_i},\ i=1,\dots,4$. The first prior is a point shrinkage prior with 
\begin{equation*}
f(\theta)=(\theta_1^2+\theta_2^2+\theta_3^2+\theta_4^2)^{-(4-2)/2}.
\end{equation*}
The second prior is ``subspace shrinkage prior 1'' with 
\begin{equation*}
f(\theta)=(\theta_1^2+\theta_2^2+\theta_3^2)^{-(4-3)/2},
\end{equation*}
which shrinks $\theta$ toward subspace $V_1=\text{span}((0,0,0,1))$. The third prior is ``subspace shrinkage prior 2'' with 
\begin{equation*}
f(\theta)=(\theta_1^2+\theta_2^2+\theta_4^2)^{-(4-3)/2},
\end{equation*}
which shrinks $\theta$ toward subspace $V_2=\text{span}((0,0,1,0))$. The fourth prior is a ``mix subspace shrinkage prior'' with 
\begin{equation*}
f(\theta)=(\theta_1^2+\theta_2^2+\theta_3^2)^{-1/2}+(\theta_1^2+\theta_2^2+\theta_4^2)^{-1/2}+(\theta_1^2+\theta_3^2+\theta_4^2)^{-1/2}+(\theta_2^2+\theta_3^2+\theta_4^2)^{-1/2}.
\end{equation*}

\begin{figure}[!t]
\centering
\includegraphics[width=0.6\columnwidth]{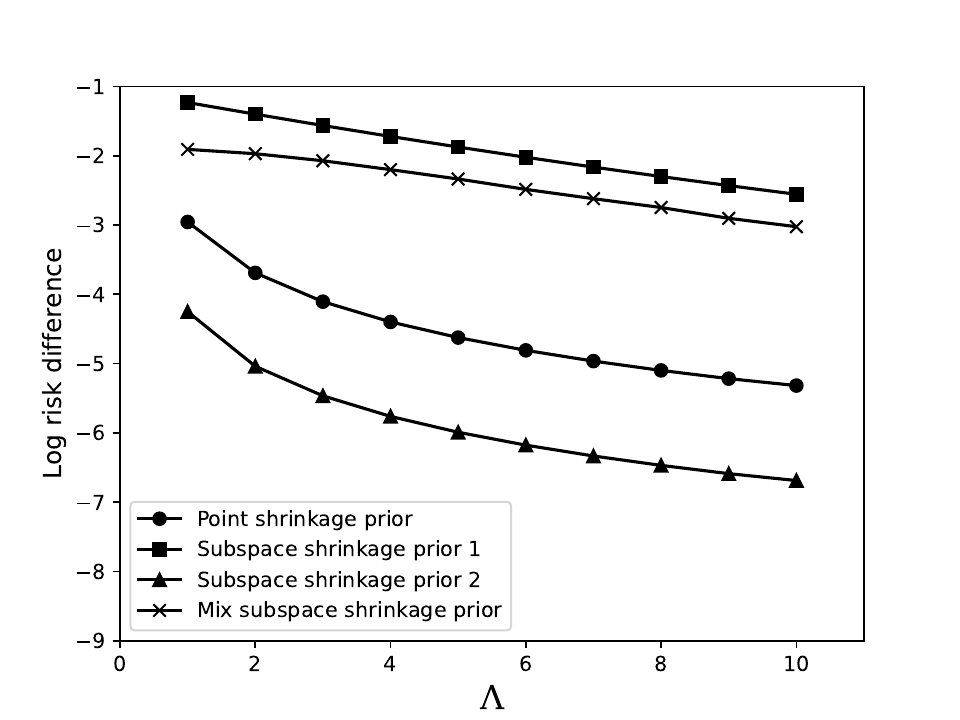}
\caption{Log value of the difference between the K-L risks of $p_{\mathrm{J}}(y\mid x)$ and $p_{f,\mathrm{J}}(y\mid x)$.}
\label{mix}
\end{figure}

Fig.~\ref{mix} shows the differences between the risks of the Bayesian predictive distributions based on the four priors and the Jeffreys prior when $\lambda=\Lambda\times(1,1,1,100)/20$. When $\Lambda$ is small, $\lambda$ is close to $\vec{0}$ and the Bayesian predictive distributions based on the point and subspace shrinkage priors perform well. Note that $\theta$ here is close to $V_1=\text{span}((0,0,0,1))$ but not close to $V_2=\text{span}((0,0,1,0))$ because $\lambda_4$ is significantly larger than the others. Therefore, the Bayesian predictive distribution based on subspace shrinkage prior 1 performs better than the others, and the Bayesian predictive distribution based on the mix subspace shrinkage prior is second best.

Note that $\lambda=\Lambda\times(1,1,1,100)/20$. If the choice among the subspace shrinkage priors is made by random selection when it is unknown which $\lambda_i$ is large, the performance of the subspace shrinkage prior $(\theta_i^2+\theta_j^2+\theta_k^2)^{-1/2}$ has $1/4$ and $3/4$ probabilities of being the performance of subspace shrinkage priors 1 and 2, respectively. In this case, the expected performance of the randomly selected subspace shrinkage prior is inferior to that of the mix subspace shrinkage prior. For example, when $\Lambda=5$, the risk reduction $\text{E}\bigl[ D(p(y\mid\lambda),p_{\mathrm{J}}(y\mid x))-D(p(y\mid\lambda),p_{\pi}(y\mid x))\,\big|\,\lambda\bigr] $ of subspace shrinkage priors 1 and 2 and the mix subspace shrinkage prior are $0.15$, $0.002$, and $0.10$, respectively. Note that $0.1>0.15\times 1/4+0.002\times 3/4=0.039$. Therefore, using the mix subspace shrinkage prior is recommended for this type of data. A similar example is considered in the context of a real data application.

\subsection{Real data application}

We utilized real data from an official database called \textit{the number of crimes in Tokyo by type and town} \cite{tokyocrime}, which reports the annual number of crimes in Tokyo. More appropriate measures for preventing crime can be implemented if the number of future crimes can be accurately predicted using past crime data.

\begin{figure*}[!t]
\centering
\includegraphics[width=0.8\columnwidth]{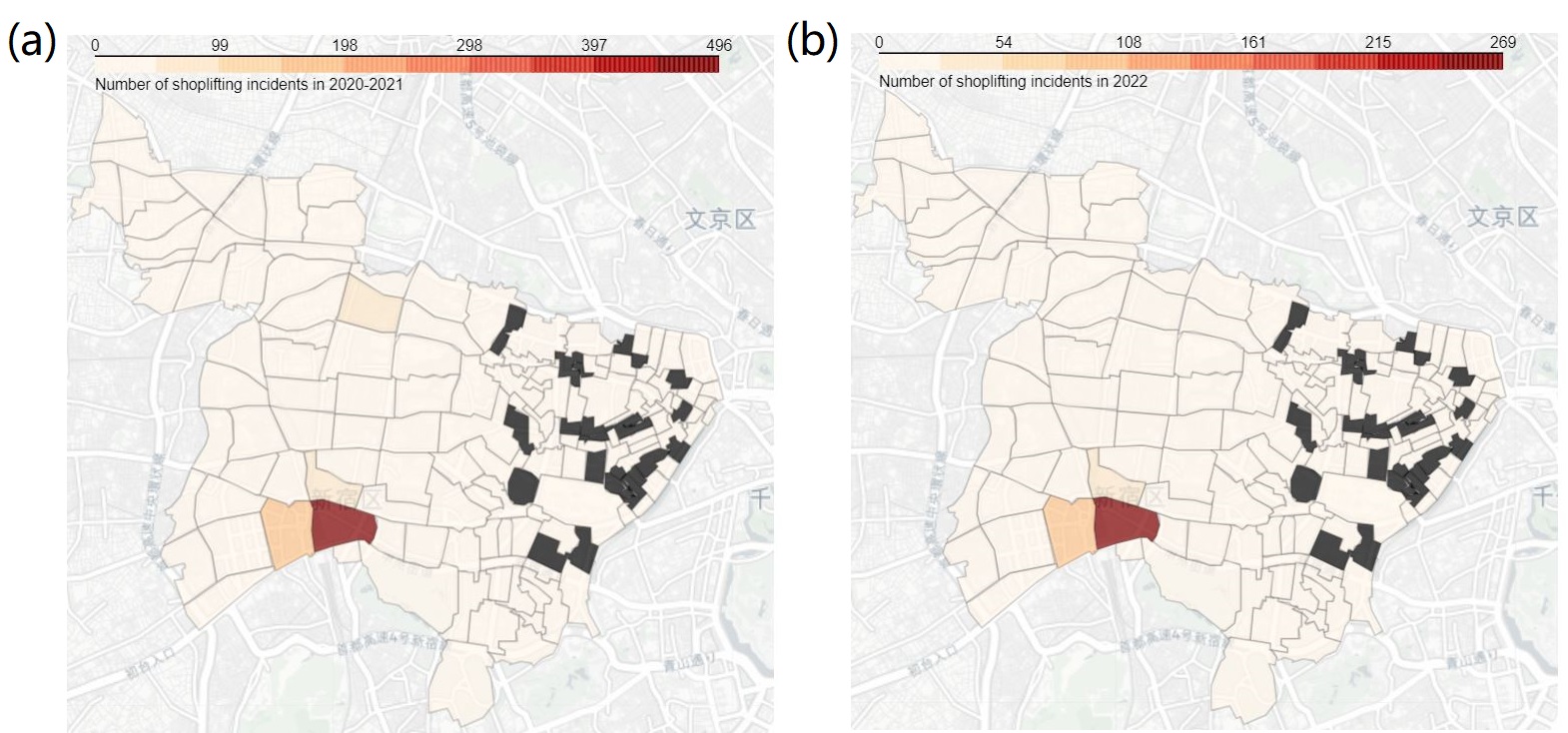}
\caption{Total number of shoplifting incidents in Shinjuku Ward for: (a) 2020 to 2021, and (b) 2022. Towns with incomplete data are shown in black.}
\label{shinjuku}
\end{figure*}

Shoplifting data for Shinjuku Ward was used from 2020--2022. After excluding towns with incomplete data, 130 towns were included. Fig.~\ref{shinjuku} shows the number of shoplifting incidents in Shinjuku Ward for 2020--2021 and 2022. Notably, one town had a significantly larger number of incidents than the others.

The experimental settings were as follows. The shoplifting data for 130 towns from 2020 to 2021 were set as observed data $x$ and the data from 2022 were set as future data $y$. The goal was to use $x$ to predict $y$. Therefore, the parameters in the prediction problem were $r=2$, $s=1$, and $d=130$. The Bayesian predictive distribution $p_{\pi}(y\mid x)$ based on prior $\pi$ was evaluated in three ways.
\begin{itemize}
\item Let $\hat y$ denote the mean vector of the predictive distribution. The K-L divergence between $\hat y$ and $y$ is expressed as $\sum_{i=1}^d \bigl(\hat y_i-y_i-y_i(\log \hat y_i-\log y_i)\bigr),$ where $0\log0=0$.
\item The weighted squared (W-S) distance between $\hat y$ and $y$ is expressed as $\sum_{i=1}^d (\hat y_i-y_i)^2/(y_i+1).$
\item The log-likelihood of $y$ in the predictive distribution is expressed as $\log p_{\pi}(y\mid x)$.
\end{itemize}

We compare the performance of the Jeffreys prior and the priors discussed in Section~\ref{sec:examples}, that is, a point shrinkage prior which shrinks $\theta$ toward the origin: 
\begin{equation*}
\Big(\sum_{i=1}^d\theta_i^2\Big)^{-(d-2)/2}\mathrm{d}\theta;
\end{equation*}
130 subspace shrinkage priors which shrink $\theta$ toward the subspaces $\{\theta\mid \theta_i=0,\forall i\neq j\},\ j=1,\dots,130$, i.e.,
\begin{equation*}
\Big(\sum_{i\neq j}\theta_i^2\Big)^{-(d-3)/2}\mathrm{d}\theta,\ j=1,\dots,130;
\end{equation*}
and the mix subspace shrinkage prior 
\begin{equation*}
\sum_{j=1}^d\Big(\sum_{i\neq j}\theta_i^2\Big)^{-(d-3)/2}\mathrm{d}\theta.
\end{equation*}
Table~\ref{tab:table1} presents a summary of the comparisons. As shown, the point and subspace shrinkage priors outperform the Jeffreys prior. For all the evaluation methods, the Bayesian predictive distribution based on the mix subspace shrinkage prior achieves the best scores.

\begin{table}[!t]
\caption{Comparison of Bayesian Predictive Distributions Based on Different Priors for Shoplifting Data Using the K-L Divergence, W-S Distance, and Predictive Log-Likelihood}
\label{tab:table1}
\centering
\begin{tabular}{c c c c c} 
 \hline
 &Jeffreys & Point shrinkage & Subspace shrinkage & Mix subspace shrinkage \\
 \hline
K-L divergence & 107.6 & 104.5 & Max:105.1, Min:100.9, Mean:104.5 & 100.9\\ 
 W-S distance & 259.4 & 235.5 & Max:240.0, Min:218.9, Mean:235.6 & 218.9\\
 Log-likelihood & $-$169.5 & $-$167.4 & Max:$-$164.9, Min:$-$167.8, Mean:$-$167.4 & $-$164.9 \\
 \hline
\end{tabular}
\end{table}

\section{Discussion}
\label{sec:discussion}

A point shrinkage prior $\pi_{f, \mathrm{J}}(\lambda)$ with $f(\theta)=\big(\sum_{i=1}^d\theta_i^2\big)^{-\alpha}$ was proposed in Komaki \cite{komaki2004simultaneous} and \cite{komaki2015simultaneous}. Additionally, the Bayesian predictive distribution based on the Jeffreys prior was shown to be dominated by that based on the point shrinkage prior if $0<\alpha\le d/2-1$. However, the underlying nature of the result was unclear. For example, the condition $0<\alpha\le d/2-1$ corresponds exactly to $f(\theta)$ being a superharmonic function, but the superharmonic function was not mentioned in the results of these two papers. 

In this study, the relationship between the superharmonic functions and improved predictive distribution are demonstrated. We show that the Bayesian predictive distribution based on the Jeffreys prior could be improved by prior $f(\theta)\mathrm{d}\theta$, if $f(\theta)=\sum_{a\in\{1,-1\}^d}h(a_1\theta_1,a_2\theta_2,\dots,a_d\theta_d)$ and $h$ is a superharmonic function. Sufficient conditions are also provided for independent Poisson processes with different durations. The results indicate the similarity and difference between the multivariate Poisson and normal models.

The results of this study help to discover different classes of priors that dominate the Jeffreys prior, such as the point shrinkage prior with $f(\theta)=\sum_{a\in\{1,-1\}^d}\big(\sum_{i=1}^d(a_i\theta_i-\eta_i)^2\big)^{-\alpha}$ and the subspace shrinkage prior with $f(\theta)=\sum_{a\in\{1,-1\}^d}(s_V(a\theta))^{-\alpha}$. The point shrinkage prior $f(\theta)=\big(\sum_{i=1}^d\theta_i^2\big)^{-\alpha}$ in the previous study performs well only when $\sum_{i=1}^d\lambda_i$ is small. In other cases, other priors discussed in Section~\ref{sec:examples} are more effective. In Experiment 1, the point shrinkage prior performs well when $\theta$ is close to a point (i.e., $\lambda$ is close to a point). In Experiment 2, the subspace shrinkage prior performs well when $\theta$ is close to a subspace (e.g., ${\lambda_1=\lambda_2=\cdots=\lambda_d}$). In Experiment 3 and the real data application, there is a $\lambda_i$ that is far larger than the others. For this type of data, it is shown that the mix subspace shrinkage prior with $f(\theta)=\sum_{j=1}^d\big(\sum_{k\neq j}\theta_k^2\big)^{(3-d)/2}$ performs well. Moreover, the mix subspace shrinkage prior does not require knowledge of the index $i$ of the large $\lambda_i$. These improved predictions based on different priors are nearly minimax, i.e., their K-L risk is less than 1.04 times the minimax lower bound.

On the basis of the results of this study, other types of prior distributions that dominate the Jeffreys prior will be constructed in the future. Following Stein \cite{stein1981estimation}, one may construct functions $f$ satisfying the conditions of Theorems~\ref{Theorem 1} and \ref{Theorem 2} through
\begin{equation*}
f(\theta) = \int \Vert\theta - \eta\Vert^{2-d}\mathrm{d}\rho(\eta),
\end{equation*}
where $\rho$ is a finite measure symmetric about each coordinate axis (equivalently, $\rho$ is invariant under all sign changes of coordinates). Different types of prior distributions are expected to perform well under different values of $\lambda$. It is also of interest to determine how to relax the sufficient conditions in Theorems~\ref{Theorem 1} and \ref{Theorem 2}, thereby broadening the class of applicable priors.

\section*{Acknowledgments}
The first author is grateful for support from the China Scholarship Council. We thank Keisuke Yano for his helpful comments. This work was supported by JSPS KAKENHI Grant Number 22H00510, and AMED Grant Numbers JP23dm0207001 and JP23dm0307009.

\appendices

\section{Proof of Theorem~\ref{Theorem 1}}
\label{app:proof theorem1}

\subsection{Proof of Part 1}
\label{app:proof part1}

We show this in two steps. Let $\mathrm{Po}(\lambda)$ denote the Poisson distribution with parameter $\lambda$. In Step~1, we show that the difference between the K-L risks of Bayesian predictive distributions based on $\pi_{f,\mathrm{J}}$ and $\pi_{\mathrm{J}}$ is 
\begin{equation}
\mathrm{E} \Bigl( \log(F(Z,r)) \, \Big| \, Z_i\sim\mathrm{Po}(r\lambda_i),i=1,\dots,d \Bigr)
-\mathrm{E} \Bigl( \log(F(Z,r+s)) \, \Big| \, Z_i\sim\mathrm{Po}((r+s)\lambda_i),i=1,\dots,d \Bigr).
\end{equation}
In Step~2, we prove that the partial derivative of $\mathrm{E} \bigl( \log(F(Z,r)) \, \big| \, Z_i\sim\mathrm{Po}(r\lambda_i),i=1,\dots,d \bigr)$ with respect to $r$ is nonnegative, and is positive when $F(z,r)$ is not a constant function in $z$. The details of each step are presented below.

\subsubsection*{Step 1}
Let $\bar\theta_j \coloneqq  \sqrt{\bar\lambda_j}$ $(j=1,\dots,d)$. The difference between the K-L risks of Bayesian predictive distributions based on $\pi_{f,\mathrm{J}}$ and $\pi_{\mathrm{J}}$ is
\begin{align}
&\mathrm{E} \biggl(
\log\frac{p_{\mathrm{J}}(y\mid x)}{p_{f,\mathrm{J}}(y\mid x)} 
\, \bigg| \, \lambda
\biggr) \notag\\
&=\mathrm{E}
\Biggl( \log\frac{\int p(x\mid\bar\lambda)\pi_{f,\mathrm{J}}(\bar\lambda) \mathrm{d}\bar\lambda}{\int p(x\mid\bar\lambda)\pi_{\mathrm{J}}(\bar\lambda) \mathrm{d}\bar\lambda}
\, \Bigg| \, \lambda \Biggr)
-\mathrm{E} \Biggl( \log\frac{\int p(x,y\mid\bar\lambda)\pi_{f,\mathrm{J}}(\bar\lambda) \mathrm{d}\bar\lambda}{\int p(x,y\mid\bar\lambda)\pi_{\mathrm{J}}(\bar\lambda) \mathrm{d}\bar\lambda}
\, \Bigg| \, \lambda \Biggr) \notag\\
&=\mathrm{E} \biggl( \log\int f(\bar\theta_1,\dots,\bar\theta_d)\prod_{i=1}^d
\frac{r^{x_i+1/2}\bar\lambda_i^{x_i-1/2}
\exp(-r\bar\lambda_i)}{\Gamma(x_i+1/2)}\mathrm{d}\bar\lambda
\, \bigg| \, \lambda \biggr) \notag\\
& \quad -\mathrm{E} \Biggl(\log\int f(\bar\theta_1,\dots,\bar\theta_d)\prod_{i=1}^d\frac{(r+s)^{x_i+y_i+1/2}\bar\lambda_i^{x_i+y_i-1/2}\exp\{-(r+s)\bar\lambda_i\}}{\Gamma(x_i+y_i+1/2)}\mathrm{d}\bar\lambda \, \Bigg| \, \lambda \Biggr) \notag\\
& = \mathrm{E} \biggl( \log(F(Z,r)) \, \bigg| \, Z_i\sim\mathrm{Po}(r\lambda_i),i=1,\dots,d \biggr)
-\mathrm{E} \biggl( \log(F(Z,r+s)) \, \bigg| \, Z_i\sim\mathrm{Po}((r+s)\lambda_i),i=1,\dots,d \biggr).
\label{eq:negativediff}
\end{align}

\subsubsection*{Step 2}
The risk difference \eqref{eq:negativediff} is nonpositive if
\begin{equation}
\mathrm{E} \biggl( \log F(Z,r) \, \bigg| \, Z_i\sim\mathrm{Po}(r\lambda_i),i=1,\dots,d \biggr) =\sum_z\log(F(z,r))\Biggl\{
\prod_{i=1}^d \frac{(r\lambda_i)^{z_i}\exp(-r\lambda_i)}{z_i!}
\Biggr\}
\label{eq:increasing}
\end{equation}
is a nondecreasing function of $r$.

Differentiating $F(z,r)$ with respect to $r$ yields
\begin{equation}
\frac{\partial F}{\partial r}(z,r)=\int f(\bar\theta_1,\dots,\bar\theta_d) \prod_{i=1}^d
 \frac{r^{z_i+1/2}\bar\lambda_i^{z_i-1/2}\exp(-r\bar\lambda_i)}{\Gamma(z_i+1/2)}
 \Bigl( \sum_{i=1}^d \frac{z_i+1/2}{r} - \sum_{i=1}^d \bar\lambda_i \Bigr) \mathrm{d}\bar\lambda.
\end{equation}
 
Differentiating \eqref{eq:increasing} with respect to $r$ gives
\begin{align}
&\sum_z \log(F(z,r))
\Biggl\{
\prod_{i=1}^d \frac{(r\lambda_i)^{z_i}\exp(-r\lambda_i)}{z_i!}
\Biggr\}
\Bigl(\sum_{i=1}^d \frac{z_i}{r} - \sum_{i=1}^d \lambda_i \Bigr) \notag \\
&\quad + \sum_z 
 \Biggl\{
 \frac{ 
 \int f(\bar\theta_1,\dots,\bar\theta_d) \prod_{i=1}^d
 \frac{r^{z_i+1/2}\bar\lambda_i^{z_i-1/2}\exp(-r\bar\lambda_i)}{\Gamma(z_i+1/2)}
 \Bigl( \sum_{i=1}^d \frac{z_i+1/2}{r} - \sum_{i=1}^d \bar\lambda_i \Bigr) \mathrm{d}\bar\lambda}
 {F(z,r)}
 \Biggr\}
 \prod_{i=1}^d \frac{(r\lambda_i)^{z_i} \exp(-r\lambda_i)}{z_i!}.
 \label{eq:differential}
\end{align}

We note that
\begin{align}
&\sum_z \Bigl\{ \log F(z_1,z_2,\dots,z_d,r) \Bigr\}
\Biggl\{ \prod_{i=1}^d\frac{(r\lambda_i)^{z_i} \exp(-r\lambda_i)}{z_i!} \Biggr\}
\sum_{i=1}^d \lambda_i
\notag \\
&= \sum_{i=1}^d \sum_z \biggl[ \frac{z_i}{r} \Bigl\{ \log F(z-\delta_i,r) \Bigr\}
\prod_{j=1}^d\frac{(r\lambda_j)^{z_j}\exp(-r\lambda_j)}{z_j!} \biggr],
\label{eq:transform1}
\end{align}
where $\delta_{ij}$ is defined as 1 if $i = j$ and 0 if $i\neq j$, $\delta_i$ is defined as the $d$-dimensional vector whose $i^{\mathrm{th}}$ element is 1 and all other elements are 0, and $F(z-\delta_i,r)$ is defined as $1$ if $z_i=0$.

Furthermore, we note that
\begin{equation}
  \int f(\bar\theta_1,\dots,\bar\theta_d)\prod_{j=1}^d\frac{r^{z_j+1/2}\bar\lambda_j^{z_j-1/2}
  \exp(-r\bar\lambda_j)}{\Gamma(z_j+1/2)}
  \bar\lambda_i \mathrm{d}\bar\lambda  = F(z+\delta_i,r)\frac{z_i+1/2}{r}.  
  \label{eq:transform2}
\end{equation}

Thus, from \eqref{eq:transform1} and \eqref{eq:transform2}, the partial derivative \eqref{eq:differential} of \eqref{eq:increasing} with respect to $r$ is
\begin{align}
&\sum_z\sum_{i=1}^d \biggl[ \frac{z_i}{r} \biggl\{ \log F(z,r) \biggr\}
\prod_{j=1}^d\frac{(r\lambda_j)^{z_j}\exp(-r\lambda_j)}{z_j!} \biggr]
-\sum_z \sum_{i=1}^d \biggl[ \frac{z_i}{r} \biggl\{ \log F(z-\delta_i,r) \biggr\}
\prod_{j=1}^d \frac{(r\lambda_j)^{z_j} \exp(-r\lambda_j)}{z_j!} \biggr] \notag\\
&\quad+ \sum_z \Biggl\{
\frac{F(z,r) \sum_{i=1}^d (z_i+1/2)/r
-\sum_{i=1}^d F(z+\delta_i,r) (z_i+1/2)/r}{F(z,r)} \Biggr\}
\prod_{j=1}^d\frac{(r\lambda_j)^{z_j}\exp(-r\lambda_j)}{z_j!} \notag\\
&= \sum_z
\Biggl[ \sum_{i=1}^d \frac{z_i}{r} \biggl\{ \log 
\frac{F(z,r)}{F(z-\delta_i,r)}\biggr\}
\prod_{j=1}^d\frac{(r\lambda_j)^{z_j}\exp(-r\lambda_j)}{z_j!}
\Biggr] \notag\\
&\quad+\sum_z \Biggl[ \sum_{i=1}^d
\frac{z_i+1/2}{r}
\biggl\{ 1-\frac{F(z+\delta_i,r)}{F(z,r)} \biggr\}
\prod_{j=1}^d
\frac{(r\lambda_j)^{z_j}\exp(-r\lambda_j)}{z_j!} \Biggr].
\label{eq:differential2}
\end{align}

By the inequality $\log \epsilon \geq 1-1/\epsilon$ for all $\epsilon > 0$, with equality if and only if $\epsilon = 1$, \eqref{eq:differential2} is larger than or equal to
\begin{equation}
\sum_z \Biggl[ \sum_{i=1}^d \frac{z_i}{r} \biggl\{ 1-\frac{F(z-\delta_i,r)}{F(z,r)} \biggr\} +\sum_{i=1}^d\frac{z_i+1/2}{r} \biggl\{ 1-\frac{F(z+\delta_i,r)}{F(z,r)} \biggr\}\Biggr]
\prod_{j=1}^d \frac{(r\lambda_j)^{z_j}\exp(-r\lambda_j)}{z_j!}. 
\label{eq:Fpositive}
\end{equation}

From \eqref{eq:Fineq}, we know that \eqref{eq:Fpositive} is nonnegative. Thus, \eqref{eq:differential2} is nonnegative, and \eqref{eq:increasing} is a nondecreasing function. This proves that $p_{f,\mathrm{J}}(y\mid x)$ weakly dominates $p_{\mathrm{J}}(y\mid x)$.

Furthermore, if $F(z,r)$ is not a constant function of $z$ on $\mathbb{N}^d$, then $F(z-\delta_i,r) \not\equiv F(z,r)$ for some $i$. Therefore, the inequality $\log \epsilon \geq 1-1/\epsilon$ is strict for some terms in the sum, and \eqref{eq:differential2} is strictly positive. Thus, \eqref{eq:increasing} is strictly increasing, and $p_{f,\mathrm{J}}(y\mid x)$ dominates $p_{\mathrm{J}}(y\mid x)$.

\subsection{Proof of Part 2}
\label{app:proof part2}

We prove that \eqref{eq:Fineq} is satisfied if $f$ satisfies the conditions of the second half of Theorem~\ref{Theorem 1}. Let $\theta_j \coloneqq  \sqrt{\lambda_j}$ $(j=1,\dots,d)$. We show this in three steps. \eqref{eq:Fineq} is obtained by combining Steps~2 and~3.

In Step~1, through integration by parts on $\theta_i$, we prove that 
\begin{equation}
F(z+\delta_i,r)-F(z,r) = 2^{d-1}\int  \frac{\partial f}{\partial\theta_i}(\theta)  \prod_{j=1}^d
\frac{r^{z_j+1/2}\theta_j^{2z_j+\delta_{ji}}\exp(-r\theta_j^2)}
{\Gamma(z_j+\delta_{ji}+1/2)}\mathrm{d}\theta.
    \label{eq:F-1st-diff}
\end{equation}

In Step~2, by performing integration by parts on each $\theta_i$ again and using the condition \eqref{eq:t1-condition-d} of the second derivative, we prove that
\begin{align}
&\sum_{i=1}^d z_i \{ F(z,r)-F(z-\delta_i,r) \}
 +\sum_{i=1}^d (z_i+1/2) \{F(z,r)-F(z+\delta_i,r)\} \notag\\
&\ge \sum_{i=1}^d \frac{2^{d-2}}{r} \Biggl[ \int_{[0,\infty)^{d-1}} \frac{\partial f}{\partial\theta_i}(\theta)
\biggl\{\prod_{j\neq i}\frac{r^{z_j+1/2}\theta_j^{2z_j}\exp(-r\theta_j^2)}{\Gamma(z_j+1/2)} \mathrm{d}\theta_j\biggr\} \frac{r^{z_i+1/2}\theta_i^{2z_i}\exp(-r\theta_i^2)}{\Gamma(z_i+1/2)} \Biggr]^{\infty}_{\theta_i = 0}. 
\label{eq:thm1-step2}
\end{align}

In Step~3, using the condition \eqref{eq:t1-condition-e} of the derivative on the boundary, we show that \eqref{eq:thm1-step2} $\ge 0$.
The details of each step are presented below.

\subsubsection*{Step 1}
Using the substitution $\theta=\sqrt{\lambda}$ and the definition of the function $F$, we obtain
\begin{align*}
&F(z+\delta_i,r)-F(z,r) \notag\\
&=-2^{d-1}\int\biggl\{
\int f(\theta)\prod_{j\neq i}\frac{r^{z_j+1/2}\theta_j^{2z_j}\exp(-r\theta_j^2)}
{\Gamma(z_j+1/2)} \mathrm{d} \theta_j 
\biggr\}
\frac{r^{z_i+1/2}\theta_i^{2z_i+1}\frac{\partial}{\partial\theta_i}\exp(-r\theta_i^2)}
{\Gamma(z_i+3/2)}\mathrm{d}\theta_i \notag\\
& \quad -2^{d-1}\int\biggl\{
\int f(\theta)\prod_{j\neq i}\frac{r^{z_j+1/2}\theta_j^{2z_j}\exp(-r\theta_j^2)}
{\Gamma(z_j+1/2)} \mathrm{d} \theta_j 
\biggr\}
\frac{r^{z_i+1/2}(\frac{\partial}{\partial\theta_i}\theta_i^{2z_i+1})\exp(-r\theta_i^2)}{\Gamma(z_i+3/2)}\mathrm{d}\theta_i.
\end{align*}

By performing integration by parts on $\theta_i$, we obtain
\begin{align}
&F(z+\delta_i,r)-F(z,r) \notag\\
&= 2^{d-1} \int_{[0,\infty)^d} \biggl\{ \frac{\partial}{\partial\theta_i}f(\theta) \biggr\}
\prod_{j=1}^d\frac{r^{z_j+1/2}
\theta_j^{2z_j+\delta_{ji}}\exp(-r\theta_j^2)}{\Gamma(z_j+\delta_{ji}+1/2)} \mathrm{d}\theta \notag\\
& \quad -2^{d-1} \Biggl[ \biggl\{
\int_{[0,\infty)^{d-1}} f(\theta)\prod_{j\neq i}\frac{r^{z_j+1/2}\theta_j^{2z_j}\exp(-r\theta_j^2)}
{\Gamma(z_j+1/2)} \mathrm{d} \theta_j 
\biggr\}
\frac{r^{z_i+1/2}\theta_i^{2z_i+1}\exp(-r\theta_i^2)}
{\Gamma(z_i+3/2)} \Biggr]^{\infty}_{\theta_i = 0}. 
\label{eq:F-diff-1}
\end{align}

As $\theta_i \to \infty$, the boundary term vanishes due to exponential decay $\exp(-r\theta_i^2)$ dominating polynomial growth $\theta_i^{2z_i+1}$. As $\theta_i \to 0$, the boundary term vanishes because $\theta_i^{2z_i+1} \to 0$ for $z_i \geq 0$. Thus, using \eqref{eq:F-diff-1}, we obtain \eqref{eq:F-1st-diff}.

\subsubsection*{Step 2}
Using \eqref{eq:F-1st-diff}, we obtain
\begin{align}
&\sum_{i=1}^d  z_i \{ F(z,r)-F(z-\delta_i,r) \}
 +\sum_{i=1}^d (z_i+1/2) \{F(z,r)-F(z+\delta_i,r)\} \notag\\
&=\sum_{i=1}^d 2^{d-1} \int\biggl[  \frac{\partial f }{\partial\theta_i}(\theta)  z_i
\prod_{j=1}^d\frac{r^{(z_j-\delta_{ji})+1/2}\theta_j^{2(z_j-\delta_{ji})+\delta_{ji}}\exp(-r\theta_j^2)}
{\Gamma(z_j+1/2)} \notag\\
& \qquad-  \frac{\partial f }{\partial\theta_i} (\theta)  (z_i+1/2)
\prod_{j=1}^d
\frac{r^{z_j+1/2}\theta_j^{2z_j+\delta_{ji}}\exp(-r\theta_j^2)}{\Gamma(z_j+\delta_{ji}+1/2)} \biggr]\mathrm{d}\theta \notag\\
&=\sum_{i=1}^d \frac{2^{d-2}}{r} \Biggl\{ \int_{[0,\infty)^d}  \frac{\partial f}{\partial\theta_i}(\theta)
\biggl\{\prod_{j\neq i}\frac{r^{z_j+1/2}\theta_j^{2z_j}\exp(-r\theta_j^2)}{\Gamma(z_j+1/2)} \mathrm{d}\theta_j\biggr\}
\frac{r^{z_i+1/2}\frac{\partial\theta_i^{2z_i}}{\partial\theta_i}\exp(-r\theta_i^2)}{\Gamma(z_i+1/2)}\mathrm{d}\theta_i \notag\\
& \qquad+ \int_{[0,\infty)^d}\frac{\partial f}{\partial\theta_i}(\theta)
\biggl\{\prod_{j\neq i}\frac{r^{z_j+1/2}\theta_j^{2z_j}\exp(-r\theta_j^2)}{\Gamma(z_j+1/2)}\mathrm{d}\theta_j\biggr\}
\frac{r^{z_i+1/2}\theta_i^{2z_i}\frac{\partial}{\partial\theta_i}\exp(-r\theta_i^2)}
{\Gamma(z_i+1/2)}\mathrm{d}\theta_i \Biggr\}. 
\label{eq:F-diff-2-1}
\end{align}

By performing integration by parts on each parameter again, we find that \eqref{eq:F-diff-2-1} equals
\begin{align}
&\sum_{i=1}^d \frac{2^{d-2}}{r} \Biggl\{
-\int_{[0,\infty)^d}  \frac{\partial}{\partial\theta_i}\Bigl\{\frac{\partial f}{\partial\theta_i}(\theta)\Bigr\}
\prod_{j\neq i}\frac{r^{z_j+1/2}\theta_j^{2z_j}\exp(-r\theta_j^2)}{\Gamma(z_j+1/2)}\mathrm{d}\theta_j \frac{r^{z_i+1/2}\theta_i^{2z_i}\exp(-r\theta_i^2)}{\Gamma(z_i+1/2)}\mathrm{d}\theta_i \notag \\
&\qquad+ \Biggl[\int_{[0,\infty)^{d-1}}  \frac{\partial f}{\partial\theta_i}(\theta) 
\biggl\{\prod_{j\neq i}\frac{r^{z_j+1/2}\theta_j^{2z_j}\exp(-r\theta_j^2)}{\Gamma(z_j+1/2)}\mathrm{d}\theta_j\biggr\} \frac{r^{z_i+1/2}\theta_i^{2z_i}\exp(-r\theta_i^2)}{\Gamma(z_i+1/2)} \Biggr]^{\infty}_{\theta_i = 0} \Biggr\}. 
\label{eq:F-diff-2}
\end{align}

From condition \eqref{eq:t1-condition-d}, we have
\begin{align}
&\sum_{i=1}^d \Biggl\{
-\int_{[0,\infty)^d} \frac{\partial^2 f}{\partial\theta_i^2}(\theta)
\prod_{j\neq i}\frac{r^{z_j+1/2}\theta_j^{2z_j}\exp(-r\theta_j^2)}{\Gamma(z_j+1/2)}\mathrm{d}\theta_j \frac{r^{z_i+1/2}\theta_i^{2z_i}\exp(-r\theta_i^2)}{\Gamma(z_i+1/2)}\mathrm{d}\theta_i \Biggr\} \notag\\
&=- \int_{[0,\infty)^d}  \sum_{i=1}^d \frac{\partial^2 f}{\partial\theta_i^2}(\theta)
\prod_{j=1}^d \frac{r^{z_j+1/2}\theta_j^{2z_j} \exp(-r\theta_j^2)}{\Gamma(z_j+1/2)}\mathrm{d}\theta \ge 0. 
\label{eq:F-2nd-diff-1}
\end{align}

Using \eqref{eq:F-diff-2-1}, \eqref{eq:F-diff-2}, and \eqref{eq:F-2nd-diff-1}, we obtain the inequality \eqref{eq:thm1-step2}.

\subsubsection*{Step 3}
We evaluate the boundary terms in \eqref{eq:thm1-step2}. As $\theta_i \to \infty$, the term vanishes due to condition \eqref{eq:t1-condition-b}. As $\theta_i \to 0$, using condition \eqref{eq:t1-condition-e} that $\frac{\partial f}{\partial\theta_i}\big|_{\theta_i=0} \leq 0$, we have: if $z_i > 0$, the factor $\theta_i^{2z_i} \to 0$ kills the boundary term; if $z_i = 0$, the boundary term becomes $\frac{\partial f}{\partial\theta_i}\big|_{\theta_i=0} \cdot (\text{positive factor}) \leq 0$. Therefore, \eqref{eq:thm1-step2} $\geq 0$, completing the proof.

\hfill\IEEEQED

\bibliographystyle{IEEEtran}
\bibliography{scholar}

\end{document}


\title{Supplement to ``Improved nearly minimax prediction for independent Poisson processes under Kullback--Leibler loss''}
\author[1]{Xiao Li}
\author[1,2]{Fumiyasu Komaki}
\affil[1]{\small Department of Mathematical Informatics, Graduate School of Information Science and Technology,\par The University of Tokyo, 7-3-1 Hongo, Bunkyo-ku, Tokyo, 113-0033, Japan}
\affil[2]{\small RIKEN Center for Brain Science, 2-1 Hirosawa, Wako City, Saitama, 351-0198, Japan}
\date{}
\maketitle

\begin{appendices}
\section{Proofs for results in the main paper}
\localtableofcontents

\subsection{Detailed Proof of Theorem 1}

Theorem 1 in the main paper is copied as follows.
\begin{theorem}\label{Theorem 1}\leavevmode
\begin{enumerate}[label=\arabic*)]
\item
Let $\theta_i\coloneqq \sqrt{\lambda_i},\ i=1,\dots,d$. The Bayesian predictive distribution $p_{f,\mathrm{J}}(y\mid x)$ weakly dominates the Bayesian predictive distribution $p_{\mathrm{J}}(y\mid x)$ if for every $z\in \mathbb{N}^d$, $r>0$,
\begin{align}
    &\sum\limits_{i=1}^dz_i(F(z,r)-F(z-\delta_i,r))+\sum\limits_{i=1}^d(z_i+1/2)(F(z,r)-F(z+\delta_i,r))\ge 0, \label{Fineq}
\end{align}
where $\delta_i$ denotes the $d$-dimensional vector whose $i^{\text{th}}$ element is 1, all other elements are 0, and $F(z-\delta_i,r)$ is defined as $1$ if $z_i=0$. 
Furthermore, if $F(z,r)$ is not a constant function of $z$ on $\mathbb{N}^d$, then $p_{f,\mathrm{J}}(y\mid x)$ dominates $p_{\mathrm{J}}(y\mid x)$.

\item
Condition \eqref{Fineq} is satisfied if $f\in\mathbf{C}^2([0,\infty)^d)$ satisfies the regularity conditions provided below and the following conditions:
\begin{subequations}
\begin{align}
&\sum\limits_{i=1}^{d}\frac{\partial^2 f}{\partial\theta_i^2}(\theta) \le0;
\label{t1-condition-d}
\end{align}
and
\begin{align}
&\left.\frac{\partial f}{\partial\theta_i}(\theta)\right|_{\theta_i=0}\le0, \label{t1-condition-e}
\end{align}
for every $i$.
\end{subequations}
\end{enumerate}
\end{theorem}

\noindent
The regularity conditions are expressed as:
\begin{subequations}
\begin{align}
&\int f(\theta) \prod\limits_{j=1}^{d}\theta_j^{2z_j}\exp(-r\theta_j^2)\du\theta<\infty, \label{t1-condition-a}\\
&\int \Bigl\lvert \frac{\partial f}{\partial\theta_i}(\theta) \Bigr\rvert \prod\limits_{j=1}^{d}\theta_j^{2z_j}\exp(-r\theta_j^2)\du\theta<\infty, \label{t1-condition-b}\\
\text{and} &\int \Bigl\lvert \frac{\partial^2 f}{\partial\theta_i^2}(\theta) \Bigr\rvert \prod\limits_{j=1}^{d}\theta_j^{2z_j}\exp(-r\theta_j^2)\du\theta<\infty, \label{t1-condition-c}
\end{align}
for every $z\in\mathbb{N}^d$, $r>0$, and $i$.
\end{subequations}

\vspace{0.2cm}
\noindent
\textit{\textbf{Proof of part 1.}}
We show this in two steps. Let $\text{Po}(\lambda)$ denote the Poisson distribution with parameter $\lambda$. In Step 1, we show that the difference between the K-L risks of Bayesian predictive distributions based on $\pi_{f,\mathrm{J}}$ and $\pi_{\mathrm{J}}$ is 
$$\text{E} \biggl( \log(F(Z,r)) \, \bigg| \, Z_i\sim\text{Po}(r\lambda_i),i=1,\dots,d \biggr)
-\text{E} \biggl( \log(F(Z,r+s)) \, \bigg| \, Z_i\sim\text{Po}((r+s)\lambda_i),i=1,\dots,d \biggr).$$
In Step 2, we prove that the partial derivative of $\text{E} \biggl( \log(F(Z,r)) \, \bigg| \, Z_i\sim\text{Po}(r\lambda_i),i=1,\dots,d \biggr)$ with respect to $r$ is nonnegative, and is positive when $F(z,r)$ is not a constant function in $z$. The details of each step are presented below.

\noindent
\textit{\textbf{Step 1.}} Let $\bar\theta_j \coloneqq  \sqrt{\bar\lambda_j}$ $(j=1,\dots,d)$. The difference between the K-L risks of Bayesian predictive distributions based on $\pi_{f,\mathrm{J}}$ and $\pi_{\mathrm{J}}$ is
\begin{align}
&\text{E} \biggl(
\log\frac{p_{\mathrm{J}}(y\mid x)}{p_{f,\mathrm{J}}(y\mid x)} 
\, \bigg| \, \lambda
\biggr) \notag\\
&=\text{E}
\Biggl( \log\frac{\int p(x\mid\bar\lambda)\pi_{f,\mathrm{J}}(\bar\lambda) \du\bar\lambda}{\int p(x\mid\bar\lambda)\pi_{\mathrm{J}}(\bar\lambda) \du\bar\lambda}
\, \Bigg| \, \lambda \Biggr)
-\text{E} \Biggl( \log\frac{\int p(x,y\mid\bar\lambda)\pi_{f,\mathrm{J}}(\bar\lambda) \du\bar\lambda}{\int p(x,y\mid\bar\lambda)\pi_{\mathrm{J}}(\bar\lambda) \du\bar\lambda}
\, \Bigg| \, \lambda \Biggr) \notag\\
&=\text{E} \biggl( \log\int f(\bar\theta_1,\dots,\bar\theta_d)\prod_{i=1}^d
\frac{r^{x_i+1/2}\bar\lambda_i^{x_i-1/2}
\exp(-r\bar\lambda_i)}{\Gamma(x_i+1/2)}\du\bar\lambda
\, \bigg| \, \lambda \biggr) \notag\\
& \quad -\text{E} \Biggl(\log\int f(\bar\theta_1,\dots,\bar\theta_d)\prod\limits_{i=1}^d\frac{(r+s)^{x_i+y_i+1/2}\bar\lambda_i^{x_i+y_i-1/2}\exp\{-(r+s)\bar\lambda_i\}}{\Gamma(x_i+y_i+1/2)}\du\bar\lambda \, \Bigg| \, \lambda \Biggr) \notag\\
& = \text{E} \biggl( \log(F(Z,r)) \, \bigg| \, Z_i\sim\text{Po}(r\lambda_i),i=1,\dots,d \biggr)
-\text{E} \biggl( \log(F(Z,r+s)) \, \bigg| \, Z_i\sim\text{Po}((r+s)\lambda_i),i=1,\dots,d \biggr).
\label{negativediff}
\end{align}
From Lemma~\ref{Lemma 1.1}, we know that the risk difference \eqref{negativediff} is finite.

\vspace{0.2cm}

\noindent
\textit{\textbf{Step 2.}} The risk difference \eqref{negativediff} is nonpositive if
\begin{align}
\text{E} \biggl( \log F(Z,r) \, \bigg| \, Z_i\sim\text{Po}(r\lambda_i),i=1,\dots,d \biggr) =\sum_z\log(F(z,r))\Biggl\{
\prod_{i=1}^d \frac{(r\lambda_i)^{z_i}\exp(-r\lambda_i)}{z_i!}
\Biggr\}
\label{increasing}
\end{align}
is a nondecreasing function of $r$.

If we exchange the integration and differentiation in $\frac{\partial F}{\partial r}(z,r)$, we have
 $$\frac{\partial F}{\partial r}(z,r)=\int f(\bar\theta_1,\dots,\bar\theta_d) \prod_{i=1}^d
 \frac{r^{z_i+1/2}\bar\lambda_i^{z_i-1/2}\exp(-r\bar\lambda_i)}{\Gamma(z_i+1/2)}
 \Bigl( \sum_{i=1}^d \frac{z_i+1/2}{r} - \sum_{i=1}^d \bar\lambda_i \Bigr) \du\bar\lambda.$$
 From Lemma~\ref{Lemma 1.2}, we can exchange the integration and differentiation in $\frac{\partial F}{\partial r}(z,r)$.
 
If we differentiate \eqref{increasing} item-by-item, the partial differential function of \eqref{increasing} with respect to $r$ is 
\begin{align}
&\sum_z \log(F(z,r))
\Biggl\{
\prod_{i=1}^d \frac{(r\lambda_i)^{z_i}\exp(-r\lambda_i)}{z_i!}
\Biggr\}
\Bigl(\sum_{i=1}^d \frac{z_i}{r} - \sum_{i=1}^d \lambda_i \Bigr) \notag \\
&\quad + \sum_z 
 \Biggl\{
 \frac{ 
 \int f(\bar\theta_1,\dots,\bar\theta_d) \prod_{i=1}^d
 \frac{r^{z_i+1/2}\bar\lambda_i^{z_i-1/2}\exp(-r\bar\lambda_i)}{\Gamma(z_i+1/2)}
 \Bigl( \sum_{i=1}^d \frac{z_i+1/2}{r} - \sum_{i=1}^d \bar\lambda_i \Bigr) \du\bar\lambda}
 {F(z,r)}
 \Biggr\}
 \prod_{i=1}^d \frac{(r\lambda_i)^{z_i} \exp(-r\lambda_i)}{z_i!}.
 \label{differential}
\end{align}
From Lemma~\ref{Lemma 1.3}, we can differentiate \eqref{increasing} by terms under the condition \eqref{Fineq}.
%

We note that
\begin{align}
&\sum_z \Bigl\{ \log F(z_1,z_2,\dots,z_d,r) \Bigr\}
\Biggl\{ \prod_{i=1}^d\frac{(r\lambda_i)^{z_i} \exp(-r\lambda_i)}{z_i!} \Biggr\}
\sum_{i=1}^d \lambda_i
\notag \\
&= \sum_{i=1}^d \sum_z \biggl[ \frac{z_i}{r} \Bigl\{ \log F(z-\delta_i,r) \Bigr\}
\prod\limits_{j=1}^d\frac{(r\lambda_j)^{z_j}\exp(-r\lambda_j)}{z_j!} \biggr],
\label{transform1}
\end{align}
where $\delta_{ij}$ is defined as 1 if $i = j$ and 0 if $i\neq j$, $\delta_i$ is defined as the $d$-dimensional vector whose $i^{\text{th}}$ element is 1 and all other elements are 0, and $F(z-\delta_i,r)$ is defined as $1$ if $z_i=0$.
Furthermore, we note that
\begin{align}
  \int f(\bar\theta_1,\dots,\bar\theta_d)\prod_{j=1}^d\frac{r^{z_j+1/2}\bar\lambda_j^{z_j-1/2}
  \exp(-r\bar\lambda_j)}{\Gamma(z_j+1/2)}
  \bar\lambda_i \du\bar\lambda  = F(z+\delta_i,r)\frac{z_i+1/2}{r}.  \label{transform2}
\end{align}

\par Thus, from \eqref{transform1} and \eqref{transform2}, the partial differential function \eqref{differential}
of \eqref{increasing} with respect to $r$ is
\begin{align}
&\sum_z\sum_{i=1}^d \biggl[ \frac{z_i}{r} \biggl\{ \log F(z,r) \biggr\}
\prod\limits_{j=1}^d\frac{(r\lambda_j)^{z_j}\exp(-r\lambda_j)}{z_j!} \biggr]
-\sum_z \sum_{i=1}^d \biggl[ \frac{z_i}{r} \biggl\{ \log F(z-\delta_i,r) \biggr\}
\prod_{j=1}^d \frac{(r\lambda_j)^{z_j} \exp(-r\lambda_j)}{z_j!} \biggr] \notag\\
&~~~+ \sum_z \Biggl\{
\frac{F(z,r) \sum_{i=1}^d (z_i+1/2)/r
-\sum_{i=1}^d F(z+\delta_i,r) (z_i+1/2)/r}{F(z,r)} \Biggr\}
\prod_{j=1}^d\frac{(r\lambda_j)^{z_j}\exp(-r\lambda_j)}{z_j!} \notag\\
&= \sum_z
\Biggl[ \sum_{i=1}^d \frac{z_i}{r} \biggl\{ \log 
\frac{F(z,r)}{F(z-\delta_i,r)}\biggr\}
\prod\limits_{j=1}^d\frac{(r\lambda_j)^{z_j}\exp(-r\lambda_j)}{z_j!}
\Biggr]  \notag\\
&~~~+\sum_z \Biggl[ \sum_{i=1}^d
\frac{z_i+1/2}{r}
\biggl\{ 1-\frac{F(z+\delta_i,r)}{F(z,r)} \biggr\}
\prod_{j=1}^d
\frac{(r\lambda_j)^{z_j}\exp(-r\lambda_j)}{z_j!} \Biggr].
\label{differential2}
\end{align}
By the inequality $\log \epsilon \geq 1-\frac{1}{\epsilon}$ for all $\epsilon > 0$, with equality if and only if $\epsilon = 1$, \eqref{differential2} is larger than or equal to
\begin{align}
\sum_z \Biggl[ \sum_{i=1}^d \frac{z_i}{r} \biggl\{ 1-\frac{F(z-\delta_i,r)}{F(z,r)} \biggr\} +\sum_{i=1}^d\frac{z_i+1/2}{r} \biggl\{ 1-\frac{F(z+\delta_i,r)}{F(z,r)} \biggr\}\Biggr]
\prod\limits_{j=1}^d \frac{(r\lambda_j)^{z_j}\exp(-r\lambda_j)}{z_j!}. \label{Fpositive}
\end{align}

From \eqref{Fineq}, we know that \eqref{Fpositive} is nonnegative. Thus, \eqref{differential2} is nonnegative, and \eqref{increasing} is a nondecreasing function. This proves that $p_{f,\mathrm{J}}(y\mid x)$ weakly dominates $p_{\mathrm{J}}(y\mid x)$.

Furthermore, if $F(z,r)$ is not a constant function of $z$ on $\mathbb{N}^d$, then $F(z-\delta_i,r) \not\equiv F(z,r)$ for some $i$. Therefore, the inequality $\log \epsilon \geq 1-\frac{1}{\epsilon}$ is strict for some terms in the sum, and \eqref{differential2} is strictly positive. Thus, \eqref{increasing} is strictly increasing, and $p_{f,\mathrm{J}}(y\mid x)$ dominates $p_{\mathrm{J}}(y\mid x)$.

\vspace{0.3cm}

\noindent
\textit{\textbf{Proof of part 2.}} We prove that \eqref{Fineq} is satisfied if $f$ satisfies the conditions of the second half of Theorem~\ref{Theorem 1}. Let $\theta_j \coloneqq  \sqrt{\lambda_j}$ $(j=1,\dots,d)$. We show this in three steps. \eqref{Fineq} is obtained by combining Steps 2 and 3.
In Step 1, through integration by parts on $\theta_i$, we prove that 
\begin{align}
F(z+\delta_i,r)-F(z,r) = 2^{d-1}\int  \frac{\partial f}{\partial\theta_i}(\theta)  \prod\limits_{j=1}^d
\frac{r^{z_j+1/2}\theta_j^{2z_j+\delta_{ji}}\exp(-r\theta_j^2)}
{\Gamma(z_j+\delta_{ji}+1/2)}\du\theta.
    \label{F-1st-diff}
\end{align}
In Step 2, by performing integration by parts on each $\theta_i$ again and using the condition \eqref{t1-condition-d} of the second derivative, we prove that
\begin{align}
&\sum_{i=1}^d z_i \{ F(z,r)-F(z-\delta_i,r) \}
 +\sum_{i=1}^d (z_i+1/2) \{F(z,r)-F(z+\delta_i,r)\} \notag\\
&\ge \sum_{i=1}^d \frac{2^{d-2}}{r}\mathop{\lim}_{a\to0\atop b\to\infty}\mathop{\lim}_{u\to0\atop v\to\infty} \Biggl\{ \Biggl[\int_{[a,b]^{d-1}} \frac{\partial f}{\partial\theta_i}(\theta)
\biggl\{\prod_{j\neq i}\frac{r^{z_j+1/2}\theta_j^{2z_j}\exp(-r\theta_j^2)}{\Gamma(z_j+1/2)} \du\theta_j\biggr\} \frac{r^{z_i+1/2}\theta_i^{2z_i}\exp(-r\theta_i^2)}{\Gamma(z_i+1/2)} \Biggr]^{v}_{\theta_i = u}\Biggr\}. \label{thm1-step2}
\end{align}
In Step 3, using the condition \eqref{t1-condition-e} of the derivative on the boundary, we show that \eqref{thm1-step2} $\ge 0$.
The details of each step are presented below.

\noindent
\textit{\textbf{Step 1.}} Using the substitution $\theta=\sqrt{\lambda}$ and the definition of the function $F$, we obtain
\begin{align*}
&F(z+\delta_i,r)-F(z,r) \notag\\
&=-2^{d-1}\int\biggl\{
\int f(\theta)\prod\limits_{j\neq i}\frac{r^{z_j+1/2}\theta_j^{2z_j}\exp(-r\theta_j^2)}
{\Gamma(z_j+1/2)} \du \theta_j 
\biggr\}
\frac{r^{z_i+1/2}\theta_i^{2z_i+1}\frac{\partial}{\partial\theta_i}\exp(-r\theta_i^2)}
{\Gamma(z_i+3/2)}\du\theta_i \notag\\
& ~~ -2^{d-1}\int\biggl\{
\int f(\theta)\prod\limits_{j\neq i}\frac{r^{z_j+1/2}\theta_j^{2z_j}\exp(-r\theta_j^2)}
{\Gamma(z_j+1/2)} \du \theta_j 
\biggr\}
\frac{r^{z_i+1/2}(\frac{\partial}{\partial\theta_i}\theta_i^{2z_i+1})\exp(-r\theta_i^2)}{\Gamma(z_i+3/2)}\du\theta_i.
\end{align*}
By performing integration by parts on $\theta_i$, we obtain
\begin{align}
&F(z+\delta_i,r)-F(z,r) \notag\\
&=\mathop{\lim}_{a\to0\atop b\to\infty}\mathop{\lim}_{u\to0\atop v\to\infty} \biggl[2^{d-1} \int_{[u,v]\times[a,b]^{d-1}} \biggl\{ \frac{\partial}{\partial\theta_i}f(\theta) \biggr\}
\prod\limits_{j=1}^d\frac{r^{z_j+1/2}
\theta_j^{2z_j+\delta_{ji}}\exp(-r\theta_j^2)}{\Gamma(z_j+\delta_{ji}+1/2)} \du\theta \notag\\
& ~~~~ -2^{d-1} \Biggl[ \biggl\{
\int_{[a,b]^{d-1}} f(\theta)\prod\limits_{j\neq i}\frac{r^{z_j+1/2}\theta_j^{2z_j}\exp(-r\theta_j^2)}
{\Gamma(z_j+1/2)} \du \theta_j 
\biggr\}
\frac{r^{z_i+1/2}\theta_i^{2z_i+1}\exp(-r\theta_i^2)}
{\Gamma(z_i+3/2)} \Biggr]^{v}_{\theta_i = u}\biggl]. \label{F-diff-1}
\end{align}
Here, we use auxiliary variables $a,b,u,v$ and Lemma~\ref{Lemma 1.4} to ensure that the above equation \eqref{F-diff-1} holds.

\par Because of Lemma~\ref{Lemma 1.5},
\begin{equation}
\mathop{\lim}_{u\to0\atop v\to\infty} -2^{d-1} \Biggl[ \biggl\{
\int_{[a,b]^{d-1}} f(\theta)\prod\limits_{j\neq i}\frac{r^{z_j+1/2}\theta_j^{2z_j}\exp(-r\theta_j^2)}
{\Gamma(z_j+1/2)} \du \theta_j 
\biggr\}
\frac{r^{z_i+1/2}\theta_i^{2z_i+1}\exp(-r\theta_i^2)}
{\Gamma(z_i+3/2)} \Biggr]^{v}_{\theta_i = u}=0. \label{1st-partial-3}
\end{equation}

Thus, using \eqref{F-diff-1}, Lemma~\ref{Lemma 1.4}, and \eqref{1st-partial-3},  we obtain \eqref{F-1st-diff}.

\vspace{0.2cm}

\noindent
\textit{\textbf{Step 2.}} Using \eqref{F-1st-diff} and auxiliary variables $a,b,u,v$, we obtain
\begin{align}
&\sum_{i=1}^d  z_i \{ F(z,r)-F(z-\delta_i,r) \}
 +\sum_{i=1}^d (z_i+1/2) \{F(z,r)-F(z+\delta_i,r)\} \notag\\
&=\sum_{i=1}^d 2^{d-1} \int\biggl[  \frac{\partial f }{\partial\theta_i}(\theta)  z_i
\prod_{j=1}^d\frac{r^{(z_j-\delta_{ji})+1/2}\theta_j^{2(z_j-\delta_{ji})+\delta_{ji}}\exp(-r\theta_j^2)}
{\Gamma(z_j+1/2)} \notag\\
& \qquad-  \frac{\partial f }{\partial\theta_i} (\theta)  (z_i+1/2)
\prod_{j=1}^d
\frac{r^{z_j+1/2}\theta_j^{2z_j+\delta_{ji}}\exp(-r\theta_j^2)}{\Gamma(z_j+\delta_{ji}+1/2)} \biggr]\du\theta \notag\\
&=\sum_{i=1}^d \frac{2^{d-2}}{r} \mathop{\lim}_{a\to0\atop b\to\infty}\mathop{\lim}_{u\to0\atop v\to\infty} \Biggl\{ \int_u^v\int_{[a,b]^{d-1}}  \frac{\partial f}{\partial\theta_i}(\theta)
\biggl\{\prod_{j\neq i}\frac{r^{z_j+1/2}\theta_j^{2z_j}\exp(-r\theta_j^2)}{\Gamma(z_j+1/2)} \du\theta_j\biggr\}
\frac{r^{z_i+1/2}\frac{\partial\theta_i^{2z_i}}{\partial\theta_i}\exp(-r\theta_i^2)}{\Gamma(z_i+1/2)}\du\theta_i \notag\\
& \qquad+ \int_u^v\int_{[a,b]^{d-1}}\frac{\partial f}{\partial\theta_i}(\theta)
\biggl\{\prod_{j\neq i}\frac{r^{z_j+1/2}\theta_j^{2z_j}\exp(-r\theta_j^2)}{\Gamma(z_j+1/2)}\du\theta_j\biggr\}
\frac{r^{z_i+1/2}\theta_i^{2z_i}\frac{\partial}{\partial\theta_i}\exp(-r\theta_i^2)}
{\Gamma(z_i+1/2)}\du\theta_i \Biggr\}. \label{F-diff-2-1}
\end{align}

By performing integration by parts on each parameter again, we find that \eqref{F-diff-2-1} equals
\begin{align}
&\sum_{i=1}^d \frac{2^{d-2}}{r}\mathop{\lim}_{a\to0\atop b\to\infty}\mathop{\lim}_{u\to0\atop v\to\infty} \Biggl\{
-\int_u^v\int_{[a,b]^{d-1}}  \frac{\partial}{\partial\theta_i}\Bigl\{\frac{\partial f}{\partial\theta_i}(\theta)\Bigr\}
\prod_{j\neq i}\frac{r^{z_j+1/2}\theta_j^{2z_j}\exp(-r\theta_j^2)}{\Gamma(z_j+1/2)}\du\theta_j \frac{r^{z_i+1/2}\theta_i^{2z_i}\exp(-r\theta_i^2)}{\Gamma(z_i+1/2)}\du\theta_i \notag \\
&\qquad+ \biggl[\int_{[a,b]^{d-1}}  \frac{\partial f}{\partial\theta_i}(\theta) 
\biggl\{\prod_{j\neq i}\frac{r^{z_j+1/2}\theta_j^{2z_j}\exp(-r\theta_j^2)}{\Gamma(z_j+1/2)}\du\theta_j\biggr\} \frac{r^{z_i+1/2}\theta_i^{2z_i}\exp(-r\theta_i^2)}{\Gamma(z_i+1/2)} \biggr]^{v}_{\theta_i = u} \Biggr\}. \label{F-diff-2}
\end{align}

From \eqref{t1-condition-c}, we have
\begin{equation}
\int\biggl\lvert \frac{\partial}{\partial\theta_i}
 \Bigl\{ \frac{\partial f}{\partial\theta_i}(\theta) 
 \Bigr\} \biggr\rvert
\prod\limits_{j=1}^d\frac{r^{z_j+1/2}\theta_j^{2z_j} \exp(-r\theta_j^2)}{\Gamma(z_j+1/2)}\du\theta <\infty,\ \forall i. \label{2nd-diff-finite}
\end{equation}

From \eqref{2nd-diff-finite} and condition \eqref{t1-condition-d}, we have
\begin{align}
&\sum_{i=1}^d \mathop{\lim}_{a\to0\atop b\to\infty}\mathop{\lim}_{u\to0\atop v\to\infty} \Biggl\{
-\int_u^v\int_{[a,b]^{d-1}} \frac{\partial^2 f}{\partial\theta_i^2}(\theta)
\prod_{j\neq i}\frac{r^{z_j+1/2}\theta_j^{2z_j}\exp(-r\theta_j^2)}{\Gamma(z_j+1/2)}\du\theta_j \frac{r^{z_i+1/2}\theta_i^{2z_i}\exp(-r\theta_i^2)}{\Gamma(z_i+1/2)}\du\theta_i \Biggr\} \notag\\
&=- \int  \sum_{i=1}^d \frac{\partial^2 f}{\partial\theta_i^2}(\theta)
\prod_{j=1}^d \frac{r^{z_j+1/2}\theta_j^{2z_j} \exp(-r\theta_j^2)}{\Gamma(z_j+1/2)}\du\theta \ge 0. \label{F-2nd-diff-1}
\end{align}

Using \eqref{F-diff-2-1}, \eqref{F-diff-2}, and \eqref{F-2nd-diff-1}, we obtain the inequality \eqref{thm1-step2}.

\vspace{0.2cm}

\noindent
\textit{\textbf{Step 3.}}
From Lemma~\ref{Lemma 1.6}, we know that \eqref{thm1-step2} $\ge 0$, which completes the proof.
\qed

\subsection{Detailed Proof of Theorem 2}
Theorem 2 in the main paper is copied as follows.

\begin{theorem}\label{Theorem 2}\leavevmode
\begin{enumerate}[label=\arabic*)]
\item
Let $\theta_i\coloneqq \sqrt{\lambda_i/\gamma_i},\ i=1,\dots,d$. The Bayesian predictive distribution $p_{f,\mathrm{J}}(y\mid x)$ weakly dominates the Bayesian predictive distribution $p_{\mathrm{J}}(y\mid x)$ if for every $z\in \mathbb{N}^d$, $r\in\mathbb{R}_+^d$,
\begin{equation}
    \sum\limits_{i=1}^d\gamma_ir_iz_i \Bigl\{ F(z,r)-F(z-\delta_i,r)\Bigr\} +\sum\limits_{i=1}^d\gamma_ir_i(z_i+1/2) \Bigl\{F(z,r)-F(z+\delta_i,r) \Bigr\} \ge 0.
\label{Fineq3}
\end{equation} 
where $F(z-\delta_i,r)\coloneqq 1$ if $z_i=0$.
Furthermore, if $F(z,r)$ is not a constant function of $z$ on $\mathbb{N}^d$, then $p_{f,\mathrm{J}}(y\mid x)$ dominates $p_{\mathrm{J}}(y\mid x)$.

\item
Condition \eqref{Fineq3} is satisfied if $f\in\mathbf{C}^2([0,\infty)^d)$ satisfies regularity conditions \eqref{t1-condition-a}, \eqref{t1-condition-b}, and \eqref{t1-condition-c}, and conditions \eqref{t1-condition-d} and \eqref{t1-condition-e}.
\end{enumerate}
\end{theorem}

\vspace{0.2cm}
\noindent
\textit{\textbf{Proof of part 1.}} For every $i$ and $\tau\in[0,1]$, let $\frac{1}{t_i(\tau)}\coloneqq \frac{1}{r_i}(1-\tau)+\frac{1}{r_i+s_i}\tau$. Then, $t_i(\tau)$ is a smooth monotonically increasing function of $\tau\in[0, 1]$ satisfying $t_i(0) = r_i$, $t_i(1) = r_i + s_i$, and $\dot t_i/t_i=\gamma_it_i$. Let $\text{Po}(\lambda)$ denote the Poisson distribution with parameter $\lambda$.

We prove the first half of Theorem~\ref{Theorem 2} in two steps. In Step 1, we show that the difference between the K-L risks of Bayesian predictive distributions based on $\pi_{f,\mathrm{J}}$ and $\pi_{\mathrm{J}}$ is $$\text{E} \biggl( \log F(Z,t(0))  \, \bigg| \, Z_i\sim\text{Po}(t_i(0)\lambda_i),i=1,\dots,d \biggr)
-\text{E} \biggl( \log F(Z,t(1)) \, \bigg| \, Z_i\sim\text{Po}(t_i(1)\lambda_i),i=1,\dots,d \biggr).$$
In Step 2, we prove that the partial derivative of $\text{E} \biggl( \log F(Z,t(\tau)) \, \bigg| \, Z_i\sim\text{Po}(t_i(\tau)\lambda_i),i=1,\dots,d \biggr)$ with respect to $\tau$ is nonnegative, and is positive when $F(z,r)$ is not a constant function in $z$. The details of each step are presented below.

\noindent
\textit{\textbf{Step 1.}} Let $\bar\theta_j\coloneqq \sqrt{\frac{\bar\lambda_j}{\gamma_j}},\ j=1,\dots,d$. The difference between the K-L risks of Bayesian predictive distributions based on $\pi_{f,\mathrm{J}}$ and $\pi_{\mathrm{J}}$ is
\begin{align}
&\text{E} \biggl(
\log\frac{p_{\mathrm{J}}(y\mid x)}{p_{f,\mathrm{J}}(y\mid x)} 
\, \bigg| \, \lambda
\biggr) \notag\\
&=\text{E}
\Biggl( \log\frac{\int p(x\mid\bar\lambda)\pi_{f,\mathrm{J}}(\bar\lambda) \du\bar\lambda}{\int p(x\mid\bar\lambda)\pi_{\mathrm{J}}(\bar\lambda) \du\bar\lambda}
\, \Bigg| \, \lambda \Biggr)
-\text{E} \Biggl( \log\frac{\int p(x,y\mid\bar\lambda)\pi_{f,\mathrm{J}}(\bar\lambda) \du\bar\lambda}{\int p(x,y\mid\bar\lambda)\pi_{\mathrm{J}}(\bar\lambda) \du\bar\lambda}
\, \Bigg| \, \lambda \Biggr) \notag\\
&=\text{E} \biggl( \log\int f(\bar\theta_1,\dots,\bar\theta_d)\prod_{i=1}^d
\frac{r_i^{x_i+1/2}\bar\lambda_i^{x_i-1/2}
\exp(-r_i\bar\lambda_i)}{\Gamma(x_i+1/2)}\du\bar\lambda
\, \bigg| \, \lambda \biggr) \notag\\
& \quad -\text{E} \Biggl(\log\int f(\bar\theta_1,\dots,\bar\theta_d)\prod\limits_{i=1}^d\frac{(r_i+s_i)^{x_i+y_i+1/2}\bar\lambda_i^{x_i+y_i-1/2}\exp\{-(r_i+s_i)\bar\lambda_i\}}{\Gamma(x_i+y_i+1/2)}\du\bar\lambda \, \Bigg| \, \lambda \Biggr) \notag\\
& = \text{E} \biggl( \log F(Z,t(0))  \, \bigg| \, Z_i\sim\text{Po}(t_i(0)\lambda_i),i=1,\dots,d \biggr)
-\text{E} \biggl( \log F(Z,t(1)) \, \bigg| \, Z_i\sim\text{Po}(t_i(1)\lambda_i),i=1,\dots,d \biggr).
\label{negativediff-2}
\end{align}

From Lemma~\ref{Lemma 2.1}, we know that the risk difference \eqref{negativediff-2} is finite.

\vspace{0.2cm}

\noindent
\textit{\textbf{Step 2.}} The risk difference \eqref{negativediff-2} is nonpositive if
\begin{align}
\text{E} \biggl( \log F(Z,t(\tau)) \, \bigg| \, Z_i\sim\text{Po}(t_i(\tau)\lambda_i),i=1,\dots,d \biggr) =\sum_z\Bigl\{ \log F(z,t(\tau)) \Bigr\}
\prod_{i=1}^d \frac{(t_i(\tau)\lambda_i)^{z_i}\exp(-t_i(\tau)\lambda_i)}{z_i!}
\label{increasing-t2}
\end{align}
is a nondecreasing function of $\tau$.

Using $\dot t_i/t_i=\gamma_it_i$, we have
\begin{equation}
 \frac{\partial F}{\partial \tau}(z,t(\tau))=\int f(\bar\theta_1,\dots,\bar\theta_d) \prod_{i=1}^d
 \frac{t_i(\tau)^{z_i+1/2}\bar\lambda_i^{z_i-1/2}\exp(-t_i(\tau)\bar\lambda_i)}{\Gamma(z_i+1/2)}
 \biggl\{ \sum_{j=1}^d (z_j+1/2)\gamma_jt_j(\tau) - \sum_{k=1}^d \bar\lambda_k\gamma_kt_k^2(\tau) \biggr\} \du\bar\lambda.
 \label{Fdiff-t2}
\end{equation}
 From Lemma~\ref{Lemma 2.2}, we can exchange the integration and differentiation in $\frac{\partial F}{\partial \tau}(z,t(\tau))$.
 
From $\dot t_i/t_i=\gamma_it_i$ and \eqref{Fdiff-t2}, the partial differential function of \eqref{increasing-t2} with respect to $\tau$ is 
\begin{align}
&\sum_z \Bigl\{ \log F(z,t(\tau)) \Bigr\} 
\Biggl\{
\prod_{i=1}^d \frac{(t_i(\tau)\lambda_i)^{z_i}\exp(-t_i(\tau)\lambda_i)}{z_i!}
\Biggr\}
\Bigl(\sum_{j=1}^d z_j\gamma_jt_j(\tau) - \sum_{k=1}^d \lambda_k\gamma_kt_k^2(\tau) \Bigr) \notag \notag\\
&\quad + \sum_z 
 \Biggl[
 \frac{ 
 \int f(\bar\theta_1,\dots,\bar\theta_d) \prod_{i=1}^d
 \frac{t_i(\tau)^{z_i+1/2}\bar\lambda_i^{z_i-1/2}\exp(-t_i(\tau)\bar\lambda_i)}{\Gamma(z_i+1/2)}
 \Bigl\{ \sum_{j=1}^d (z_j+1/2)\gamma_jt_j(\tau) - \sum_{k=1}^d \bar\lambda_k\gamma_kt_k^2(\tau) \Bigr\} \du\bar\lambda}
 { F(z,t(\tau))}
 \Biggr] \notag\\
 & ~~~~~
 \Biggl\{
\prod_{l=1}^d \frac{(t_l(\tau)\lambda_l)^{z_l}\exp(-t_l(\tau)\lambda_l)}{z_l!}
\Biggr\}.
 \label{differential-t2}
\end{align}
From Lemma~\ref{Lemma 2.3}, we can differentiate \eqref{increasing-t2} by terms under condition \eqref{Fineq3}.

We note that
\begin{align}
&\sum_z \biggl\{ \log F(z,t(\tau)) \biggr\}
\Biggl\{ \prod_{i=1}^d\frac{(t_i(\tau)\lambda_i)^{z_i} \exp(-t_i(\tau)\lambda_i)}{z_i!} \Biggr\}
\sum_{j=1}^d \lambda_j\gamma_jt_j^2(\tau)
\notag \\
&= \sum_{i=1}^d \sum_z \biggl[ z_i\gamma_it_i(\tau) \Bigl\{ \log F(z-\delta_i,t(\tau)) \Bigr\}
\prod\limits_{j=1}^d\frac{(t_j(\tau)\lambda_j)^{z_j}\exp(-t_j(\tau)\lambda_j)}{z_j!} \biggr],
\label{transform1-t2}
\end{align}
where $\delta_{ij}$ is defined as 1 if $i = j$ and 0 if $i\neq j$, $\delta_i$ is defined as the $d$-dimensional vector whose $i^{\text{th}}$ element is 1 and all other elements are 0, and $F(z-\delta_i,t(\tau))$ is defined as $1$ if $z_i=0$.
Furthermore, we note that
\begin{align}
  \int f(\bar\theta_1,\dots,\bar\theta_d)
\prod_{j=1}^d\frac{t_j(\tau)^{z_j+1/2}\bar\lambda_j^{z_j-1/2}
  \exp(-t_j(\tau)\bar\lambda_j)}{\Gamma(z_j+1/2)}
  \bar\lambda_i\gamma_it_i^2(\tau) \du\bar\lambda  = F(z+\delta_i,t(\tau))(z_i+1/2)\gamma_it_i(\tau). \label{transform2-t2}
\end{align}

\par Thus, from \eqref{transform1-t2} and \eqref{transform2-t2}, the partial differential function \eqref{differential-t2}
of \eqref{increasing-t2} with respect to $\tau$ is
\begin{align}
&\sum_z\sum_{i=1}^d \biggl[ z_i\gamma_it_i(\tau) \bigl\{ \log F(z,t(\tau)) \bigr\}
\prod\limits_{j=1}^d\frac{(t_j(\tau)\lambda_j)^{z_j}\exp(-t_j(\tau)\lambda_j)}{z_j!} \biggr] \notag\\
&~~~-\sum_z \sum_{i=1}^d \biggl[ z_i\gamma_it_i(\tau) \bigl\{ \log F(z-\delta_i,t(\tau)) \bigr\}
\prod_{j=1}^d \frac{(t_j(\tau)\lambda_j)^{z_j} \exp(-t_j(\tau)\lambda_j)}{z_j!} \biggr] \notag\\
&~~~+ \sum_z \Biggl\{
\frac{F(z,t(\tau)) \sum_{i=1}^d (z_i+1/2)\gamma_it_i(\tau)
-\sum_{i=1}^d F(z+\delta_i,t(\tau)) (z_i+1/2)\gamma_it_i(\tau)}{F(z,t(\tau))} \Biggr\}  \biggl\{ \prod_{j=1}^d\frac{(t_j(\tau)\lambda_j)^{z_j}\exp(-t_j(\tau)\lambda_j)}{z_j!}\biggr\} \notag\\
&= \sum_z
\Biggl[ \sum_{i=1}^d z_i\gamma_it_i(\tau) \biggl\{ \log 
\frac{F(z,t(\tau))}{F(z-\delta_i,t(\tau))}\biggr\}
\prod\limits_{j=1}^d\frac{(t_j(\tau)\lambda_j)^{z_j}\exp(-t_j(\tau)\lambda_j)}{z_j!}
\Biggr] \notag \\
&~~~+\sum_z \Biggl[ \sum_{i=1}^d
(z_i+1/2)\gamma_it_i(\tau)
\biggl\{ 1-\frac{F(z+\delta_i,t(\tau))}{F(z,t(\tau))} \biggr\}
\prod_{j=1}^d
\frac{(t_j(\tau)\lambda_j)^{z_j}\exp(-t_j(\tau)\lambda_j)}{z_j!} \Biggr].
\label{differential2-t2}
\end{align}
By the inequality $\log \epsilon \geq 1-\frac{1}{\epsilon}$ for all $\epsilon > 0$, with equality if and only if $\epsilon = 1$, \eqref{differential2-t2} is larger than or equal to
\begin{align}
&\sum_z \Biggl[ \sum_{i=1}^d z_i\gamma_it_i(\tau) \biggl\{ 1-\frac{F(z-\delta_i,t(\tau))}{F(z,t(\tau))} \biggr\}  \notag \\
&~~~+\sum_{i=1}^d(z_i+1/2)\gamma_it_i(\tau) \biggl\{ 1-\frac{F(z+\delta_i,t(\tau))}{F(z,t(\tau))} \biggr\}\Biggr]  \biggl\{ \prod_{j=1}^d\frac{(t_j(\tau)\lambda_j)^{z_j}\exp(-t_j(\tau)\lambda_j)}{z_j!}\biggr\}. \label{Fpositive-t2}
\end{align}

From \eqref{Fineq3}, we know that \eqref{Fpositive-t2} is nonnegative. Thus, \eqref{differential2-t2} is nonnegative, and \eqref{increasing-t2} is a nondecreasing function. This proves that $p_{f,\mathrm{J}}(y\mid x)$ weakly dominates $p_{\mathrm{J}}(y\mid x)$.

Furthermore, if $F(z,r)$ is not a constant function of $z$ on $\mathbb{N}^d$, then $F(z-\delta_i,r) \not\equiv F(z,r)$ for some $i$. Therefore, the inequality $\log \epsilon \geq 1-\frac{1}{\epsilon}$ is strict for some terms in the sum, and \eqref{differential2-t2} is strictly positive when $\tau=0$. Thus, \eqref{increasing-t2} is strictly increasing at $\tau=0$, and $p_{f,\mathrm{J}}(y\mid x)$ dominates $p_{\mathrm{J}}(y\mid x)$.

\vspace{0.3cm}

\noindent
\textit{\textbf{Proof of part 2.}} Next, we prove that if $f$ satisfies the conditions of the second half of Theorem~\ref{Theorem 2}, \eqref{Fineq3} is satisfied.
Let $\theta_j\coloneqq\sqrt{\frac{\lambda_j}{\gamma_j}},\ j=1,\dots,d$. We show this in three steps. \eqref{Fineq3} is obtained by combining Steps 2 and 3.
In Step 1, through integration by parts on $\theta_i$, we prove that 
\begin{align}
F(z+\delta_i,r)-F(z,r)= 2^{d-1}\int  \frac{\partial f}{\partial\theta_i}(\theta) \prod\limits_{j=1}^d
\frac{(\gamma_jr_j)^{z_j+1/2}\theta_j^{2z_j+\delta_{ji}}\exp(-\gamma_jr_j\theta_j^2)}
{\Gamma(z_j+\delta_{ji}+1/2)}\du\theta.
    \label{F-1st-diff-t2}
\end{align}
In Step 2, by performing integration by parts on each $\theta_i$ again and using the condition \eqref{t1-condition-d} of the second derivative, we prove that
\begin{align}
    &\sum_{i=1}^d  \gamma_ir_iz_i \{ F(z,r)-F(z-\delta_i,r) \}
 +\sum_{i=1}^d \gamma_ir_i(z_i+1/2) \{F(z,r)-F(z+\delta_i,r)\} \notag\\
 &\ge \sum_{i=1}^d 2^{d-2}\mathop{\lim}_{a\to0\atop b\to\infty}\mathop{\lim}_{u\to0\atop v\to\infty} \Biggl\{
\biggl[ \int_{[a,b]^{d-1}}  \frac{\partial f}{\partial\theta_i}(\theta)
\biggl\{ \prod_{j\neq i}\frac{(\gamma_jr_j)^{z_j+1/2}\theta_j^{2z_j}\exp(-\gamma_jr_j\theta_j^2)}{\Gamma(z_j+1/2)}\du\theta_j\biggr\} \frac{(\gamma_ir_i)^{z_i+1/2}\theta_i^{2z_i}\exp(-\gamma_ir_i\theta_i^2)}{\Gamma(z_i+1/2)} \biggr]^{v}_{\theta_i = u} \Biggr\}. \label{thm2-step2}
\end{align}
In Step 3, using the condition \eqref{t1-condition-e} of the derivative on the boundary, we show that \eqref{thm2-step2} $\ge 0$.
The details of each step are presented below.

\noindent
\textit{\textbf{Step 1.}} Using the substitution $\lambda_j=\gamma_j\theta_j^2$, $\du\lambda_j=2\theta_j\gamma_j\du\theta_j$ and the definition of the function $F$, we obtain 
\begin{align}
	&F(z+\delta_i,r)-F(z,r) \notag\\
	&=-2^{d-1}\int f(\theta)\prod_{j\neq i}\frac{(\gamma_jr_j)^{z_j+1/2}\theta_j^{2z_j}\exp(-\gamma_jr_j\theta_j^2)}
{\Gamma(z_j+1/2)}\du\theta_j \frac{(\gamma_ir_i)^{z_i+1/2}\theta_i^{2z_i+1}\frac{\partial}{\partial\theta_i}\exp(-\gamma_ir_i\theta_i^2)}
{\Gamma(z_i+3/2)}\du\theta_i \notag\\
	&\qquad-2^{d-1}\int f(\theta)\prod_{j\neq i}\frac{(\gamma_jr_j)^{z_j+1/2}\theta_j^{2z_j}\exp(-\gamma_jr_j\theta_j^2)}
{\Gamma(z_j+1/2)}\du\theta_j \frac{(\gamma_ir_i)^{z_i+1/2}(\frac{\partial}{\partial\theta_i}\theta_i^{2z_i+1})
	\exp(-\gamma_ir_i\theta_i^2)}{\Gamma(z_i+3/2)}\du\theta_i. \notag
\end{align}

By performing integration by parts on $\theta_i$, we obtain
\begin{align}
&F(z+\delta_i,r)-F(z,r) \notag\\
&=\mathop{\lim}_{a\to0\atop b\to\infty}\mathop{\lim}_{u\to0\atop v\to\infty} \Biggl( 2^{d-1} \int_{[u,v]\times[a,b]^{d-1}}  \frac{\partial f }{\partial\theta_i}(\theta) 
\prod\limits_{j=1}^d\frac{(\gamma_jr_j)^{z_j+1/2}
\theta_j^{2z_j+\delta_{ji}}\exp(-\gamma_jr_j\theta_j^2)}{\Gamma(z_j+\delta_{ji}+1/2)} \du\theta \notag\\
& ~~~~ -2^{d-1} \Biggl[ \biggl\{
\int_{[a,b]^{d-1}} f(\theta)\prod\limits_{j\neq i}\frac{(\gamma_jr_j)^{z_j+1/2}\theta_j^{2z_j}\exp(-\gamma_jr_j\theta_j^2)}
{\Gamma(z_j+1/2)} \du \theta_j 
\biggr\}
\frac{(\gamma_ir_i)^{z_i+1/2}\theta_i^{2z_i+1}\exp(-\gamma_ir_i\theta_i^2)}
{\Gamma(z_i+3/2)} \Biggr]^{v}_{\theta_i = u}\Biggl) . \label{F-diff-1-t2}
\end{align}
Here, we use auxiliary variables $a,b,u,v$ and Lemma~\ref{Lemma 2.4} to ensure that the above equation \eqref{F-diff-1-t2} holds.

\par Because of Lemma~\ref{Lemma 2.5},
\begin{equation}
\mathop{\lim}_{u\to0\atop v\to\infty} -2^{d-1} \Biggl[ \biggl\{
\int_{[a,b]^{d-1}} f(\theta)\prod\limits_{j\neq i}\frac{(\gamma_jr_j)^{z_j+1/2}\theta_j^{2z_j}\exp(-\gamma_jr_j\theta_j^2)}
{\Gamma(z_j+1/2)} \du \theta_j 
\biggr\}
\frac{(\gamma_ir_i)^{z_i+1/2}\theta_i^{2z_i+1}\exp(-\gamma_ir_i\theta_i^2)}
{\Gamma(z_i+3/2)} \Biggr]^{v}_{\theta_i = u}=0. \label{1st-partial-3-t2}
\end{equation}

Thus, using \eqref{F-diff-1-t2}, Lemma~\ref{Lemma 2.4}, and \eqref{1st-partial-3-t2}, we obtain \eqref{F-1st-diff-t2}.

\vspace{0.2cm}

\noindent
\textit{\textbf{Step 2.}} Using \eqref{F-1st-diff-t2} and auxiliary variables $a,b,u,v$, we obtain
\begin{align}
&\sum_{i=1}^d  \gamma_ir_iz_i \{ F(z,r)-F(z-\delta_i,r) \}
 +\sum_{i=1}^d \gamma_ir_i(z_i+1/2) \{F(z,r)-F(z+\delta_i,r)\} \notag\\
&=\sum_{i=1}^d 2^{d-1} \int\biggl\{  \frac{\partial f }{\partial\theta_i}(\theta)  \gamma_ir_iz_i
\prod_{j=1}^d\frac{(\gamma_jr_j)^{(z_j-\delta_{ji})+1/2}\theta_j^{2(z_j-\delta_{ji})+\delta_{ji}}\exp(-\gamma_jr_j\theta_j^2)}
{\Gamma(z_j+1/2)} \notag\\
& \qquad-  \frac{\partial f }{\partial\theta_i} (\theta)  \gamma_ir_i(z_i+1/2)
\prod_{j=1}^d
\frac{(\gamma_jr_j)^{z_j+1/2}\theta_j^{2z_j+\delta_{ji}}\exp(-\gamma_jr_j\theta_j^2)}{\Gamma(z_j+\delta_{ji}+1/2)} \biggr\}\du\theta \notag\\
&=\sum_{i=1}^d 2^{d-2} \mathop{\lim}_{a\to0\atop b\to\infty}\mathop{\lim}_{u\to0\atop v\to\infty} \Biggl\{ \int_u^v\int_{[a,b]^{d-1}}  \frac{\partial f}{\partial\theta_i}(\theta)
\biggl\{\prod_{j\neq i}\frac{(\gamma_jr_j)^{z_j+1/2}\theta_j^{2z_j}\exp(-\gamma_jr_j\theta_j^2)}{\Gamma(z_j+1/2)} \du\theta_j\biggr\}
\frac{(\gamma_ir_i)^{z_i+1/2}\frac{\partial\theta_i^{2z_i}}{\partial\theta_i}\exp(-\gamma_ir_i\theta_i^2)}{\Gamma(z_i+1/2)}\du\theta_i \notag\\
& \qquad+ \int_u^v\int_{[a,b]^{d-1}}\frac{\partial f}{\partial\theta_i}(\theta)
\biggl\{\prod_{j\neq i}\frac{(\gamma_jr_j)^{z_j+1/2}\theta_j^{2z_j}\exp(-\gamma_jr_j\theta_j^2)}{\Gamma(z_j+1/2)}\du\theta_j\biggr\}
\frac{(\gamma_ir_i)^{z_i+1/2}\theta_i^{2z_i}\frac{\partial}{\partial\theta_i}\exp(-\gamma_ir_i\theta_i^2)}
{\Gamma(z_i+1/2)}\du\theta_i \Biggr\}. \label{F-diff-2-1-t2}
\end{align}

By performing integration by parts on each parameter again, we find that \eqref{F-diff-2-1-t2} equals
\begin{align}
&\sum_{i=1}^d 2^{d-2}\mathop{\lim}_{a\to0\atop b\to\infty}\mathop{\lim}_{u\to0\atop v\to\infty} \Biggl\{
-\int_u^v\int_{[a,b]^{d-1}}  \frac{\partial}{\partial\theta_i}\Bigl\{\frac{\partial f}{\partial\theta_i}(\theta)\Bigr\}
\prod_{j\neq i}\frac{(\gamma_jr_j)^{z_j+1/2}\theta_j^{2z_j}\exp(-\gamma_jr_j\theta_j^2)}{\Gamma(z_j+1/2)}\du\theta_j \frac{(\gamma_ir_i)^{z_i+1/2}\theta_i^{2z_i}\exp(-\gamma_ir_i\theta_i^2)}{\Gamma(z_i+1/2)}\du\theta_i \notag \\
&\qquad+ \biggl[ \int_{[a,b]^{d-1}}  \frac{\partial f}{\partial\theta_i}(\theta) 
\biggl\{\prod_{j\neq i}\frac{(\gamma_jr_j)^{z_j+1/2}\theta_j^{2z_j}\exp(-\gamma_jr_j\theta_j^2)}{\Gamma(z_j+1/2)}\du\theta_j\biggr\} \frac{(\gamma_ir_i)^{z_i+1/2}\theta_i^{2z_i}\exp(-\gamma_ir_i\theta_i^2)}{\Gamma(z_i+1/2)} \biggr]^{v}_{\theta_i = u} \Biggr\}. \label{F-diff-2-t2}
\end{align}

From \eqref{t1-condition-c}, we have
\begin{equation}
\int\biggl\lvert \frac{\partial}{\partial\theta_i}
 \Bigl\{ \frac{\partial f}{\partial\theta_i}(\theta) 
 \Bigr\} \biggr\rvert
\prod\limits_{j=1}^d\frac{(\gamma_jr_j)^{z_j+1/2}\theta_j^{2z_j} \exp(-\gamma_jr_j\theta_j^2)}{\Gamma(z_j+1/2)}\du\theta <\infty,\ \forall i. \label{2nd-diff-finite-t2}
\end{equation}

From \eqref{2nd-diff-finite-t2} and condition \eqref{t1-condition-d}, we have 
\begin{align}
    & \sum_{i=1}^d \mathop{\lim}_{a\to0\atop b\to\infty}\mathop{\lim}_{u\to0\atop v\to\infty} \Biggl[
-\int_u^v\int_{[a,b]^{d-1}}  \frac{\partial}{\partial\theta_i}\Bigl\{\frac{\partial f}{\partial\theta_i}(\theta)\Bigr\} 
\prod_{j\neq i}\frac{(\gamma_jr_j)^{z_j+1/2}\theta_j^{2z_j}\exp(-\gamma_jr_j\theta_j^2)}{\Gamma(z_j+1/2)}\du\theta_j\notag\\
&~~~~~~~~~~~~~~~~~~ 
\frac{(\gamma_ir_i)^{z_i+1/2}\theta_i^{2z_i}\exp(-\gamma_ir_i\theta_i^2)}{\Gamma(z_i+1/2)}\du\theta_i \Biggl] \notag\\
&=- \int  \sum_{i=1}^d \frac{\partial}{\partial\theta_i} \biggl\{ \frac{\partial f}{\partial\theta_i}(\theta)  \biggr\} 
\prod_{j=1}^d \frac{(\gamma_jr_j)^{z_j+1/2}\theta_j^{2z_j} \exp(-\gamma_jr_j\theta_j^2)}{\Gamma(z_j+1/2)}\du\theta \ge 0. \label{F-2nd-diff-2}
\end{align}

Using \eqref{F-diff-2-1-t2}, \eqref{F-diff-2-t2}, and \eqref{F-2nd-diff-2}, we obtain the inequality \eqref{thm2-step2}.

\vspace{0.2cm}

\noindent
\textit{\textbf{Step 3.}}
From Lemma~\ref{Lemma 2.6}, we know that \eqref{thm2-step2} $\ge 0$, which completes the proof.
\qed

\subsection{Lemmas used in the proofs of Theorem 1 and Theorem 2}

\noindent
\begin{lemma} \label{Lemma 1.1}
Under condition \eqref{Fineq}, for a given $r$, there exists $\epsilon>0$ such that
\begin{align}
&\sum_z\Biggl[ \max_{t\in[r-\epsilon,r+\epsilon]}\biggl\{ \Bigl| \log F(z,t) \Bigr| (\sum_{i=1}^dz_i+1)
\prod_{i=1}^d \frac{(t\lambda_i)^{z_i}\exp(-t\lambda_i)}{z_i!}
 \biggr\} \Biggr] <\infty. \label{lemma1-item}
\end{align}
In particular, $
 \text{E} \biggl( \Bigl| \log F(Z,r) \Bigr| \, \bigg| \, Z_i\sim\text{Po}(r\lambda_i),i=1,\dots,d \biggr)<\infty.
$
 \end{lemma}

\noindent
\textbf{Proof}\\
From \eqref{Fineq}, we know that 
\begin{equation}
\sum\limits_{i=1}^dz_i(F(z,t)-F(z-\delta_i,t))+\sum\limits_{i=1}^d(z_i+1/2)(F(z,t)-F(z+\delta_i,t))\ge 0,\ \forall z. \label{Fineq-L1}
\end{equation}
Using \eqref{Fineq-L1}, we have $F(z+\delta_i,t)<\frac{2}{1/2}\sum_{i=1}^d(z_i+1/2)F(z,t)=4\sum_{i=1}^d(z_i+1/2)F(z,t),\ \forall z$. Hence,
\begin{equation}
    F(z,t)\le\Bigl(4\sum_{i=1}^d(z_i+1/2)\Bigr)^{\sum_{i=1}^d z_i}F(\vec{0},t). \label{F-upper-L1}
\end{equation}
From \eqref{Fineq-L1}, we have $F(z-\delta_i,t)<2\sum_{i=1}^d(z_i+1/2)F(z,t)$ if $z_i>0$. Hence,
\begin{equation}
    F(z,t)\ge\Bigl(2\sum_{i=1}^d(z_i+1/2)\Bigr)^{-\sum_{i=1}^d z_i}F(\vec{0},t). \label{F-lower-L1}
\end{equation}
We choose $\epsilon=r/2$. Then, we have
\begin{align}
    \max_{t\in[r/2,3r/2]}F(\vec{0},t) 
    \le\int f(\bar\theta_1,\dots,\bar\theta_d)\prod\limits_{i=1}^d\frac{(3r/2)^{1/2}\bar\lambda_i^{-1/2}\exp(-r/2\bar\lambda_i)}{\Gamma(1/2)}\du\bar\lambda =3^{d/2}F(\vec{0},r/2) \label{F-upper-L1-2}.
\end{align}
We also have
\begin{align}
    \min_{t\in[r/2,3r/2]}F(\vec{0},t) 
    \ge\int f(\bar\theta_1,\dots,\bar\theta_d)\prod\limits_{i=1}^d\frac{(r/2)^{1/2}\bar\lambda_i^{-1/2}\exp(-3r/2\bar\lambda_i)}{\Gamma(1/2)}\du\bar\lambda =(1/3)^{d/2}F(\vec{0},3r/2) \label{F-lower-L1-2}.
\end{align}
From \eqref{F-upper-L1}, \eqref{F-lower-L1}, \eqref{F-upper-L1-2}, and \eqref{F-lower-L1-2}, there exists a constant $C$ such that $$\max_{t\in[r/2,3r/2]}\Bigl| \log F(z,t) \Bigr|<C(\sum_{i=1}^d z_i+C)^2,\ \forall z.$$
Thus, the left-hand side of \eqref{lemma1-item} is less than 
\begin{align*}
 & \sum_z\Biggl[ C(\sum_{i=1}^d z_i+C)^2(\sum_{i=1}^dz_i+1)
\prod_{i=1}^d \frac{(3r/2\lambda_i)^{z_i}\exp(-r/2\lambda_i)}{z_i!}
 \Biggr] \\
 &=\exp(r\sum_{i=1}^d\lambda_i)\text{E}\Bigl( C(\sum_{i=1}^d Z_i+C)^2(\sum_{i=1}^d Z_i+1)
\Big|\ Z_i\sim\text{Po}(3r/2\lambda_i) \Bigr) <\infty. \notag
\end{align*}
\qed

\vspace{0.3cm}

\noindent
\begin{lemma} \label{Lemma 1.2}
If the function $F(z,r)$ defined in Theorem~\ref{Theorem 1} is finite, then for any given $r$ and $z$, there exists $\epsilon>0$ such that
\begin{align}
\int \max_{t\in[r-\epsilon,r+\epsilon]}\Biggl\{ f(\bar\theta_1,\dots,\bar\theta_d) \prod_{i=1}^d
 \frac{t^{z_i+1/2}\bar\lambda_i^{z_i-1/2}\exp(-t\bar\lambda_i)}{\Gamma(z_i+1/2)}
 \Bigl| \sum_{i=1}^d \frac{z_i+1/2}{t} - \sum_{i=1}^d \bar\lambda_i \Bigr| \Biggr\} \du\bar\lambda<\infty. \label{lemma2-item}
\end{align}
 \end{lemma}
 
 \noindent
 \textbf{Proof}\\
 We choose $\epsilon=r/2$. Then, the left-hand side of \eqref{lemma2-item} is not greater than
 \begin{align*}
 &\int \max_{t\in[r/2,2r]} \Biggl\{ f(\bar\theta_1,\dots,\bar\theta_d) \prod_{i=1}^d
 \frac{t^{z_i+1/2}\bar\lambda_i^{z_i-1/2}\exp(-t\bar\lambda_i)}{\Gamma(z_i+1/2)}
 \Bigl| \sum_{i=1}^d \frac{z_i+1/2}{t} - \sum_{i=1}^d \bar\lambda_i \Bigr| \Biggr\} \du\bar\lambda\\
 &\le\int f(\bar\theta_1,\dots,\bar\theta_d) \prod_{i=1}^d
 \frac{(2r)^{z_i+1/2}\bar\lambda_i^{z_i-1/2}\exp(-\frac{r}{2}\bar\lambda_i)}{\Gamma(z_i+1/2)}
 \Bigl( \sum_{i=1}^d \frac{2z_i+1}{r} + \sum_{i=1}^d \bar\lambda_i \Bigr) \du\bar\lambda \\
 &=\frac{\sum_{i=1}^d(z_i+1/2)}{r}2^{2\sum_{i=1}^d(z_i+1/2)+1}F(z,r/2)+\frac{1}{r}2^{2\sum_{i=1}^d(z_i+1/2)+1}\sum_{j=1}^{d}F(z+\delta_j,r/2)(z_j+1/2).
 \end{align*}
 Because $F(z,r/2)$ and $F(z+\delta_j,r/2)$ are finite, the proof is complete.
\qed
 
\vspace{0.3cm}

\noindent
\begin{lemma} \label{Lemma 1.3}
Under condition \eqref{Fineq}, for a given $r$, there exists $\epsilon>0$ such that
\begin{align*}
&\sum_z \max_{t\in[r-\epsilon,r+\epsilon]} 
\Biggl[ \biggl| \log(F(z,t)) \biggr|
\Bigl\{ \prod_{i=1}^d \frac{(t\lambda_i)^{z_i}\exp(-t\lambda_i)}{z_i!} \Bigr\}
\Bigl| \sum_{i=1}^d \frac{z_i}{t} - \sum_{i=1}^d \lambda_i \Bigr| \Biggr] \notag \\
&+ \sum_z \max_{t\in[r-\epsilon,r+\epsilon]}
 \Biggl[
 \frac{  \Biggl|
 \int f(\bar\theta_1,\dots,\bar\theta_d) \prod_{i=1}^d
 \frac{t^{z_i+1/2}\bar\lambda_i^{z_i-1/2}\exp(-t\bar\lambda_i)}{\Gamma(z_i+1/2)}
 \Bigl( \sum_{i=1}^d \frac{z_i+1/2}{t} - \sum_{i=1}^d \bar\lambda_i \Bigr) \du\bar\lambda \Biggr| }
 { F(z,t)}
 \prod_{i=1}^d \frac{(t\lambda_i)^{z_i} \exp(-t\lambda_i)}{z_i!} 
 \Biggr] <\infty.
\end{align*}
\end{lemma}

\noindent
\textbf{Proof}\\
We know that
\begin{align*}
  \int f(\bar\theta_1,\dots,\bar\theta_d)
      \prod_{j=1}^d\frac{t^{z_j+1/2}\bar\lambda_j^{z_j-1/2}
  \exp(-t\bar\lambda_j)}{\Gamma(z_j+1/2)}
  \bar\lambda_i \du\bar\lambda \notag= F(z+\delta_i,t)\frac{z_i+1/2}{t}.
\end{align*}
Therefore, the lemma is equivalent to
\begin{align}
&\sum_z \max_{t\in[r-\epsilon,r+\epsilon]}
\Biggl[ \biggl| \log F(z,t) \biggr|
\Bigl\{
\prod_{i=1}^d \frac{(t\lambda_i)^{z_i}\exp(-t\lambda_i)}{z_i!}
\Bigr\}
\Bigl| \sum_{i=1}^d \frac{z_i}{t} - \sum_{i=1}^d \lambda_i \Bigr| \Biggr] \notag \\
&\quad + \sum_z \max_{t\in[r-\epsilon,r+\epsilon]}
 \Biggl[
 \frac{  \Biggl|
 F(z,t)\sum_{i=1}^d \frac{z_i+1/2}{t}-\sum_{i=1}^dF(z+\delta_i,t)\frac{z_i+1/2}{t} \Biggr| }
 {F(z,t)}
 \prod_{i=1}^d \frac{(t\lambda_i)^{z_i} \exp(-t\lambda_i)}{z_i!} 
 \Biggr] <\infty. \label{finite-goal}
\end{align}

From Lemma~\ref{Lemma 1.1}, we have $\epsilon\in(0,r)$ such that the first term on the left-hand side of \eqref{finite-goal} is not greater than
\begin{align*}
\sum_z \max_{t\in[r-\epsilon,r+\epsilon]}\Biggl[ \biggl| \log F(z,t) \biggr|
\biggl\{
\prod_{i=1}^d \frac{(t\lambda_i)^{z_i}\exp(-t\lambda_i)}{z_i!}
\biggr\}\frac{ \sum_{i=1}^d z_i +1}{\min\{ r-\epsilon,1/\sum_{i=1}^d \lambda_i \}} \Biggl] <\infty.
\end{align*}

Therefore, we need only prove that the second term on the left-hand side of \eqref{finite-goal} is finite. Using \eqref{Fineq}, we have
$\sum_{i=1}^d
(z_i+1/2)F(z+\delta_i,t)
 <2\sum_{i=1}^d(z_i+1/2)F(z,t).$
Therefore, the second term on the left-hand side of \eqref{finite-goal} is less than
\begin{align*}
& \sum_z \max_{t\in[r-\epsilon,r+\epsilon]}
 \Biggl[
 \frac{  
 F(z,t)\sum_{i=1}^d \frac{z_i+1/2}{t}}
 { F(z,t)}
 \prod_{i=1}^d \frac{(t\lambda_i)^{z_i} \exp(-t\lambda_i)}{z_i!} 
 \Biggr] \notag \\
 &\le \frac{\exp(2\epsilon\sum_{i=1}^d\lambda_i)}{r-\epsilon}\sum_z \biggl\{
 \sum_{i=1}^d (z_i+1/2)
 \prod_{i=1}^d \frac{((r+\epsilon)\lambda_i)^{z_i} \exp(-(r+\epsilon)\lambda_i)}{z_i!} \biggr\} <\infty. 
\end{align*}
\qed

\vspace{0.3cm}

\noindent
\begin{lemma} \label{Lemma 1.4}
Under conditions \eqref{t1-condition-a} and \eqref{t1-condition-b}, we have
\begin{align}
    &~~~~~F(z,r)<\infty,~~~~F(z+\delta_i,r)<\infty, \notag\\
    \int \biggl\{& \int f(\theta)\prod\limits_{j\neq i}
\frac{r^{z_j+1/2}\theta_j^{2z_j}\exp(-r\theta_j^2)}
{\Gamma(z_j+1/2)}\du\theta_j \biggr\} 
\frac{r^{z_i+1/2}\theta_i^{2z_i+1}
\exp(-r\theta_i^2)}{\Gamma(z_i+3/2)}\du\theta_i
<\infty, \label{lemma4-item1}\\
    \text{and}~~~~&\int  \Bigl\lvert\frac{\partial}{\partial\theta_i}f(\theta) \Bigl\rvert \prod\limits_{j=1}^d \frac{r^{z_j+1/2}
\theta_j^{2z_j+\delta_{ji}}\exp(-r\theta_j^2)}{\Gamma(z_j+\delta_{ji}+1/2)} \du\theta < \infty. \label{lemma4-item2}
\end{align}
\end{lemma}

\noindent
\textbf{Proof}\\
From \eqref{t1-condition-a} and the substitution $\lambda=\theta^2$, we have $F(z,r)<\infty$ and
\begin{equation}
\int f(\theta)\prod\limits_{j=1}^d
\frac{r^{z_j+1/2}\theta_j^{2z_j}\exp(-r\theta_j^2)}
{\Gamma(z_j+1/2+\delta_{ij})}\du\theta<\infty. \label{L3-infty-1}
\end{equation}
Substituting $z+\delta_i$ for $z$ yields $F(z+\delta_i,r)<\infty$ and
\begin{equation}
\int f(\theta)\prod\limits_{j=1}^d
\frac{r^{z_j+1/2}\theta_j^{2z_j+2\delta_{ij}}\exp(-r\theta_j^2)}
{\Gamma(z_j+1/2+\delta_{ij})}\du\theta<\infty. \label{L3-infty-2}
\end{equation}
Using \eqref{L3-infty-1}, \eqref{L3-infty-2}, and $\theta_i^{2z_i}+\theta_i^{2z_i+2}\ge\theta_i^{2z_i+1}$, we obtain \eqref{lemma4-item1}.
From \eqref{t1-condition-b}, we have
\begin{equation}
\int \Bigl\lvert\frac{\partial}{\partial\theta_i}f(\theta)\Bigl\rvert\prod\limits_{j=1}^d
\frac{r^{z_j+1/2}\theta_j^{2z_j}\exp(-r\theta_j^2)}
{\Gamma(z_j+1/2+\delta_{ij})}\du\theta<\infty. \label{L3-infty-3}
\end{equation}
From \eqref{t1-condition-b}, we have
\begin{equation}
\int \Bigl\lvert\frac{\partial}{\partial\theta_i}f(\theta)\Bigl\rvert\prod\limits_{j=1}^d
\frac{r^{z_j+1/2}\theta_j^{2z_j+2\delta_{ij}}\exp(-r\theta_j^2)}
{\Gamma(z_j+1/2+\delta_{ij})}\du\theta<\infty. \label{L3-infty-4}
\end{equation}
Using \eqref{L3-infty-3}, \eqref{L3-infty-4}, and $\theta_i^{2z_i}+\theta_i^{2z_i+2}\ge\theta_i^{2z_i+1}$, we obtain \eqref{lemma4-item2}.
\qed

\vspace{0.3cm}

\noindent
\begin{lemma} \label{Lemma 1.5}
If $f\in\mathbf{C}^2([0,\infty)^d)$ and conditions \eqref{t1-condition-a} and \eqref{t1-condition-b} are satisfied, then
\begin{align}
\biggl\{ \int_{[a,b]^{d-1}} f(\theta)\prod\limits_{j\neq i}
\frac{r^{z_j+1/2}\theta_j^{2z_j}\exp(-r\theta_j^2)}
{\Gamma(z_j+1/2)}\du\theta_j \biggr\}
\frac{r^{z_i+1/2}\theta_i^{2z_i+1}
\exp(-r\theta_i^2)}{\Gamma(z_i+3/2)} \label{lemma5-item}
\end{align} converges to $0$ as $\theta_i\to 0$ or $\theta_i\to\infty$.
\end{lemma}

\noindent
\textbf{Proof}
By performing integration by parts on $\theta_i$, we obtain
\begin{align}
&\int_u^v \biggl\{
\int_{[a,b]^{d-1}} f(\theta)\prod\limits_{j\neq i}\frac{r^{z_j+1/2}\theta_j^{2z_j}\exp(-r\theta_j^2)}
{\Gamma(z_j+1/2)} \du \theta_j 
\biggr\}
\frac{r^{z_i+1/2}\theta_i^{2z_i+1}\frac{\partial}{\partial\theta_i}\exp(-r\theta_i^2)}
{\Gamma(z_i+3/2)}\du\theta_i \notag\\
& ~~ +\int_u^v \biggl\{
\int_{[a,b]^{d-1}} f(\theta)\prod\limits_{j\neq i}\frac{r^{z_j+1/2}\theta_j^{2z_j}\exp(-r\theta_j^2)}
{\Gamma(z_j+1/2)} \du \theta_j 
\biggr\}
\frac{r^{z_i+1/2}(\frac{\partial}{\partial\theta_i}\theta_i^{2z_i+1})\exp(-r\theta_i^2)}
{\Gamma(z_i+3/2)}\du\theta_i \notag\\
& ~~ +\int_{[u,v]\times[a,b]^{d-1}} \biggl\{ \frac{\partial}{\partial\theta_i}f(\theta) \biggr\}
\prod\limits_{j=1}^d\frac{r^{z_j+1/2}
\theta_j^{2z_j+\delta_{ji}}\exp(-r\theta_j^2)}{\Gamma(z_j+\delta_{ji}+1/2)} \du\theta \notag\\
&=\Biggl[ \biggl\{
\int_{[a,b]^{d-1}} f(\theta)\prod\limits_{j\neq i}\frac{r^{z_j+1/2}\theta_j^{2z_j}\exp(-r\theta_j^2)}
{\Gamma(z_j+1/2)} \du \theta_j 
\biggr\}
\frac{r^{z_i+1/2}\theta_i^{2z_i+1}\exp(-r\theta_i^2)}
{\Gamma(z_i+3/2)} \Biggr]^{v}_{\theta_i = u}. \label{L4}
\end{align}
From Lemma~\ref{Lemma 1.4}, we know that all three terms on the left-hand side of \eqref{L4} converge as $v\to\infty$. Hence, the right-hand side of \eqref{L4} converges as $v\to\infty$. From Lemma~\ref{Lemma 1.4}, we have $$
\int \biggl\{ \int_{[a,b]^{d-1}} f(\theta)\prod\limits_{j\neq i}
\frac{r^{z_j+1/2}\theta_j^{2z_j}\exp(-r\theta_j^2)}
{\Gamma(z_j+1/2)}\du\theta_j \biggr\}
\frac{r^{z_i+1/2}\theta_i^{2z_i+1}
\exp(-r\theta_i^2)}{\Gamma(z_i+3/2)}\du\theta_i
<\infty.$$ Hence, \eqref{lemma5-item} converges to $0$ as $\theta_i\to\infty$.

Because $f\in \mathbf{C}^2([0,\infty)^d)$, $f$ is bounded on $[0,1]\times[a,b]^{d-1}$. Thus, $$\biggl\{ \int_{[a,b]^{d-1}} f(\theta)\prod\limits_{j\neq i}
\frac{r^{z_j+1/2}\theta_j^{2z_j}\exp(-r\theta_j^2)}
{\Gamma(z_j+1/2)}\du\theta_j \biggr\}$$ is bounded for $\theta_i\le1$.
Thus, \eqref{lemma5-item} converges to $0$ as $\theta_i\to 0$.
\qed

\vspace{0.3cm}

\noindent
\begin{lemma} \label{Lemma 1.6}
If $f\in\mathbf{C}^2([0,\infty)^d)$ and conditions \eqref{t1-condition-a}, \eqref{t1-condition-b}, \eqref{t1-condition-c}, and \eqref{t1-condition-e} are satisfied, then
\begin{align}
&\mathop{\lim}_{u\to0\atop v\to\infty} \Biggl[\int_{[a,b]^{d-1}} \frac{\partial f}{\partial\theta_i}(\theta)
\biggl\{\prod_{j\neq i}\frac{r^{z_j+1/2}\theta_j^{2z_j}\exp(-r\theta_j^2)}{\Gamma(z_j+1/2)}\du\theta_j\biggr\} \frac{r^{z_i+1/2}\theta_i^{2z_i}\exp(-r\theta_i^2)}{\Gamma(z_i+1/2)} \Biggr]^{v}_{\theta_i = u} \ge0. \notag
\end{align}
\end{lemma}

\noindent
\textbf{Proof}\\

From \eqref{t1-condition-b}, we have
$$\int \Bigl\lvert\frac{\partial f}{\partial\theta_i}(\theta)\Bigl\rvert\prod\limits_{j=1}^{d}\theta_j^{2z_j}\exp(-r\theta_j^2)\du\theta<\infty.$$
Hence,
\begin{equation}
\int \int_{[a,b]^{d-1}}  \frac{\partial f}{\partial\theta_i}(\theta)
\biggl\{\prod_{j\neq i}\frac{r^{z_j+1/2}\theta_j^{2z_j}\exp(-r\theta_j^2)}{\Gamma(z_j+1/2)}\du\theta_j\biggr\} \frac{r^{z_i+1/2}\theta_i^{2z_i}\exp(-r\theta_i^2)}{\Gamma(z_i+1/2)}\du\theta_i<\infty.
\label{2st-partial-3}
\end{equation}

By performing integration by parts on $\theta_i$, we obtain
\begin{align}
& \int_u^v\int_{[a,b]^{d-1}}  \frac{\partial f}{\partial\theta_i}(\theta)
\biggl\{\prod_{j\neq i}\frac{r^{z_j+1/2}\theta_j^{2z_j}\exp(-r\theta_j^2)}{\Gamma(z_j+1/2)} \du\theta_j\biggr\}
\frac{r^{z_i+1/2}\frac{\partial\theta_i^{2z_i}}{\partial\theta_i}\exp(-r\theta_i^2)}{\Gamma(z_i+1/2)}\du\theta_i \notag\\
& \qquad+ \int_u^v\int_{[a,b]^{d-1}} \frac{\partial f}{\partial\theta_i}(\theta)
\biggl\{\prod_{j\neq i}\frac{r^{z_j+1/2}\theta_j^{2z_j}\exp(-r\theta_j^2)}{\Gamma(z_j+1/2)}\du\theta_j\biggr\}
\frac{r^{z_i+1/2}\theta_i^{2z_i}\frac{\partial}{\partial\theta_i}\exp(-r\theta_i^2)}
{\Gamma(z_i+1/2)}\du\theta_i \notag\\
& \qquad+\int_u^v\int_{[a,b]^{d-1}}  \frac{\partial}{\partial\theta_i}\Bigl\{\frac{\partial f}{\partial\theta_i}(\theta)\Bigr\}
\biggl\{\prod_{j\neq i}\frac{r^{z_j+1/2}\theta_j^{2z_j}\exp(-r\theta_j^2)}{\Gamma(z_j+1/2)}\du\theta_j\biggr\} \frac{r^{z_i+1/2}\theta_i^{2z_i}\exp(-r\theta_i^2)}{\Gamma(z_i+1/2)}\du\theta_i \label{L5}\\
&= \biggl[\int_{[a,b]^{d-1}} \frac{\partial f}{\partial\theta_i}(\theta)
\biggl\{\prod_{j\neq i}\frac{r^{z_j+1/2}\theta_j^{2z_j}\exp(-r\theta_j^2)}{\Gamma(z_j+1/2)}\du\theta_j\biggr\} \frac{r^{z_i+1/2}\theta_i^{2z_i}\exp(-r\theta_i^2)}{\Gamma(z_i+1/2)} \biggr]^{v}_{\theta_i = u}. \notag
\end{align}

From \eqref{F-1st-diff} and $z_i \{ F(z,r)-F(z-\delta_i,r) \}<\infty$, we know that the first term in \eqref{L5} converges as $u\to0$ or $v\to\infty$.
From \eqref{F-1st-diff} and $(z_i+1/2) \{F(z,r)-F(z+\delta_i,r)\}<\infty$, we know that the second term in \eqref{L5} converges as $u\to0$ or $v\to\infty$.
From \eqref{2nd-diff-finite}, we know that the third term in \eqref{L5} converges as $u\to0$ or $v\to\infty$.
Thus, all three terms in \eqref{L5} converge as $u\to0$ or $v\to\infty$. Hence, $$\biggl[\int_{[a,b]^{d-1}}  \frac{\partial f}{\partial\theta_i}(\theta)
\biggl\{\prod_{j\neq i}\frac{r^{z_j+1/2}\theta_j^{2z_j}\exp(-r\theta_j^2)}{\Gamma(z_j+1/2)}\du\theta_j\biggr\} \frac{r^{z_i+1/2}\theta_i^{2z_i}\exp(-r\theta_i^2)}{\Gamma(z_i+1/2)} \biggr]^{v}_{\theta_i = u}$$
converges as $u\to0$ or $v\to\infty$. Accordingly, using \eqref{2st-partial-3}, we obtain
\begin{equation}
\lim_{\theta_i\to\infty} \int_{[a,b]^{d-1}}  \frac{\partial f}{\partial\theta_i}(\theta)
\biggl\{\prod_{j\neq i}\frac{r^{z_j+1/2}\theta_j^{2z_j}\exp(-r\theta_j^2)}{\Gamma(z_j+1/2)}\du\theta_j\biggr\} \frac{r^{z_i+1/2}\theta_i^{2z_i}\exp(-r\theta_i^2)}{\Gamma(z_i+1/2)}
= 0. \label{2st-partial-1}
\end{equation}

From \eqref{t1-condition-e}, we have
\begin{equation}
\lim_{\theta_i\to0}  \int_{[a,b]^{d-1}}  \frac{\partial f}{\partial\theta_i}(\theta)
\biggl\{\prod_{j\neq i}\frac{r^{z_j+1/2}\theta_j^{2z_j}\exp(-r\theta_j^2)}{\Gamma(z_j+1/2)}\du\theta_j\biggr\} \frac{r^{z_i+1/2}\theta_i^{2z_i}\exp(-r\theta_i^2)}{\Gamma(z_i+1/2)}\le0 \label{2st-partial-2}
\end{equation}

Using \eqref{2st-partial-1} and \eqref{2st-partial-2}, we complete the proof.
\qed

\vspace{0.3cm}

\noindent
\begin{lemma} \label{Lemma 2.1}
Under condition \eqref{Fineq3}, for a given $r$, there exists $\epsilon>0$ such that
\begin{align}
&\sum_z\Biggl[ \max_{t\in \bar B(r,\epsilon)}\biggl\{ \Bigl| \log F(z,t) \Bigr| (\sum_{i=1}^dz_i+1)
\prod_{i=1}^d \frac{(t_i\lambda_i)^{z_i}\exp(-t_i\lambda_i)}{z_i!}
 \biggr\} \Biggr] <\infty, \label{lemma7-term}
\end{align}
where $\bar B(r,\epsilon)$ is a closed ball with center $r=(r_1,\dots,r_d)$ and radius $\epsilon$.
In particular, $
 \text{E} \biggl( \Bigl| \log F(Z,r) \Bigr| \, \bigg| \, Z_i\sim\text{Po}(r_i\lambda_i),i=1,\dots,d \biggr)<\infty.
$ \end{lemma}

\noindent
\textbf{Proof}\\
From \eqref{Fineq3}, we know that 
\begin{equation}
\sum\limits_{i=1}^d\gamma_it_iz_i \Bigl\{ F(z,t)-F(z-\delta_i,t)\Bigr\} +\sum\limits_{i=1}^d\gamma_it_i(z_i+1/2) \Bigl\{F(z,t)-F(z+\delta_i,t) \Bigr\} \ge 0,\ \forall z. \label{Fineq-L2}
\end{equation}
Let $\gamma_0\coloneqq \frac{\min\{\gamma_1,\dots,\gamma_d\}}{\max\{\gamma_1,\dots,\gamma_d\}}$.\\ Using \eqref{Fineq-L2}, we obtain \[F(z+\delta_i,t)<\frac{2}{1/2 \cdot \gamma_0}\frac{\max\{t_1,\dots,t_d\}}{\min\{t_1,\dots,t_d\}}\sum_{i=1}^d(z_i+1/2)F(z,t)=\frac{4}{\gamma_0}\frac{\max\{t_1,\dots,t_d\}}{\min\{t_1,\dots,t_d\}}\sum_{i=1}^d(z_i+1/2)F(z,t),\ \forall z.\] Hence,
\begin{equation}
    F(z,t)\le\Bigl(\frac{4}{\gamma_0}\frac{\max\{t_1,\dots,t_d\}}{\min\{t_1,\dots,t_d\}}\sum_{i=1}^d(z_i+1/2)\Bigr)^{\sum_{i=1}^d z_i}F(\vec{0},t). \label{F-upper-L2}
\end{equation}
Using \eqref{Fineq-L2}, we obtain $F(z-\delta_i,t)<\frac{2}{\gamma_0}\frac{\max\{t_1,\dots,t_d\}}{\min\{t_1,\dots,t_d\}}\sum_{i=1}^d(z_i+1/2)F(z,t)$ if $z_i>0$. Hence,
\begin{equation}
    F(z,t)\ge\Bigl(\frac{2}{\gamma_0}\frac{\max\{t_1,\dots,t_d\}}{\min\{t_1,\dots,t_d\}}\sum_{i=1}^d(z_i+1/2)\Bigr)^{-\sum_{i=1}^d z_i}F(\vec{0},t). \label{F-lower-L2}
\end{equation}
We choose $\epsilon=\min\{r_1,\dots,r_d\}/2$. Then, we have
\begin{align}
    \max_{t\in\bar B(r,\epsilon)}F(\vec{0},t) \le\int f(\bar\theta_1,\dots,\bar\theta_d)\prod\limits_{i=1}^d\frac{(3r_i/2)^{1/2}\bar\lambda_i^{-1/2}\exp(-r_i/2\bar\lambda_i)}{\Gamma(1/2)}\du\bar\lambda =3^{d/2}F(\vec{0},r/2) \label{F-upper-L2-2}.
\end{align}
We also have
\begin{align}
    \min_{t\in\bar B(r,\epsilon)}F(\vec{0},t) \ge\int f(\bar\theta_1,\dots,\bar\theta_d)\prod\limits_{i=1}^d\frac{(r_i/2)^{1/2}\bar\lambda_i^{-1/2}\exp(-3r_i/2\bar\lambda_i)}{\Gamma(1/2)}\du\bar\lambda =(1/3)^{d/2}F(\vec{0},3r/2) \label{F-lower-L2-2}.
\end{align}
From \eqref{F-upper-L2}, \eqref{F-lower-L2}, \eqref{F-upper-L2-2}, and \eqref{F-lower-L2-2}, there exists a constant $C$ such that $$\max_{t\in\bar B(r,\epsilon)}\Bigl| \log F(z,t) \Bigr|<C(\sum_{i=1}^d z_i+C)^2,\ \forall z.$$
Thus, the left-hand side of \eqref{lemma7-term} is less than
\begin{align*}
 & \sum_z\Biggl[ C(\sum_{i=1}^d z_i+C)^2(\sum_{i=1}^dz_i+1)
\prod_{i=1}^d \frac{(3r_i/2\lambda_i)^{z_i}\exp(-r_i/2\lambda_i)}{z_i!}
 \Biggr] \\
 &=\exp(\sum_{i=1}^dr_i\lambda_i)\text{E}\Bigl( C(\sum_{i=1}^d Z_i+C)^2(\sum_{i=1}^d Z_i+1)
\Big|\ Z_i\sim\text{Po}(3r_i/2\lambda_i) \Bigr) <\infty. \notag
\end{align*}
\qed

\vspace{0.3cm}

\noindent
\begin{lemma} \label{Lemma 2.2}
If the function $F(z,r)$ defined in Theorem~\ref{Theorem 2} is finite, then for any given $r=(r_1,\dots,r_d)$, $s=(s_1,\dots,s_d)$, and $z=(z_1,\dots,z_d)$, we have
\begin{align}
\int \max_{\tau\in[0,1]}\Biggl\{ f(\bar\theta_1,\dots,\bar\theta_d) \prod_{i=1}^d
 \frac{t_i(\tau)^{z_i+1/2}\bar\lambda_i^{z_i-1/2}\exp(-t_i(\tau)\bar\lambda_i)}{\Gamma(z_i+1/2)}
 \Bigl| \sum_{i=1}^d (z_i+1/2)\gamma_it_i(\tau) - \sum_{i=1}^d \bar\lambda_i\gamma_it_i^2(\tau) \Bigr| \Biggr\} \du\bar\lambda<\infty. \label{lemma8-item}
\end{align}
 \end{lemma}
 
 \noindent
 \textbf{Proof}\\
 We define $r_0=\min\{r_1,\dots,r_d\}$ and $s_0=\max\{r_1+s_1,\dots,r_d+s_d\}$. Because $t_i(\tau)\in[r_i,r_i+s_i]$, we have $t_i(\tau)\in[r_0,s_0]$.
 Then, the left-hand side of \eqref{lemma8-item} is not greater than
 \begin{align*}
 &\int \Biggl[ f(\bar\theta_1,\dots,\bar\theta_d) \prod_{i=1}^d
 \frac{s_0^{z_i+1/2}\bar\lambda_i^{z_i-1/2}\exp(-r_0\bar\lambda_i)}{\Gamma(z_i+1/2)}
 \Bigl\{ \sum_{i=1}^d (z_i+1/2)\gamma_is_0 + \sum_{i=1}^d \bar\lambda_i\gamma_is_0^2 \Bigr\} \Biggl] \du\bar\lambda \\
 &=\sum_{i=1}^d (z_i+1/2)\gamma_is_0\Bigl( \frac{s_0}{r_0}\Bigr) ^{\sum_{i=1}^d(z_i+1/2)}F(z,(r_0,\dots,r_0))
 +\sum_{j=1}^{d}\gamma_js_0\Bigl(\frac{s_0}{r_0}\Bigr)^{\sum_{i=1}^d(z_i+1/2)+1}F(z+\delta_j,(r_0,\dots,r_0))(z_j+1/2).
 \end{align*}
 
 Because $F(z,(r_0,\dots,r_0))$ and $F(z+\delta_j,(r_0,\dots,r_0))$ are finite, the proof is complete.
\qed
 
 \vspace{0.3cm}

\noindent
\begin{lemma} \label{Lemma 2.3}
Under condition \eqref{Fineq3}, for a given $\tau_0$, there exists $\epsilon>0$ such that 
\begin{align*}
&\sum_z \max_{|\tau-\tau_0|\le\epsilon} 
\Biggl[ \Biggl| \biggl\{ \log F(z,t(\tau)) \biggr\}
\biggl\{
\prod_{i=1}^d \frac{(t_i(\tau)\lambda_i)^{z_i}\exp(-t_i(\tau)\lambda_i)}{z_i!}
\biggr\}
\Bigl(\sum_{i=1}^d z_i\gamma_it_i(\tau) - \sum_{i=1}^d \lambda_i\gamma_it_i^2(\tau) \Bigr) \Biggr| \Biggr] \notag \\
&\quad + \sum_z \max_{|\tau-\tau_0|\le\epsilon}
 \Biggl[
 \frac{  \biggl|
 \int f(\bar\theta_1,\dots,\bar\theta_d) \prod_{i=1}^d
 \frac{t_i(\tau)^{z_i+1/2}\bar\lambda_i^{z_i-1/2}\exp(-t_i(\tau)\bar\lambda_i)}{\Gamma(z_i+1/2)}
 \Bigl( \sum_{i=1}^d (z_i+1/2)\gamma_it_i(\tau) - \sum_{i=1}^d \bar\lambda_i\gamma_it_i^2(\tau) \Bigr) \du\bar\lambda \biggr| } 
 { F(z,t(\tau))}
  \notag\\
 & ~~~~~
 \biggl\{
\prod_{i=1}^d \frac{(t_i(\tau)\lambda_i)^{z_i}\exp(-t_i(\tau)\lambda_i)}{z_i!}
\biggr\}
 \Biggr] \notag <\infty.
\end{align*}
\end{lemma}

\noindent
\textbf{Proof}\\
We know that
\begin{align*}
  &\int f(\bar\theta_1,\dots,\bar\theta_d)
\prod_{j=1}^d\frac{t_j(\tau)^{z_j+1/2}\bar\lambda_j^{z_j-1/2}
  \exp(-t_j(\tau)\bar\lambda_j)}{\Gamma(z_j+1/2)}
  \bar\lambda_i\gamma_it_i^2(\tau) \du\bar\lambda = F(z+\delta_i,t(\tau))(z_i+1/2)\gamma_it_i(\tau).
\end{align*}
Therefore, the lemma is equivalent to that
\begin{align}
&\sum_z \max_{|\tau-\tau_0|\le\epsilon}
\Biggl[ \biggl| \log F(z,t(\tau)) \biggr|
\Bigl\{
\prod_{i=1}^d \frac{(t_i(\tau)\lambda_i)^{z_i}\exp(-t_i(\tau)\lambda_i)}{z_i!}
\Bigr\}
\biggl| \sum_{i=1}^d z_i\gamma_it_i(\tau) - \sum_{i=1}^d \lambda_i\gamma_it_i^2(\tau) \biggr| \Biggl] \notag \\
& + \sum_z \max_{|\tau-\tau_0|\le\epsilon}
 \Biggl\{
 \frac{ \biggl|
 F(z,t(\tau))\sum\limits_{i=1}^d (z_i+1/2)\gamma_it_i(\tau)-\sum\limits_{i=1}^dF(z+\delta_i,t(\tau))(z_i+1/2)\gamma_it_i(\tau) \biggr| }
 {F(z,t(\tau))} \prod_{i=1}^d \frac{(t_i(\tau)\lambda_i)^{z_i} \exp(-t_i(\tau)\lambda_i)}{z_i!} 
 \Biggr\} \label{finite-goal-t2}
\end{align}
 is finite.

From Lemma~\ref{Lemma 2.1}, when we set $r=t(\tau_0)$, we know that there exists $\delta>0$ such that
\begin{align}
&\sum_z\Biggl[ \max_{t\in \bar B(t(\tau_0),\delta)}\biggl\{ \Bigl| \log F(z,t) \Bigr| (\sum_{i=1}^dz_i+1)
\prod_{i=1}^d \frac{(t_i\lambda_i)^{z_i}\exp(-t_i\lambda_i)}{z_i!}
 \biggr\} \Biggr] <\infty. \notag
\end{align}
Because $t(\tau)$ is continuous, there exists $\epsilon>0$ such that for any $\tau\in[\tau_0-\epsilon,\tau_0+\epsilon]$ , $t(\tau)\in \bar B(t(\tau_0),\delta)$.
Thus, the first term on the left-hand side of \eqref{finite-goal-t2} is not greater than
\begin{align}
& \sum_z \max_{t\in \bar B(t(\tau_0),\delta)}
\Biggl[ \biggl| \log F(z,t) \biggr|
\Bigl\{
\prod_{i=1}^d \frac{(t_i\lambda_i)^{z_i}\exp(-t_i\lambda_i)}{z_i!}
\Bigr\}
\biggl| \sum_{i=1}^d z_i\gamma_it_i - \sum_{i=1}^d \lambda_i\gamma_it_i^2 \biggr| \Biggl] \notag \\
&\le\sum_z \max_{t\in \bar B(t(\tau_0),\delta)}\Biggl[ \biggl| \log F(z,t) \biggr| \Bigl\{\prod_{i=1}^d \frac{(t_i\lambda_i)^{z_i}\exp(-t_i\lambda_i)}{z_i!}\Bigr\}( \sum_{i=1}^d z_i +1) \Biggl]\Bigl\{ \sum_{i=1}^d \gamma_i(\Vert t(\tau_0)\Vert+\delta)+\sum_{i=1}^d \lambda_i\gamma_i(\Vert t(\tau_0)\Vert+\delta)^2 \Bigr\} \notag\\
&<\infty. \label{L2-finite1-t2}
\end{align}

Therefore, we need only prove that the second term on the left-hand side of \eqref{finite-goal-t2} is finite. Using \eqref{Fineq3}, when we set $r=t(\tau)$, we obtain
$\sum_{i=1}^d\gamma_it_i(\tau)(z_i+1/2)F(z+\delta_i,t(\tau))
 <2\sum_{i=1}^d\gamma_it_i(\tau)(z_i+1/2)F(z,t(\tau)).$
 We define $r_0=\min\{r_1,\dots,r_d\}$ and $s_0=\max\{r_1+s_1,\dots,r_d+s_d\}$. Because $t_i(\tau)\in[r_i,r_i+s_i]$, we have $t_i(\tau)\in[r_0,s_0]$. Therefore, the second term on the left-hand side of \eqref{finite-goal-t2} is less than
 \begin{align*}
& \sum_z \max_{|\tau-\tau_0|\le\epsilon}
 \Biggl\{
 \frac{ 
 F(z,t(\tau))\sum_{i=1}^d (z_i+1/2)\gamma_it_i(\tau) }
 {F(z,t(\tau))}\prod_{i=1}^d \frac{(t_i(\tau)\lambda_i)^{z_i} \exp(-t_i(\tau)\lambda_i)}{z_i!} 
 \Biggr\} \notag\\
 &\le \exp((s_0-r_0)\sum_{i=1}^d\lambda_i)\sum_z 
 \Biggl\{
 \biggl( \sum_{i=1}^d (z_i+1/2)\gamma_is_0\biggl) \prod_{i=1}^d \frac{(s_0\lambda_i)^{z_i} \exp(-s_0\lambda_i)}{z_i!} 
 \Biggr\} <\infty.
\end{align*}
\qed

\vspace{0.3cm}

\noindent
\begin{lemma} \label{Lemma 2.4}
Under conditions \eqref{t1-condition-a} and \eqref{t1-condition-b}, we have
\begin{align}
    &~~~~~F(z,r)<\infty,~~~~F(z+\delta_i,r)<\infty, \notag\\
    \int \biggl\{& \int f(\theta)\prod\limits_{j\neq i}
\frac{(\gamma_jr_j)^{z_j+1/2}\theta_j^{2z_j}\exp(-\gamma_jr_j\theta_j^2)}
{\Gamma(z_j+1/2)}\du\theta_j \biggr\} 
\frac{(\gamma_ir_i)^{z_i+1/2}\theta_i^{2z_i+1}
\exp(-\gamma_ir_i\theta_i^2)}{\Gamma(z_i+3/2)}\du\theta_i
<\infty, \label{lemma10-item1}\\
    \text{and}~~~~&\int  \Bigl\lvert\frac{\partial}{\partial\theta_i}f(\theta) \Bigl\rvert \prod\limits_{j=1}^d \frac{(\gamma_jr_j)^{z_j+1/2}
\theta_j^{2z_j+\delta_{ji}}\exp(-\gamma_jr_j\theta_j^2)}{\Gamma(z_j+\delta_{ji}+1/2)} \du\theta < \infty. \label{lemma10-item2}
\end{align}
\end{lemma}

\noindent
\textbf{Proof}\\
From \eqref{t1-condition-a} and the substitution $\lambda=\gamma\theta^2$, we have $F(z,r)<\infty$ and
\begin{equation}
\int f(\theta)\prod\limits_{j=1}^d
\frac{(\gamma_jr_j)^{z_j+1/2}\theta_j^{2z_j}\exp(-\gamma_jr_j\theta_j^2)}
{\Gamma(z_j+1/2+\delta_{ij})}\du\theta<\infty. \label{L3-infty-1-t2}
\end{equation}
Substituting $z+\delta_i$ for $z$ yields $F(z+\delta_i,r)<\infty$ and
\begin{equation}
\int f(\theta)\prod\limits_{j=1}^d
\frac{(\gamma_jr_j)^{z_j+1/2}\theta_j^{2z_j+2\delta_{ij}}\exp(-\gamma_jr_j\theta_j^2)}
{\Gamma(z_j+1/2+\delta_{ij})}\du\theta<\infty. \label{L3-infty-2-t2}
\end{equation}
Using \eqref{L3-infty-1-t2}, \eqref{L3-infty-2-t2}, and $\theta_i^{2z_i}+\theta_i^{2z_i+2}\ge\theta_i^{2z_i+1}$, we obtain \eqref{lemma10-item1}.
From \eqref{t1-condition-b}, we have
\begin{equation}
\int \Bigl\lvert\frac{\partial}{\partial\theta_i}f(\theta)\Bigl\rvert\prod\limits_{j=1}^d
\frac{(\gamma_jr_j)^{z_j+1/2}\theta_j^{2z_j}\exp(-\gamma_jr_j\theta_j^2)}
{\Gamma(z_j+1/2+\delta_{ij})}\du\theta<\infty. \label{L3-infty-3-t2}
\end{equation}
From \eqref{t1-condition-b}, we have
\begin{equation}
\int \Bigl\lvert\frac{\partial}{\partial\theta_i}f(\theta)\Bigl\rvert\prod\limits_{j=1}^d
\frac{(\gamma_jr_j)^{z_j+1/2}\theta_j^{2z_j+2\delta_{ij}}\exp(-\gamma_jr_j\theta_j^2)}
{\Gamma(z_j+1/2+\delta_{ij})}\du\theta<\infty. \label{L3-infty-4-t2}
\end{equation}
Using \eqref{L3-infty-3-t2}, \eqref{L3-infty-4-t2}, and $\theta_i^{2z_i}+\theta_i^{2z_i+2}\ge\theta_i^{2z_i+1}$, we obtain \eqref{lemma10-item2}.
\qed

\vspace{0.3cm}

\noindent
\begin{lemma} \label{Lemma 2.5}
If $f\in\mathbf{C}^2([0,\infty)^d)$ and conditions \eqref{t1-condition-a} and \eqref{t1-condition-b} are satisfied, then
\begin{align}
\biggl\{ \int_{[a,b]^{d-1}} f(\theta)\prod\limits_{j\neq i}
\frac{(\gamma_jr_j)^{z_j+1/2}\theta_j^{2z_j}\exp(-\gamma_jr_j\theta_j^2)}
{\Gamma(z_j+1/2)}\du\theta_j \biggr\}
\frac{(\gamma_ir_i)^{z_i+1/2}\theta_i^{2z_i+1}
\exp(-\gamma_ir_i\theta_i^2)}{\Gamma(z_i+3/2)} \label{lemma11-item}
\end{align} converges to $0$ as $\theta_i\to 0$ or $\theta_i\to\infty$.
\end{lemma}

\noindent
\textbf{Proof}
By performing integration by parts on $\theta_i$, we obtain
\begin{align}
&\int_u^v \biggl\{
\int_{[a,b]^{d-1}} f(\theta)\prod\limits_{j\neq i}\frac{(\gamma_jr_j)^{z_j+1/2}\theta_j^{2z_j}\exp(-\gamma_jr_j\theta_j^2)}
{\Gamma(z_j+1/2)} \du \theta_j 
\biggr\}
\frac{(\gamma_ir_i)^{z_i+1/2}\theta_i^{2z_i+1}\frac{\partial}{\partial\theta_i}\exp(-\gamma_ir_i\theta_i^2)}
{\Gamma(z_i+3/2)}\du\theta_i \notag\\
& ~~ +\int_u^v \biggl\{
\int_{[a,b]^{d-1}} f(\theta)\prod\limits_{j\neq i}\frac{(\gamma_jr_j)^{z_j+1/2}\theta_j^{2z_j}\exp(-\gamma_jr_j\theta_j^2)}
{\Gamma(z_j+1/2)} \du \theta_j 
\biggr\}
\frac{(\gamma_ir_i)^{z_i+1/2}(\frac{\partial}{\partial\theta_i}\theta_i^{2z_i+1})\exp(-\gamma_ir_i\theta_i^2)}
{\Gamma(z_i+3/2)}\du\theta_i \notag\\
& ~~ +\int_{[u,v]\times[a,b]^{d-1}} \biggl\{ \frac{\partial}{\partial\theta_i}f(\theta) \biggr\}
\prod\limits_{j=1}^d\frac{(\gamma_jr_j)^{z_j+1/2}
\theta_j^{2z_j+\delta_{ji}}\exp(-\gamma_jr_j\theta_j^2)}{\Gamma(z_j+\delta_{ji}+1/2)} \du\theta \notag\\
&=\Biggl[ \biggl\{
\int_{[a,b]^{d-1}} f(\theta)\prod\limits_{j\neq i}\frac{(\gamma_jr_j)^{z_j+1/2}\theta_j^{2z_j}\exp(-\gamma_jr_j\theta_j^2)}
{\Gamma(z_j+1/2)} \du \theta_j 
\biggr\}
\frac{(\gamma_ir_i)^{z_i+1/2}\theta_i^{2z_i+1}\exp(-\gamma_ir_i\theta_i^2)}
{\Gamma(z_i+3/2)} \Biggr]^{v}_{\theta_i = u}. \label{L4-t2}
\end{align}
From Lemma~\ref{Lemma 2.4}, we know that all three terms on the left-hand side of \eqref{L4-t2} converge as $v\to\infty$. Hence, the right-hand side of \eqref{L4-t2} converges as $v\to\infty$. From Lemma~\ref{Lemma 2.4}, we have $$
\int \biggl\{ \int_{[a,b]^{d-1}} f(\theta)\prod\limits_{j\neq i}
\frac{(\gamma_jr_j)^{z_j+1/2}\theta_j^{2z_j}\exp(-\gamma_jr_j\theta_j^2)}
{\Gamma(z_j+1/2)}\du\theta_j \biggr\}
\frac{(\gamma_ir_i)^{z_i+1/2}\theta_i^{2z_i+1}
\exp(-\gamma_ir_i\theta_i^2)}{\Gamma(z_i+3/2)}\du\theta_i
<\infty.$$ Hence, \eqref{lemma11-item} converges to $0$ as $\theta_i\to\infty$.

Because $f\in \mathbf{C}^2([0,\infty)^d)$, $f$ is bounded on $[0,1]\times[a,b]^{d-1}$. Thus, $$ \int_{[a,b]^{d-1}} f(\theta)\prod\limits_{j\neq i}
\frac{(\gamma_jr_j)^{z_j+1/2}\theta_j^{2z_j}\exp(-\gamma_jr_j\theta_j^2)}
{\Gamma(z_j+1/2)}\du\theta_j $$ is bounded for $\theta_i\le1$.
Thus, \eqref{lemma11-item} converges to $0$ as $\theta_i\to 0$.
\qed

\vspace{0.3cm}

\noindent
\begin{lemma} \label{Lemma 2.6}
If $f\in\mathbf{C}^2([0,\infty)^d)$ and conditions \eqref{t1-condition-a}, \eqref{t1-condition-b}, \eqref{t1-condition-c}, and \eqref{t1-condition-e} are satisfied, then
\begin{align*}
&\mathop{\lim}_{u\to0\atop v\to\infty} \Biggl[\int_{[a,b]^{d-1}} \frac{\partial f}{\partial\theta_i}(\theta)
\biggl\{ \prod_{j\neq i}\frac{(\gamma_jr_j)^{z_j+1/2}\theta_j^{2z_j}\exp(-\gamma_jr_j\theta_j^2)}{\Gamma(z_j+1/2)}\du\theta_j\biggr\} \frac{(\gamma_ir_i)^{z_i+1/2}\theta_i^{2z_i}\exp(-\gamma_ir_i\theta_i^2)}{\Gamma(z_i+1/2)} \Biggr]^{v}_{\theta_i = u} \ge0. 
\end{align*}
\end{lemma}

\noindent
\textbf{Proof}\\
From \eqref{t1-condition-b}, we have
$$\int \Bigl\lvert\frac{\partial f}{\partial\theta_i}(\theta)\Bigl\rvert\prod\limits_{j=1}^{d}\theta_j^{2z_j}\exp(-\gamma_jr_j\theta_j^2)\du\theta<\infty.$$
Thus,
\begin{equation}
\int \int_{[a,b]^{d-1}}  \frac{\partial f}{\partial\theta_i}(\theta)
\biggl\{ \prod_{j\neq i}\frac{(\gamma_jr_j)^{z_j+1/2}\theta_j^{2z_j}\exp(-\gamma_jr_j\theta_j^2)}{\Gamma(z_j+1/2)}\du\theta_j\biggr\} \frac{(\gamma_ir_i)^{z_i+1/2}\theta_i^{2z_i}\exp(-\gamma_ir_i\theta_i^2)}{\Gamma(z_i+1/2)}\du\theta_i<\infty.
\label{2st-partial-3-t2}
\end{equation}

By performing integration by parts on $\theta_i$, we obtain
\begin{align}
& \int_u^v\int_{[a,b]^{d-1}}  \frac{\partial f}{\partial\theta_i}(\theta)
\biggl\{\prod_{j\neq i}\frac{(\gamma_jr_j)^{z_j+1/2}\theta_j^{2z_j}\exp(-\gamma_jr_j\theta_j^2)}{\Gamma(z_j+1/2)} \du\theta_j\biggr\}
\frac{(\gamma_ir_i)^{z_i+1/2}\frac{\partial\theta_i^{2z_i}}{\partial\theta_i}\exp(-\gamma_ir_i\theta_i^2)}{\Gamma(z_i+1/2)}\du\theta_i \notag\\
& \qquad+ \int_u^v\int_{[a,b]^{d-1}} \frac{\partial f}{\partial\theta_i}(\theta)
\biggl\{\prod_{j\neq i}\frac{(\gamma_jr_j)^{z_j+1/2}\theta_j^{2z_j}\exp(-\gamma_jr_j\theta_j^2)}{\Gamma(z_j+1/2)}\du\theta_j\biggr\}
\frac{(\gamma_ir_i)^{z_i+1/2}\theta_i^{2z_i}\frac{\partial}{\partial\theta_i}\exp(-\gamma_ir_i\theta_i^2)}
{\Gamma(z_i+1/2)}\du\theta_i \notag\\
& \qquad+\int_u^v\int_{[a,b]^{d-1}}  \frac{\partial}{\partial\theta_i}\Bigl\{\frac{\partial f}{\partial\theta_i}(\theta)\Bigr\}
\biggl\{ \prod_{j\neq i}\frac{(\gamma_jr_j)^{z_j+1/2}\theta_j^{2z_j}\exp(-\gamma_jr_j\theta_j^2)}{\Gamma(z_j+1/2)}\du\theta_j\biggr\} \frac{(\gamma_ir_i)^{z_i+1/2}\theta_i^{2z_i}\exp(-\gamma_ir_i\theta_i^2)}{\Gamma(z_i+1/2)}\du\theta_i \label{L5-t2}\\
&= \biggl[\int_{[a,b]^{d-1}}  \frac{\partial f}{\partial\theta_i}(\theta)
\biggl\{ \prod_{j\neq i}\frac{(\gamma_jr_j)^{z_j+1/2}\theta_j^{2z_j}\exp(-\gamma_jr_j\theta_j^2)}{\Gamma(z_j+1/2)}\du\theta_j\biggr\} \frac{(\gamma_ir_i)^{z_i+1/2}\theta_i^{2z_i}\exp(-\gamma_ir_i\theta_i^2)}{\Gamma(z_i+1/2)} \biggr]^{v}_{\theta_i = u}. \notag
\end{align}

From \eqref{F-1st-diff-t2} and $\gamma_ir_iz_i \{ F(z,r)-F(z-\delta_i,r) \}<\infty$, we know that the first term in \eqref{L5-t2} converges as $u\to0$ or $v\to\infty$.
From \eqref{F-1st-diff-t2} and $\gamma_ir_i(z_i+1/2) \{F(z,r)-F(z+\delta_i,r)\}<\infty$, we know that the second term in \eqref{L5-t2} converges as $u\to0$ or $v\to\infty$.
From \eqref{2nd-diff-finite-t2}, we know that the third term in \eqref{L5-t2} converges as $u\to0$ or $v\to\infty$.
Thus, all three terms in \eqref{L5-t2} converge as $u\to0$ or $v\to\infty$. Hence, $$\biggl[\int_{[a,b]^{d-1}}  \frac{\partial f}{\partial\theta_i}(\theta)
\biggl\{ \prod_{j\neq i}\frac{(\gamma_jr_j)^{z_j+1/2}\theta_j^{2z_j}\exp(-\gamma_jr_j\theta_j^2)}{\Gamma(z_j+1/2)}\du\theta_j\biggr\} \frac{(\gamma_ir_i)^{z_i+1/2}\theta_i^{2z_i}\exp(-\gamma_ir_i\theta_i^2)}{\Gamma(z_i+1/2)} \biggr]^{v}_{\theta_i = u}$$
converges as $u\to0$ or $v\to\infty$. Thus, using \eqref{2st-partial-3-t2}, we obtain
\begin{equation}
\lim_{\theta_i\to\infty} \int_{[a,b]^{d-1}}  \frac{\partial f}{\partial\theta_i}(\theta)
\biggl\{ \prod_{j\neq i}\frac{(\gamma_jr_j)^{z_j+1/2}\theta_j^{2z_j}\exp(-\gamma_jr_j\theta_j^2)}{\Gamma(z_j+1/2)}\du\theta_j\biggr\} \frac{(\gamma_ir_i)^{z_i+1/2}\theta_i^{2z_i}\exp(-\gamma_ir_i\theta_i^2)}{\Gamma(z_i+1/2)}
= 0. \label{2st-partial-1-t2}
\end{equation}

From \eqref{t1-condition-e}, we have
\begin{equation}
\lim_{\theta_i\to0} \int_{[a,b]^{d-1}}  \frac{\partial f}{\partial\theta_i}(\theta)
\biggl\{ \prod_{j\neq i}\frac{(\gamma_jr_j)^{z_j+1/2}\theta_j^{2z_j}\exp(-\gamma_jr_j\theta_j^2)}{\Gamma(z_j+1/2)}\du\theta_j\biggr\}\frac{(\gamma_ir_i)^{z_i+1/2}\theta_i^{2z_i}\exp(-\gamma_ir_i\theta_i^2)}{\Gamma(z_i+1/2)}\le0 \label{2st-partial-2-t2}
\end{equation}

Using \eqref{2st-partial-1-t2} and \eqref{2st-partial-2-t2}, we complete the proof.
\qed

\subsection{Proof of Proposition 1}

\begin{proposition}\label{Proposition 2} The Bayesian predictive distribution based on $\pi_{f,\mathrm{J}}(\lambda)$ dominates that based on the Jeffreys prior and is thus nearly minimax for the prediction of independent Poisson processes with the same duration or different durations when $$f(\theta)=\sum_{a\in\{1,-1\}^d}\Big(\sum\limits_{i=1}^d(a_i\theta_i-\eta_i)^2\Big)^{-\alpha},$$ where $0<\alpha\le (d-2)/2$.
\end{proposition}

\noindent
\textbf{Proof.} Independent Poisson processes with the same duration correspond to $r_1=r_2=\cdots=r_d$ and $s_1=s_2=\cdots=s_d$. Therefore, we need only prove that $f(\theta)=\sum_{a\in\{1,-1\}^d}\big(\sum_{i=1}^d(a_i\theta_i-\eta_i)^2\big)^{-\alpha}$ satisfies the conditions in the first half of Theorem~\ref{Theorem 2}. We show this in two steps. In Step 1, we prove that the condition \eqref{Fineq3} is satisfied. In Step 2, we prove that $F(z,r)$ is not a constant function of $z$. The details of each step are presented below.

\noindent
\textit{\textbf{Step 1.}} We first show that $f_{\epsilon}(\theta)=\sum_{a\in\{1,-1\}^d}\big(\sum_{i=1}^d(a_i\theta_i-\eta_i)^2+\epsilon\big)^{-\alpha}$ satisfies the conditions in the second half of Theorem~\ref{Theorem 2} when $\epsilon>0,\,0<\alpha\le d/2-1$. Here, $\theta_i=\sqrt{\frac{\lambda_i}{\gamma_i}},\ i=1,\dots,d$.

\eqref{t1-condition-a} is satisfied because
$
    f_{\epsilon}(\theta)\le  2^d\epsilon^{-\alpha} = \exp(\mathrm{o}(\sum_{j=1}^d\theta_j^2)).
$

\eqref{t1-condition-b} is satisfied because
\begin{align}
    \Bigl\lvert\frac{\partial f_{\epsilon}}{\partial\theta_i}(\theta)
\Bigr\rvert
    \le 2^{d+1}\alpha\epsilon^{-\alpha-1}(\theta_i+\vert\eta_i\vert)
    = \exp(\mathrm{o}(\sum_{j=1}^d\theta_j^2)). \notag
\end{align}

\eqref{t1-condition-c} is satisfied because
\begin{align}
    \Bigl\lvert\frac{\partial^2 f_{\epsilon}}{\partial\theta_i^2}(\theta)\Bigl\rvert
    \le 2^{d+2}\alpha(\alpha+1)\epsilon^{-\alpha-2}(\theta_i+\vert\eta_i\vert)^2 + 2^{d+1}\alpha\epsilon^{-\alpha-1}
    = \exp(\mathrm{o}(\sum_{j=1}^d\theta_j^2)). \notag
\end{align}

$h(\theta)=\big(\sum_{i=1}^d(\theta_i-\eta_i)^2+\epsilon\big)^{-\alpha}$ is superharmonic because
\begin{align}
    \sum\limits_{i=1}^{d} \frac{\partial^2 h}{\partial\theta_i^2}(\theta) 
    = 2\alpha\Big(\sum\limits_{i=1}^d(\theta_i-\eta_i)^2+\epsilon \Big)^{-\alpha-2}\Big\{(2(\alpha+1)-d)\sum\limits_{i=1}^d(\theta_i-\eta_i)^2-d\epsilon \Big\} \le 0. \notag
\end{align}

\par From the statements above, we know that $f_{\epsilon}(\theta)$ satisfies the conditions of the second half of Theorem~\ref{Theorem 2}. Therefore, \eqref{Fineq3} holds when $f(\theta)=\sum_{a\in\{1,-1\}^d}\big(\sum_{i=1}^d(a_i\theta_i-\eta_i)^2+\epsilon\big)^{-\alpha},\,\epsilon>0,\, 0<\alpha\le d/2-1$.

\par Next, we prove that \eqref{Fineq3} still holds when $f(\theta)=\sum_{a\in\{1,-1\}^d}\big(\sum_{i=1}^d(a_i\theta_i-\eta_i)^2\big)^{-\alpha}$. The key is to consider $\epsilon\to0$.
Let $\zeta\coloneqq \min\{\gamma_1r_1,\dots,\gamma_dr_d\}$. By using inequality $u^{2v}\exp(-2\sqrt{v}u)\le v^v\exp(-2v)$, we know that there exists constant $c$ such that
\begin{equation*}
    \frac{u^{2v}}{\Gamma(v+1/2)}\le c\exp(2\sqrt{v}u-v),\ \forall u,v\ge0. 
\end{equation*}
Therefore, there exists constant $C$ such that for every $\theta$ and $z$,
\begin{align}
\prod\limits_{j=1}^d\frac{r_j^{z_j+1/2}\lambda_j^{z_j-1/2}\exp(-r_j\lambda_j)}{\Gamma(z_j+1/2)}\du\lambda &= \prod\limits_{j=1}^d\frac{2(r_j\gamma_j)^{1/2}(\sqrt{r_j\gamma_j}\theta_j)^{2z_j}\exp(-r_j\gamma_j\theta_j^2)}{\Gamma(z_j+1/2)}\du\theta 
\notag\\&\le C\prod\limits_{j=1}^d\frac{\exp(-(\sqrt{r_j\gamma_j}\theta_j-\sqrt{z_j})^2)}{\sqrt{2\pi\zeta^{-1}}}\du\theta. \label{p2-ineq}
\end{align}
Therefore, we have
\begin{align}
    \int(\sum\limits_{j=1}^d(a_j\theta_j-\eta_j)^2)^{-\alpha}\prod\limits_{j=1}^d\frac{r_j^{z_j+1/2}\lambda_j^{z_j-1/2}\exp(-r_j\lambda_j)}{\Gamma(z_j+1/2)}\du\lambda &\le \int(\sum\limits_{j=1}^d(a_j\theta_j-\eta_j)^2)^{-\alpha}C\prod\limits_{j=1}^d\frac{\exp(-\zeta(\theta_j-\sqrt{\frac{z_j}{r_j\gamma_j}})^2/2)}{\sqrt{2\pi\zeta^{-1}}}\du\theta \notag\\
    & =\text{E}\Bigl(C\Vert\theta-a\eta\Vert^{-2\alpha}\mid \theta\sim\mathrm{N}_d(\sqrt{\frac{z}{r\gamma}},\zeta^{-1}I_{d})\Bigr) \notag\\&<\infty, \label{p2-ineq2}
\end{align} 
where $\sqrt{\frac{z}{r\gamma}}$ denotes the vector whose i-th element is $\sqrt{\frac{z_i}{r_i\gamma_i}}$ and $a\eta$ denotes the vector whose i-th element is $a_i\eta_i.$
Therefore, using the dominated convergence theorem, for every $a$, we have $$\int(\sum\limits_{j=1}^d(a_j\theta_j-\eta_j)^2+\epsilon)^{-\alpha}\prod\limits_{j=1}^d\frac{r_j^{z_j+1/2}\lambda_j^{z_j-1/2}\exp(-r_j\lambda_j)}{\Gamma(z_j+1/2)}\du\lambda \to \int(\sum\limits_{j=1}^d(a_j\theta_j-\eta_j)^2)^{-\alpha}\prod\limits_{j=1}^d\frac{r_j^{z_j+1/2}\lambda_j^{z_j-1/2}\exp(-r_j\lambda_j)}{\Gamma(z_j+1/2)}\du\lambda$$
when $\epsilon\to0$. From the definition of $F$ in \eqref{Fineq3}, we know that \eqref{Fineq3} still holds.

\vspace{0.2cm}

\noindent
\textit{\textbf{Step 2.}} 
$F(z,r)$ is not a constant function of $z$ because when $z_1\to\infty$, from the inequality \eqref{p2-ineq2},
\begin{align*}
    \int(\sum\limits_{j=1}^d(a_j\theta_j-\eta_j)^2)^{-\alpha}\prod\limits_{j=1}^d\frac{r_j^{z_j+1/2}\lambda_j^{z_j-1/2}\exp(-r_j\lambda_j)}{\Gamma(z_j+1/2)}\du\lambda &\le \text{E}\Bigl(C\Vert\theta-a\eta\Vert^{-2\alpha}\mid \theta\sim\mathrm{N}_d(\sqrt{\frac{z}{r\gamma}},\zeta^{-1}I_{d})\Bigr)\to 0.
\end{align*} 

Therefore, $f(\theta)=\sum_{a\in\{1,-1\}^d}\big(\sum_{i=1}^d(a_i\theta_i-\eta_i)^2\big)^{-\alpha}$ satisfies the conditions in the first half of Theorem~\ref{Theorem 2}. Thus, the Bayesian predictive distribution $p_{f,\mathrm{J}}(y\mid x)$ dominates the Bayesian predictive distribution $p_{\mathrm{J}}(y\mid x)$ for the prediction of independent Poisson processes with the same or different durations.
\qed

\subsection{Proof of Proposition 2}
\begin{proposition}\label{proposition 3} The Bayesian predictive distribution based on $\pi_{f,\mathrm{J}}(\lambda)$ dominates that based on the Jeffreys prior and is thus nearly minimax for the prediction of independent Poisson processes with the same duration or different durations when $$f(\theta)=\sum_{a\in\{1,-1\}^d}(s_V(a_1\theta_1,a_2\theta_2,\dots,a_d\theta_d))^{-\alpha}\eqqcolon\sum_{a\in\{1,-1\}^d}(s_V(a\theta))^{-\alpha},$$ where $0<\alpha\le (d-k-2)/2$.
\end{proposition}

\noindent
\textbf{Proof.} We need only prove that $f(\theta)=\sum_{a\in\{1,-1\}^d}(s_V(a\theta))^{-\alpha}$ satisfies the conditions in the first half of Theorem~\ref{Theorem 2}. We show this in two steps. In Step 1, we prove that the condition \eqref{Fineq3} is satisfied. In Step 2, we prove that $F(z,r)$ is not a constant function of $z$.

\noindent
\textit{\textbf{Step 1.}} 
We construct a $\mathbf{C}^2$ prior $$f_{\epsilon}(\theta)=\sum_{a\in\{1,-1\}^d}(s_V(a\theta)+\epsilon)^{-\alpha},\ 0<\epsilon<1.$$ We first show that $f_{\epsilon}(\theta)$ satisfies the conditions of the second half of Theorem~\ref{Theorem 2}.

\par First, we have $
f_{\epsilon}(\theta)\le 2^d\epsilon^{-\alpha} 
=\exp(\mathrm{o}(\sum_{j=1}^d\theta_j^2)), $
\begin{align}
\Bigl\lvert \frac{\partial f_{\epsilon}}{\partial\theta_i}(\theta) \Bigr\rvert
\le \epsilon^{-\alpha-1}\alpha 2^{d+1}(d-k)(\sum\limits_{j=1}^{d}\theta_j)
= \exp(\mathrm{o}(\sum_{j=1}^d\theta_j^2)), \notag
\end{align}
and
\begin{align}
\Bigl\lvert \frac{\partial^2 f_{\epsilon}}{\partial\theta_i^2}(\theta) \Bigr\rvert
\le \epsilon^{-\alpha-2}\alpha(\alpha+1)2^{d+2}(d-k)^2(\sum\limits_{j=1}^{d}\theta_j)^2+\epsilon^{-\alpha-1}\alpha2^{d+1}(d-k)
= \exp(\mathrm{o}(\sum_{j=1}^d\theta_j^2)). \notag
\end{align}

\par Second, we know that $h(\theta)=(s_V(\theta)+\epsilon)^{-\alpha}$ is a superharmonic function because
\begin{align}
\sum\limits_{i=1}^{d} \frac{\partial^2 h}{\partial\theta_i^2}(\theta) 
&= \alpha(\sum\limits_{j=1}^{d-k}\langle \theta,v_j \rangle^2+\epsilon)^{-\alpha-2}\biggl\{ (4(\alpha+1)-2(d-k))\sum\limits_{j=1}^{d-k}\langle \theta,v_j \rangle^2-2(d-k)\epsilon\biggr\} <0. \notag
\end{align}

\par Therefore, from the second half of Theorem~\ref{Theorem 2}, we know that \eqref{Fineq3} holds when $f(\theta)=\sum_{a\in\{1,-1\}^d}(s_V(a\theta)+\epsilon)^{-\alpha},\ 0<\alpha\le(d-k-2)/2$. Next, we prove that \eqref{Fineq3} still holds when $f(\theta)=\sum_{a\in\{1,-1\}^d}(s_V(a\theta))^{-\alpha}$. The key is to consider $\epsilon\to0$.

From \eqref{p2-ineq}, we have
\begin{align}
    \int(s_V(a\theta))^{-\alpha}\prod\limits_{i=1}^d\frac{r_i^{z_i+1/2}\lambda_i^{z_i-1/2}\exp(-r_i\lambda_i)}{\Gamma(z_i+1/2)}\du\lambda 
    &\le \int(\sum\limits_{i=1}^{d-k}\langle \theta,av_i \rangle ^2)^{-\alpha}C\prod\limits_{j=1}^d\frac{\exp(-\zeta(\theta_j-\sqrt{\frac{z_j}{r_j\gamma_j}})^2/2)}{\sqrt{2\pi\zeta^{-1}}}\du\theta \notag\\
    &=C\text{E}\Big[(\sum\limits_{i=1}^{d-k}\langle \theta,av_i \rangle ^2)^{-\alpha}\mid \theta\sim\mathrm{N}_d(\sqrt{\frac{z}{r\gamma}},\zeta^{-1}I_{d})\Big] \notag\\
    &=C\text{E}\Big[\Vert\mu\Vert^{-2\alpha}\mid \mu\sim\mathrm{N}_{d-k}((av_1,\dots,av_{d-k})^\top\sqrt{\frac{z}{r\gamma}},\zeta^{-1}I_{d-k})\Big]<\infty. \label{p3-Fvalue}
\end{align}

Therefore, using the dominated convergence theorem, for every $a$,
$$\int (s_V(a\theta)+\epsilon)^{-\alpha}\prod\limits_{i=1}^d\frac{r_i^{z_i+1/2}\lambda_i^{z_i-1/2}\exp(-r_i\lambda_i)}{\Gamma(z_i+1/2)}\du\lambda \to \int (s_V(a\theta))^{-\alpha}\prod\limits_{i=1}^d\frac{r_i^{z_i+1/2}\lambda_i^{z_i-1/2}\exp(-r_i\lambda_i)}{\Gamma(z_i+1/2)}\du\lambda$$
when $\epsilon\to0$. From the definition of $F$ in \eqref{Fineq3}, \eqref{Fineq3} still holds when $f(\theta)=\sum_{a\in\{1,-1\}^d}(s_V(a\theta))^{-\alpha}$.

\vspace{0.2cm}

\noindent
\textit{\textbf{Step 2.}} 
From \eqref{p3-Fvalue}, we note that $F(z,r)\to0$ when $\min_{a\in\{1,-1\}^d}\vert\langle \sqrt{\frac{z}{r\gamma}},av_1 \rangle\vert \to\infty$. Because there exists $z\in\mathbb{N}^d$ such that $\min_{a\in\{1,-1\}^d}\vert\langle \sqrt{\frac{z}{r\gamma}},av_1 \rangle\vert>0$, there exists $z\in\mathbb{N}^d$ such that $\min_{a\in\{1,-1\}^d}\vert\langle \sqrt{\frac{z}{r\gamma}},av_1 \rangle\vert>c$ for any constant $c$.

From the first half of Theorem~\ref{Theorem 2}, we know that the Bayesian predictive distribution based on $\pi_{f,\mathrm{J}}(\lambda)$ dominates that based on the Jeffreys prior.
\qed

\subsection{Proof of Proposition 3}
\begin{proposition}\label{proposition 4} The Bayesian predictive distribution based on $\pi_{f,\mathrm{J}}(\lambda)$ dominates that based on the Jeffreys prior and is thus nearly minimax for the prediction of independent Poisson processes with the same duration or different durations when $$f(\theta)=\sum_{a\in\{1,-1\}^d}(s_{V_1}(a\theta))^{-\alpha_1}+\sum_{a\in\{1,-1\}^d}(s_{V_2}(a\theta))^{-\alpha_2}+\cdots+\sum_{a\in\{1,-1\}^d}(s_{V_n}(a\theta))^{-\alpha_n},$$ where $0<\alpha_i\le (d-k_i-2)/2$ and $k_i$ denotes the dimension of $V_i$, $i=1,2,\dots,n$.
\end{proposition}

\noindent
\textbf{Proof.} From Step 1 of the proof of Proposition~\ref{proposition 3} and additivity of $F$, we know that the combination
$$f(\theta)=\sum_{a\in\{1,-1\}^d}(s_{V_1}(a\theta))^{-\alpha_1}+\sum_{a\in\{1,-1\}^d}(s_{V_2}(a\theta))^{-\alpha_2}+\cdots+\sum_{a\in\{1,-1\}^d}(s_{V_n}(a\theta))^{-\alpha_n}$$
satisfies the inequality condition \eqref{Fineq3}. 
From Step 2 of the proof of Proposition~\ref{proposition 3}, $F(z,r)\to0$ when $$\min_{a\in\{1,-1\}^d,\ 1\le i\le n}\vert\langle \sqrt{\frac{z}{r\gamma}},av_{i1} \rangle\vert \to\infty,$$ where $v_{i1}$ is the first vector of a standard orthonormal basis of the complementary space of $V_i$. Thus, $F(z,r)$ is not a constant function of $z$. In summary, $f(\theta)$ satisfies the conditions in the first half of Theorem~\ref{Theorem 2}.
Therefore, the Bayesian predictive distribution $p_{f,\mathrm{J}}(y\mid x)$ dominates that based on the Jeffreys prior.\qed

\subsection{Proof of Theorem 3}

Theorem 3 in the main paper is copied as follows.
\begin{theorem}\label{Theorem 3}\leavevmode
\begin{enumerate}[label=\arabic*)]
\item
For any $\lambda$, the K-L risk of $p_{\mathrm{J}}(y\mid x)$ is less than $0.52\sum_{i=1}^d\log\big((r_i+s_i)/r_i\big)$.
\item
For any predictive distribution $q(y\mid x)$ and positive number $\epsilon$, there exists $\lambda$ such that the K-L risk of $q(y\mid x)$ is greater than $0.5\sum_{i=1}^d\log\big((r_i+s_i)/r_i\big)-\epsilon$.
\end{enumerate}
\end{theorem}

\noindent
\textbf{Proof.} We use the following lemma for the proof of the theorem. The proof of the lemma is presented at Li \cite{li2024nearly}. Let $\text{Po}(\lambda)$ denote the Poisson distribution with parameter $\lambda$.
\noindent
\begin{lemma}\label{Lemma 13}
For any $\lambda>0$, $f(\lambda)=\lambda\text{E}\Big(\log\big(\frac{x+0.5}{\lambda}\big)\mid x\sim\text{Po}(\lambda)\Big)>-0.02$. Furthermore, $\lim_{\lambda\to\infty}f(\lambda)=0$.
\end{lemma}

\vspace{0.3cm}

\noindent
\textit{\textbf{Proof of part 1.}}
Similar to Eq.(A.1) in the main paper for the case of same duration, we have
$$p_{\mathrm{J}}(y\mid x)=\prod_{i=1}^d\Big(\frac{r_i}{r_i+s_i}\Big)^{x_i+1/2}\Big(\frac{s_i}{r_i+s_i}\Big)^{ y_i}\frac{\Gamma(x_i+y_i+1/2)}{\Gamma(x_i+1/2)y_i!}.$$ Therefore, the K-L risk $\text{E}\big(D(p(y \mid \lambda),p_{\mathrm{J}}(y\mid x))\big)$ is given by
\begin{align}
    &\text{E}\Big(\log p(y\mid \lambda)-\log p_{\mathrm{J}}(y\mid x) \mid \lambda\Big)\notag\\
    &=\sum_{i=1}^d\biggl(-s_i\lambda_i+s_i\lambda_i\log(s_i\lambda_i)-(r_i\lambda_i+\frac{1}{2})\log\Big(\frac{r_i}{r_i+s_i}\Big)-s_i\lambda_i\log\Big(\frac{s_i}{r_i+s_i}\Big)\notag\\
    &\qquad-\text{E}\Bigl(\log\Gamma(x_i+y_i+\frac{1}{2})-\log\Gamma(x_i+\frac{1}{2})\;\big|\;\lambda_i \Bigr) \biggr) \label{t3-1}
\end{align}
Considering the functions
$$F_i(t):=t\lambda_i\log\lambda_i+\frac{1}{2}\log t+\lambda_i(t\log t-t)-\text{E}\Bigl(\log\Gamma(x+\frac{1}{2})\;\big|\;x\sim\text{Po}(t\lambda_i)\Bigr),$$
the K-L risk \eqref{t3-1} is equal to
\begin{equation}
 \sum_{i=1}^d\Big(F_i(r_i+s_i)-F_i(r_i)\Big)=\sum_{i=1}^d\Big(\int_{r_i}^{r_i+s_i} F_i^{\prime}(t)\du t\Big), \label{t3-1.1}   
\end{equation}
where
$$F_i^{\prime}(t)=\frac{1}{2t}-\text{E}\Bigl(\lambda_i\log\Big(\frac{x+0.5}{t\lambda_i}\Big)\;\big|\;x\sim \text{Po}(t\lambda_i)\Bigr).$$
From Lemma \ref{Lemma 13}, $F_i^{\prime}(t)<0.52/t.$ Thus, \eqref{t3-1.1} is equal to $\sum_{i=1}^d\Big(\int_{r_i}^{r_i+s_i} F_i^{\prime}(t)\du t\Big)<0.52\sum_{i=1}^d\log\big((r_i+s_i)/r_i\big)$.

\vspace{0.3cm}

\noindent
\textit{\textbf{Proof of part 2.}}
We only need to show that $0.5\sum_{i=1}^d\log\big((r_i+s_i)/r_i\big)$ is the Bayes risk limit of a sequence of Bayes rules $p_{\pi_n}$ with $$\pi_n(\lambda)=\prod_{i=1}^d\lambda_i^{-1/2}\exp(-\frac{\lambda_i}{n})\frac{1}{n^{1/2}\Gamma(1/2)}.$$

Similar to Eq.(A.3) in the main paper for the case of same duration, we have
$$p_{\pi_n}(y\mid x)=\prod_{i=1}^d\Big(\frac{r_i+1/n}{r_i+s_i+1/n}\Big)^{x_i+1/2}\Big(\frac{s_i}{r_i+s_i+1/n}\Big)^{y_i}\frac{\Gamma(x_i+y_i+1/2)}{\Gamma(x_i+1/2)y_i!}.$$
The Bayes risk of $p_{\pi_n}$ is 
\begin{align}
&\text{E}\Big(\text{E}\Big(\sum\limits_yp(y\mid\lambda)
\log \frac{p(y\mid\lambda)}{p_{\pi_n}(y\mid x)}\Big)\;\Big|\; \lambda_i\sim\Gamma(\frac{1}{2},\frac{1}{n})\Big) \notag\\
&=\text{E}\Big(\text{E}\Big(\sum\limits_yp(y\mid\lambda)
\log \frac{p(y\mid\lambda)}{p_{\mathrm{J}}(y\mid x)}\Big) \;\Big|\;\lambda_i\sim \Gamma(\frac{1}{2},\frac{1}{n})\Big)+\text{E}\Big(\text{E}\Big(\sum\limits_yp(y\mid\lambda)
\log \frac{p_{\mathrm{J}}(y\mid x)}{p_{\pi_n}(y\mid x)}\Big)\;\Big|\;\lambda_i\sim \Gamma(\frac{1}{2},\frac{1}{n})\Big). \label{t3-2}
\end{align}
We use \eqref{t3-1.1} in the proof of part 1 to show that the left term in \eqref{t3-2} is equal to
\begin{align}
   &\text{E}\bigg(\sum_{i=1}^d \int_{r_i}^{r_i+s_i} \biggl(\frac{1}{2t}-\text{E}\Bigl(\lambda_i\log(\frac{x_i+0.5}{t\lambda_i})\;\Big|\;x_i\sim \text{Po}(t\lambda_i)\Bigr) \biggr) \du t\;\Big|\;\lambda_i\sim \Gamma(\frac{1}{2},\frac{1}{n})\bigg) \notag\\
   &= 0.5\sum_{i=1}^d\log\Big(\frac{r_i+s_i}{r_i}\Big)-\sum_{i=1}^d\int_{r_i}^{r_i+s_i}\frac{1}{t}\text{E}\bigg(\text{E}\Bigl(t\lambda_i\log(\frac{x_i+0.5}{t\lambda_i})\;\Big|\;x_i\sim \text{Po}(t\lambda_i)\Bigr)\;\Big|\;\lambda_i\sim \Gamma(\frac{1}{2},\frac{1}{n})\bigg)\du t \notag\\
   &= 0.5\sum_{i=1}^d\log\Big(\frac{r_i+s_i}{r_i}\Big)-\sum_{i=1}^d\int_{r_i}^{r_i+s_i}\frac{1}{t}\text{E}\bigg(f(t\lambda_i)\;\Big|\;\lambda_i\sim \Gamma(\frac{1}{2},\frac{1}{n})\bigg)\du t \label{t3-3}.
\end{align}
According to $\lim_{\lambda\to\infty}f(\lambda)=0$ from Lemma \ref{Lemma 13}, \eqref{t3-3} converges to $0.5\sum_{i=1}^d\log\big((r_i+s_i)/r_i\big)$ when $n\to\infty$.
Based on $$\frac{p_{\mathrm{J}}(y\mid x)}{p_{\pi_n}(y\mid x)}=\prod_{i=1}^d\Big(\frac{r_i}{r_i+1/n}\Big)^{x_i+1/2}\Big(\frac{r_i+s_i+1/n}{r_i+s_i}\Big)^{x_i+y_i+1/2},$$ when $n\to\infty$, the right term in \eqref{t3-2} is equal to
\begin{align}
&\sum_{i=1}^d\text{E}\bigg(\text{E}\Big((x_i+1/2)\log\Big(\frac{r_i}{r_i+1/n}\Big)+(x_i+y_i+1/2)\log\Big(\frac{r_i+s_i+1/n}{r_i+s_i}\Big)\Big)\;\Big|\;\lambda_i\sim \Gamma(\frac{1}{2},\frac{1}{n})\bigg) \notag\\
&=\sum_{i=1}^d\text{E}\Big((r_i\lambda_i+1/2)\log\Big(\frac{r_i}{r_i+1/n}\Big)+((r_i+s_i)\lambda_i+1/2)\log\Big(\frac{r_i+s_i+1/n}{r_i+s_i}\Big)\;\Big|\;\lambda_i\sim \Gamma(\frac{1}{2},\frac{1}{n})\Big)  \notag\\
   &\to \sum_{i=1}^dn/2\Big(r_i\log\Big(\frac{r_i}{r_i+1/n}\Big)+(r_i+s_i)\log\Big(\frac{r_i+s_i+1/n}{r_i+s_i}\Big)\Big) \to 0. \notag
\end{align}
Therefore, \eqref{t3-2} converges to $0.5\sum_{i=1}^d\log\big((r_i+s_i)/r_i\big)$ when $n\to\infty$, which completes the proof.
\qed

\end{appendices}

\bibliographystyle{IEEEtran}
\bibliography{scholar}